\documentclass{article}
\usepackage{graphicx} 
\usepackage{float}
\usepackage[left=1in,top=1in,right=1in,bottom=1in,nohead]{geometry}
\usepackage{multirow}
\usepackage{url,subfig,latexsym,amsmath,amssymb, amsfonts,enumerate,
  epsfig, graphics,times, amsthm,mathtools, bbm}
\usepackage{appendix}
\usepackage{soul}

\usepackage{tikz, pgfplots}
\usetikzlibrary{automata,positioning,arrows, intersections, calc}
\pgfplotsset{compat=1.8}

\usepackage{mathdots}
  
\usepackage{commath}
\usepackage{wrapfig}
\usepackage{xcolor}
\usepackage{nameref}
\usepackage{listings}
\usepackage[dvipsnames]{xcolor} 
\usepackage{hyperref}
\hypersetup{
    colorlinks=true,
    allcolors=MidnightBlue 
}
\usepackage{makeidx}
\usepackage{tabularx}
\usepackage{array}
    \newcolumntype{P}[1]{>{\centering\arraybackslash}p{#1}}
    \newcolumntype{M}[1]{>{\centering\arraybackslash}m{#1}}
\usepackage{comment}
\usetikzlibrary{arrows.meta,positioning}

\usetikzlibrary {graphs}
\usepackage{chemfig}

\newtheorem{thm}{Theorem}[section]
\newtheorem{cor}[thm]{Corollary}
\newtheorem{lem}[thm]{Lemma}
\newtheorem{prop}[thm]{Proposition}

\newtheorem*{thm*}{Theorem}

\theoremstyle{definition}
\newtheorem{defn}{Definition}[section]
\newtheorem{rem}{Remark}[section]

\newtheorem{cond}{Condition}[section]

\newcommand{\var}{\mathrm{Var}}

\def\F{\mathcal F}

\def\T1{T^1_{\{x_n\}}}

\def\R{\mathbb R}

\def\Z{\mathbb Z}

\def\T1{T^1_{\{x_n\}}}
\def\Td1{T^{D,1}_{\{x_n\}}}
\def\Ts1{T^{S,1}_{\{x_n\}}}

\def\E{\mathbb{E}}

\def\cov{\mathrm{Cov}}

\def\Tid1{T^{D,1}_{\{x_n+\zeta_i\}}}
\def\Tis1{T^{S,1}_{\{x_n+\zeta_i\}}}
\def\Tjs1{T^{S,1}_{\{x_n+\zeta_j\}}}
\def\Tjd1{T^{D,1}_{\{x_n+\zeta_j\}}}

\newcount\colveccount
\newcommand*\colvec[1]{
        \global\colveccount#1
        \begin{pmatrix}
        \colvecnext
}
\def\colvecnext#1{
        #1
        \global\advance\colveccount-1
        \ifnum\colveccount>0
                \\
                \expandafter\colvecnext
        \else
                \end{pmatrix}
        \fi
}

\title{
Mean-level arithmetic computation via stochastic reaction networks}

\author{
David F. Anderson\thanks{Department of Mathematics, University of
Wisconsin-Madison, USA.  anderson@math.wisc.edu.}
\and 
Seokhwan Moon\thanks{Department of Mathematics, University of
Wisconsin-Madison, USA.  moon64@wisc.edu. Corresponding author.}
}

\date{\today}

\begin{document}

\maketitle

\begin{abstract}
    Reaction networks provide a natural model for analog computation in which inputs and outputs are encoded by molecular abundances. Much of the work on analog computation with chemical reaction networks (CRNs) has focused on deterministic mass-action systems. In DNA- and protein-based circuits, however, low copy numbers can make intrinsic fluctuations significant, motivating the use of stochastic models. In this paper, we develop a framework for computation with stochastic CRNs at the level of stationary means. Specifically, we construct elementary arithmetic modules for identification, addition, and multiplication, and analyze their mean-level correctness and convergence rates. We focus on these three modules because they are simple enough to analyze explicitly while already exhibiting the new issues raised by stochastic mean-level computation. Unlike in deterministic systems, the computed value cannot be read from a single output measurement. Instead, readout requires Monte Carlo estimation, with error arising from transient bias as the system approaches stationarity and intrinsic sampling noise that persists even at stationarity. We therefore study the convergence rates of output means, variances, and autocovariances to their stationary values, as well as the stationary fluctuations governing sampling error. We also show that this readout cost can depend strongly, and sometimes unexpectedly, on the architecture of the computation. In particular, comparing linear-chain and binary-tree implementations reveals that binary trees can greatly reduce stationary variance for addition, but can instead produce exponentially larger variance for multiplication.
\end{abstract}

\section{Introduction}
\label{sec:introduction}

Reaction networks have long been studied as models of computation~\cite{adleman1994molecular, magnasco1997chemical}. One important approach uses deterministic dynamics, in which the reaction network is modeled by a system of ordinary differential equations (ODEs). In this deterministic framework, input and output values are encoded by the concentrations of designated species, often in the long time limit. This ODE-based perspective has led to constructions for arithmetic operations~\cite{anderson2025chemical, buisman2009computing}, polynomial and more general function evaluation~\cite{anderson2026computing, solanki2023computing}, and neural-network-like computations~\cite{anderson2021reaction, kang2025noise, cherry2018scaling}. Rate-independent computation has also been studied in continuous CRNs, where correctness is required across a broad class of kinetic laws~\cite{chen2023rate}.

At the same time, discrete stochastic reaction networks have also been studied as models of computation. This viewpoint is especially natural for biochemical systems operating in small volumes and at low molecular copy numbers, where intrinsic randomness cannot be neglected. In this setting, the reaction network is typically modeled as a discrete-valued continuous-time Markov chain whose state gives the copy numbers of the various species, rather than a vector of real-valued concentrations. In one major line of work that uses these discrete models,  inputs are encoded by molecular counts, and the output is read from a stable or eventual molecular configuration~\cite{soloveichik2008computation, cummings2016probability, chen2014deterministic}. A notable advantage of that viewpoint is that the output of the computation is largely independent of the particular kinetic law: the emphasis is on reachability and eventual stabilization of molecular counts, rather than on convergence of a specific system of rate equations.  It has been shown that this ``count-based'' framework can perform powerful forms of computation, including Turing-universal computation~\cite{soloveichik2008computation, cummings2016probability}.

There are also possible limitations of this count-based framework.  First, in many such constructions, the input is specified by an exact molecular count, so implementing the computation requires initializing the network with specified molecular counts. There is interesting work on how efficiently such counts can be generated by CRNs~\cite{luchsinger2023optimal}, but this may be challenging in some implementations.  Second, in many such constructions the input species are consumed or transformed as part of the computation, so the input value is not retained in a reusable form.
In contrast, the modules studied here use input species catalytically: the input species determine output intensities but are not consumed by the computational module, and are therefore reusable.

A third mode of computation comes from stochastic logic, where values in $[0,1]$ are represented by probabilities, and realized through empirical frequencies of random bit streams, often modeled as sequences of independent Bernoulli random variables \cite{gaines1967stochastic, qian2010architecture, alaghi2013survey}.  In that setting, elementary logic gates implement arithmetic operations at the level of probabilities; for example, an AND gate applied to independent bit streams with probabilities $p$ and $q$ produces a bit stream with probability $pq$.  This perspective has also been connected to chemical reaction networks through deterministic mass-action CRN implementations of stochastic-logic constructions, in which the relevant probabilities are represented by concentration fractions rather than by intrinsic stochastic dynamics~\cite{solanki2023computing}.

In this paper, we combine features of the deterministic and discrete stochastic viewpoints above.  As in ODE-based analog computation, we encode real-valued quantities through mean values of designated species.  
However, the underlying model is a discrete stochastic reaction network, and the relevant mean is taken with respect to the stationary distribution of the stochastic process. 
More specifically, we say that a chemical species, $X$, encodes a positive real number, $x$, at the mean level if its marginal under the stationary distribution satisfies
\begin{align*}
\E_\pi[X]=x .
\end{align*}
This encoding has a simple and useful feature: real-valued quantities can be represented with very small reaction networks.  For example, the birth-death system
\begin{align*}
    \emptyset \xrightarrow{x} X,
    \qquad
    X \xrightarrow{1} \emptyset,
\end{align*}
has stationary distribution $\operatorname{Poisson}(x)$
for the copy number of $X$, and hence stationary mean $x$.

\begin{figure}[h!]
\centering
    \includegraphics[width=1.0\linewidth]{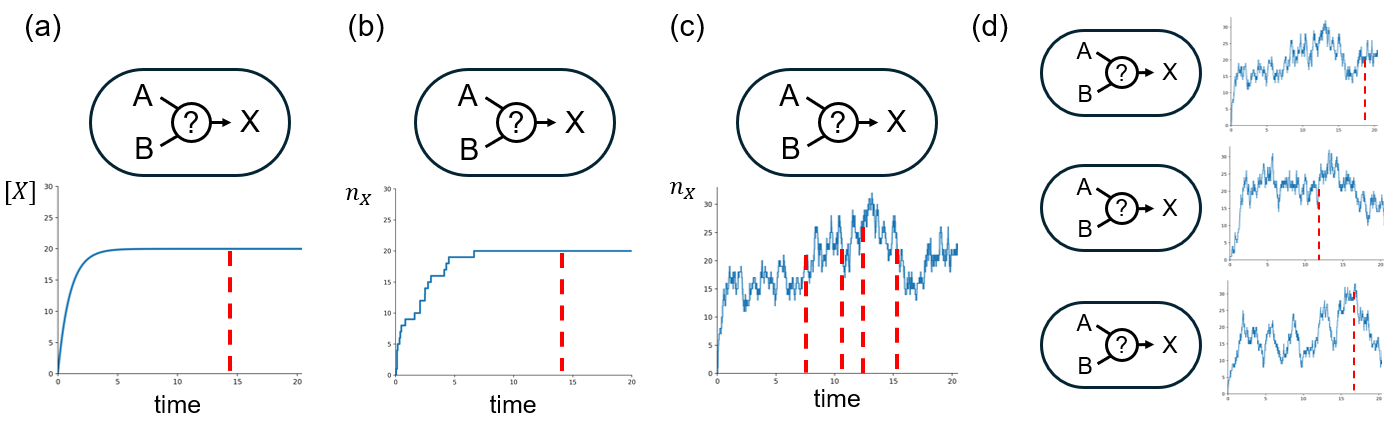}
    \caption{
     Comparison of chemical reaction network computing frameworks.
(a) In deterministic reaction network computation, the output is encoded by the equilibrium concentration of a designated species. After the system reaches equilibrium, the output concentration, $[X]$, can be read from a single measurement (red dashed line).
(b) In count-based stochastic reaction network computation, the output is encoded by a stable or eventual molecular count. After the system stabilizes, the output count, $n_X$, is read from a single measurement (red dashed line). Many such constructions  require exact molecular counts for initialization.
(c) In the stationary-mean framework studied here, the computed value is encoded by the stationary mean of an output species. After the system approaches stationarity, the output is estimated from multiple measurements of the copy number. If these measurements are taken along a single trajectory, then the sampling error depends on the variance and autocovariance of the output process.
(d) If instead the output is estimated from independent copies of the same system, then the sampling error is determined by the variance of the output species.}
    \label{fig:schematic-figure}
\end{figure}

Encoding values by stationary means changes how the result of a computation must be read compared with the deterministic and count-based viewpoints above.  Since the output of a stochastic reaction network at a particular time is a random variable, the encoded value cannot be obtained from a single concentration measurement, as in deterministic ODE-based computation (Fig.~\ref{fig:schematic-figure}a), nor from a single measurement of a stable or eventual molecular count, as in count-based computation (Fig.~\ref{fig:schematic-figure}b).  Instead, the stationary mean must be estimated from samples of the stochastic process (Fig.~\ref{fig:schematic-figure}c,d), and the accuracy of this readout depends on properties of the underlying process. Thus, in this framework, the computation is not fully specified by the stationary mean alone: one must also understand the sampling cost required to read that mean with a desired accuracy.

To make this precise, let $X$ be the species representing the result of the computation. Let $\bar X_N$ be the empirical average of $N$ observations of $X$ at times $t_1,\dots,t_N$, and let $x$ be the encoded stationary mean value: 
\begin{align*}
    \bar X_N = \frac{1}{N}\sum_{i=1}^N X(t_i), \quad x = \E_\pi [X].
\end{align*}
If the measurements come from a single trajectory (Fig.~\ref{fig:schematic-figure}c), then the mean-squared error can be decomposed as follows:  
\[
\E\left[\left(x-\bar X_N\right)^2\right]
=
\underbrace{\left(
x-\frac1N\sum_{i=1}^N \E[X(t_i)]
\right)^2}_{\text{squared bias}}
+
\underbrace{\frac1{N^2}\sum_{i,j=1}^N
\operatorname{Cov}(X(t_i),X(t_j))}_{\text{sampling error}}.
\]
On the other hand, if the observations are taken from $N$ independent copies of the network, denoted $\{X_i\}$ (Fig.~\ref{fig:schematic-figure}d), so that $\bar X_N = \frac{1}{N}\sum_{i=1}^N X_i(t_i)$, then we have
\[
\E\left[\left(x-\bar X_N\right)^2\right]
=
\underbrace{\left(
x-\frac1N\sum_{i=1}^N \E[X_i(t_i)]
\right)^2}_{\text{squared bias}}
+
\underbrace{\frac1{N^2}\sum_{i=1}^N
\var(X_i(t_i))}_{\text{sampling error}}.
\]
These two sources of error are controlled by different properties of the stochastic reaction network. The squared bias is governed by convergence of the output mean to its stationary value, while the sampling error depends on the output variance and autocovariance. For a fixed number of observations and fixed sampling lags, the sampling error generally persists even as the bias vanishes. Thus, we study both how quickly the two terms approach their stationary values and the magnitude of the sampling error that remains at stationarity. Specifically, we derive stationary variance and autocovariance formulas in Section~\ref{sec:mean_level}, and analyze convergence rates of the output mean, variance, and autocovariance in Section~\ref{sec:convergence}.

This readout perspective also affects how we compare different implementations of the same arithmetic expression.  By an architecture, we mean the particular feedforward arrangement of modules used to carry out the computation.  For example, an iterated sum or product may be implemented sequentially, through a linear chain of modules, or hierarchically, through a binary tree.  
These architectures can have the same stationary mean at the output, and therefore compute the same value at the mean level, while nevertheless having different mean convergence times, variances, and autocovariances. Thus, in stochastic mean-level computation, architectural choices can affect both the time required for the mean to approach its stationary value and the cost of reading that mean accurately. We show that these architectural effects are operation-dependent. For iterated addition, binary-tree architectures give better asymptotic upper bounds on both mean and autocovariance convergence times and smaller stationary variances than linear chains (Theorems~\ref{thm:6.5}--\ref{thm:6.9}). For iterated multiplication, both architectures have convergence-time upper bounds of the same order, but the variance comparison can reverse. For example, for the computation of $a^n$ with fixed $0<a\le1/2$, the binary-tree stationary variance is exponentially larger than the linear-chain variance as $n=2^m\to\infty$ (Theorems~\ref{thm:multiplication-binary-tree} and~\ref{thm:fixed-time-aucotovariance-bounds-for-linear-chain-multiplication}). This reversal was unexpected. Thus, architectural depth alone can be misleading when evaluating stochastic readout performance.

The goal of this paper is to develop this stationary-mean viewpoint for elementary arithmetic computations and to analyze the readout costs that arise from the stochastic dynamics.  We focus on identification, addition, and multiplication.  These operations are not intended to exhaust the class of arithmetic functions that can be implemented by stochastic reaction networks at the mean level; rather, they are simple enough to analyze explicitly while already exhibiting the central phenomena being studied. For each module, we establish mean-level correctness under the appropriate input assumptions, including independence for multiplication. We derive stationary variance and autocovariance formulas and quantify convergence of the output mean and autocovariance to their stationary values. Finally, we analyze identification chains and compare convergence-time bounds and stationary fluctuations for binary-tree and linear-chain implementations of iterated addition and multiplication.

This paper is organized as follows.
Section~\ref{sec:def_and_notation} introduces stochastic chemical reaction networks and mass-action kinetics. Section~\ref{sec:modules} presents the computational modules and architectures. Section~\ref{sec:theoretical-tools} develops conditional Poisson representations and moment formulas, and establishes conditions ensuring that stationary distributions with finite moments are preserved under feedforward composition. Section~\ref{sec:convergence-mean} establishes mean-level correctness and analyzes stationary fluctuations and convergence of the output mean and autocovariance. Section~\ref{sec:architecture} studies how the selected architectures affect convergence-time bounds and stationary fluctuations. Finally, Section~\ref{sec:discussion} discusses implications and future directions. Longer proofs are deferred to the appendices.

\section{Mathematical models}
\label{sec:def_and_notation}

In this section, we introduce the basic mathematical model, chemical reaction networks, and then specify the stochastic mass-action dynamics used throughout the paper.  Let $\mathcal{S}=\{S_1,\dots,S_d\}$ be the set of chemical species.  A reaction channel, which will correspond to a possible transition of the stochastic model below, is written
\begin{equation*}
\nu_j \to \nu'_j,\qquad j\in\{1,\dots,J\},
\end{equation*}
 where $\nu_j$ and $\nu'_j$ are non-negative integer combinations of the species. For example, $\nu_1$ could be $S_1 + S_2$ and $\nu_1'$ could be $S_3$ in which case $\nu_1 \to \nu_1'$ represents the reaction channel $S_1 + S_2 \to S_3$. Equivalently, we identify $\nu_j$ and $\nu'_j$ with vectors in $\Z_{\ge 0}^d$.  With this identification, the vector
\[
\zeta_j=\nu'_j-\nu_j \in\Z^d,
\]
is then the net change associated with reaction $j$: when reaction $j$ fires, the state of the stochastic model changes by $\zeta_j$.  The formal linear combinations $\nu_j$ and $\nu'_j$ are called complexes, and we write $\mathcal C$ and $\mathcal R$ for the sets of complexes and reactions, respectively.  Thus, a chemical reaction network is specified by the species, complexes, and reactions, often denoted by the triple $(\mathcal{S},\mathcal{C},\mathcal{R})$.

 A given CRN can be modeled in several ways, including by ordinary differential equations, stochastic differential equations driven by Brownian motion, or continuous-time Markov chains (CTMCs). In this paper, we use the continuous-time Markov chain model. This is the standard stochastic model underlying the chemical master equation and the Gillespie simulation algorithm ~\cite{gillespie1976general, gillespie1977exact, anderson2015stochastic}.  When the process is in state $n = (n_1 , \cdots, n_d) \in \mathbb{Z}_{\ge 0} ^d$, reaction $j$ occurs during the interval $[t,t+\Delta t)$, and the model is updated via the reaction vector $\zeta_j$,  with probability
 \[
    \lambda_j (n)\Delta t + o(\Delta t), \quad \text{ as } \quad \Delta t \to 0,
\]
for a specified intensity function $\lambda_j:\Z^d_{\ge 0} \to \R_{\ge 0}$. Equivalently, the process is a  continuous-time Markov chain  with  infinitesimal generator~\cite{anderson2015stochastic}
\begin{align*}
    \mathcal{A}f(n) = \sum_{j = 1} ^{J} \lambda_j (n) (f(n + \zeta_j) - f(n)).
\end{align*}

There are many ways to define the intensity function $\lambda_j (n)$. In this paper, we  consider \textit{mass-action kinetics}, where $\lambda_j (n)$ is defined as
\begin{align*}
    \lambda_j (n) = \begin{cases}
        k_j \prod_{i=1} ^d \frac{n_i !}{(n_i - \nu_{j,i})!}, & n_i \ge \nu_{j,i} \quad \text{ for all $i$},\\
        0, & \text{otherwise}.
        \end{cases}
\end{align*}
for some reaction rate constant $k_j \in \mathbb{R}_{> 0}$. Throughout this paper, we place the rate constants next to the corresponding reaction arrow, such as $\xrightarrow{k_j}$. 

As an example, consider the following reactions:
\begin{align*}
\begin{split}
    \emptyset &\xrightleftharpoons[1]{3} A,\\
    A &\xrightarrow{1} A +X,\\
    X &\xrightarrow{1} \emptyset.
\end{split}
\end{align*}
Then, the corresponding CRN can be defined as
\begin{align*}
    \begin{split}
        \mathcal{S} & =\{A, X \},\\
        \mathcal{C} &= \{\emptyset, A, X,A+X \},\\
        \mathcal{R} &= \{ \emptyset \rightarrow A,~ A \rightarrow \emptyset,~ A\rightarrow A + X,~ X \rightarrow \emptyset \}.
    \end{split}
\end{align*}
The mass-action intensities are
\begin{align*}
    \begin{split}
        \lambda_{\emptyset \xrightarrow{} A} (n_A, n_X) &= 3,\\
        \lambda_{A \xrightarrow{} \emptyset} (n_A , n_X) &= n_A,\\
        \lambda_{A \xrightarrow{} A + X} (n_A , n_X) &= n_A,\\
        \lambda_{X \xrightarrow{} \emptyset} (n_A , n_X) &= n_X.
    \end{split}
\end{align*}

Throughout this paper, uppercase letters generally denote species or the associated random variables. Lowercase letters denote reaction parameters and numbers encoded at the stationary mean level. Notation such as $n_X$ denotes the copy number of species $X$. We write $\pi$ for a stationary distribution.

\section{Basic modules and architectures}
\label{sec:modules}

In this section, we introduce the three elementary modules studied in this paper: identification, addition, and multiplication. The proofs of mean-level correctness and the analysis of readout properties are deferred to later sections. We also introduce schematic diagrams for the modules and illustrate how the same arithmetic expression can be implemented by different network architectures.

\subsection{Identification, addition, and multiplication modules}

Throughout, the input species of a computational module are understood to be generated by an upstream reaction network and are not altered by the reactions of the module itself. Thus, adjoining a computational module leaves the upstream dynamics unchanged: the input species act catalytically, affecting the intensity functions of reactions that produce the output species. We refer to this by saying that the upstream process evolves autonomously.

\subsubsection{Identification}

Let $A$ be the input species, and let $X$ be the output species. We define the identification module by
\begin{align}\label{eq : identification-module}
\begin{split}
    A  &\xrightarrow[]{1} A + X,\\
    X & \xrightarrow[]{1} \emptyset.
\end{split}
\end{align}

\subsubsection{Addition}

Let $A$ and $B$ be the input species, and let $X$ be the output species. We define the addition module by
\begin{align}\label{eq : single-addition-module}
\begin{split}
    A &\xrightarrow[]{1} A + X,\\
    B &\xrightarrow[]{1} B + X,\\
    X & \xrightarrow[]{1} \emptyset.
\end{split}
\end{align}

\subsubsection{Multiplication}

Let $A$ and $B$ be the input species, and let $X$ be the output species. We define the multiplication module by
\begin{align}\label{eq : single-multiplication-module}
\begin{split}
    A + B &\xrightarrow[]{1} A + B + X,\\
    X & \xrightarrow[]{1} \emptyset.
\end{split}
\end{align}
Independence of the input processes will be assumed when establishing multiplication of stationary means in Sections~\ref{sec:convergence-mean} and~\ref{sec:architecture}.

\subsection{Schematic representations and architectures}

Throughout this paper, we use schematic diagrams to represent computational modules. The operation performed by a module is indicated inside a circle. Lines entering the circle from the left represent the input species, while arrows exiting the circle on the right represent the output species of the computation. See Figure \ref{fig:operations} below.

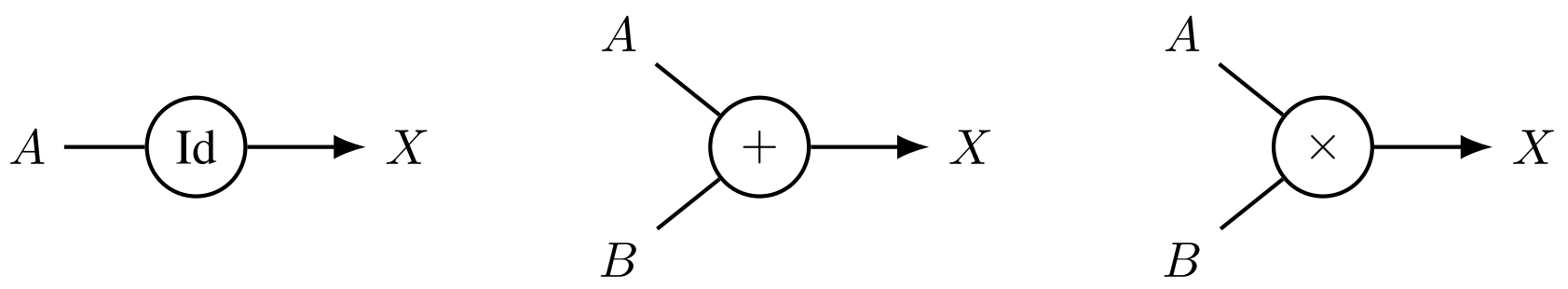
\begin{figure}[H]
\centering
\begin{tikzpicture}[
  thick, >=Latex,
  sum/.style={draw,circle,minimum size=7mm,inner sep=0pt},
  every node/.style={font=\normalsize}
]
  \node (A) at (-1.2, 2) {$A$};
  \node (X) at (1.5, 2) {$X$};

  \node[sum] (s1) at (0,2) {Id};

  \draw (A) -- (s1);
  \draw[->] (s1) -- (X);

  \node (A1) at (3, 2.8) {$A$};
  \node (B1) at (3, 1.2) {$B$};
  \node (X1) at (5.5, 2) {$X$};

  \node[sum] (s11) at (4,2) {$+$};

  \draw (A1) -- (s11);
  \draw (B1) -- (s11);
  \draw[->] (s11) -- (X1);

  \node (A2) at (7, 2.8) {$A$};
  \node (B2) at (7, 1.2) {$B$};
  \node (X2) at (9.5, 2) {$X$};

  \node[sum] (s12) at (8,2) {$\times$};

  \draw (A2) -- (s12);
  \draw (B2) -- (s12);
  \draw[->] (s12) -- (X2);

\end{tikzpicture}
\caption{Schematic representations of the identification (left), addition (center), and multiplication (right) modules. The diagrams correspond to equations \eqref{eq : identification-module}, \eqref{eq : single-addition-module}, and \eqref{eq : single-multiplication-module}, respectively.}
\label{fig:operations}
\end{figure}

\begin{figure}[H]
\centering
\begin{tikzpicture}[
  thick, >=Latex,
  sum/.style={draw,circle,minimum size=7mm,inner sep=0pt},
  every node/.style={font=\normalsize}
]
  \node (A) at (-1.4, 2) {$\emptyset \xrightleftharpoons[1]{a} A$};
  \node (X) at (1.5, 2) {$X$};

  \node[sum] (s1) at (0,2) {Id};

  \draw (A) -- (s1);
  \draw[->] (s1) -- (X);

  \node (A1) at (2.5, 2.8) {$\emptyset \xrightleftharpoons[1]{a} A$};
  \node (B1) at (2.5, 1.2) {$\emptyset \xrightleftharpoons[1]{b} B$};
  \node (X1) at (5.5, 2) {$X$};

  \node[sum] (s11) at (4,2) {$+$};

  \draw (A1) -- (s11);
  \draw (B1) -- (s11);
  \draw[->] (s11) -- (X1);

  \node (A2) at (6.5, 2.8) {$\emptyset \xrightleftharpoons[1]{a} A$};
  \node (B2) at (6.5, 1.2) {$\emptyset \xrightleftharpoons[1]{b} B$};
  \node (X2) at (9.5, 2) {$X$};

  \node[sum] (s12) at (8,2) {$\times$};

  \draw (A2) -- (s12);
  \draw (B2) -- (s12);
  \draw[->] (s12) -- (X2);

\end{tikzpicture}
\caption{Examples in which the input species are generated by birth-death systems.}
\label{fig:specified-examples}
\end{figure}
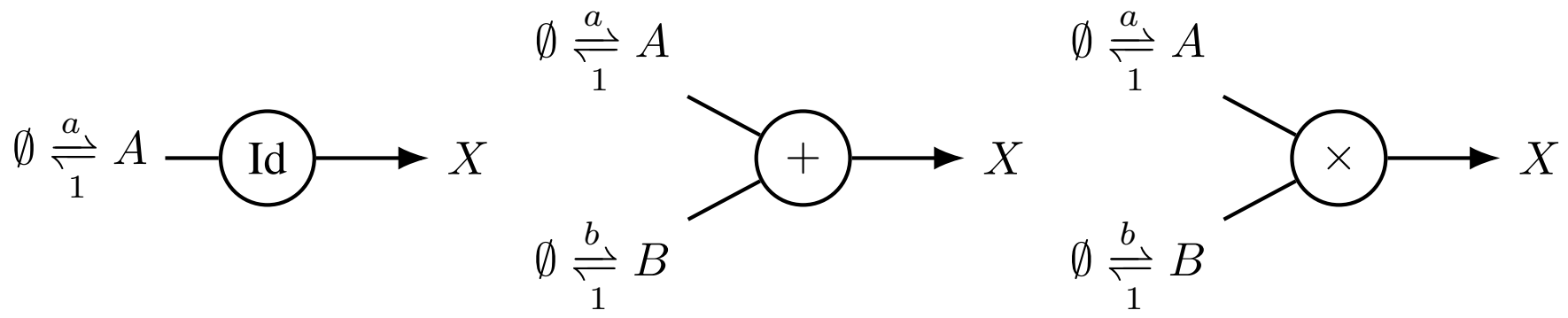

The input species need not have any particular stationary distribution. When an input is generated by a birth-death system, we indicate that system explicitly in the diagram (Fig.~\ref{fig:specified-examples}); otherwise, we leave the upstream network unspecified.

By composing these modules, one can represent the computation of more complex expressions. Suppose, for example, that $A$, $B$, and $C$ are generated by an upstream reaction network and encode $a$, $b$, and $c$, respectively, with the process $C$ independent of the joint process $(A,B)$. Then the expression $(a+b)\times c$ can be realized by the architecture shown below, with final output species $Y$.

\begin{figure}[H]
\centering
\begin{tikzpicture}[
  thick, >=Latex,
  sum/.style={draw,circle,minimum size=7mm,inner sep=0pt},
  every node/.style={font=\normalsize}
]
  \node (A) at (-1, 2.8) {$A$};
  \node (B) at (-1, 1.2) {$B$};
  \node (E) at (1.5,-0) {$C$};

  \node[sum] (s1) at (0,2) {$+$};
  \node[sum] (s3) at (2.5,1) {$\times$};

  \coordinate (Xtip) at (1.0,2.0);
  \coordinate (Xstart) at (1.8,2.0);

  \draw (A) -- (s1);
  \draw (B) -- (s1);
  \draw (E) -- (s3);

  \draw[->] (s1) -- (Xtip)
    node[pos=1, anchor=west, xshift=4pt] {$X$};
  \draw (Xstart) -- (s3);

  \draw[->] (s3) -- ++(1.5,0) node[right] {$Y$};
\end{tikzpicture}
\caption{Schematic representation of the computation of $(a+b) \times c$. Species $A,B,$ and $C$ encode the values $a,b,$ and $c$, respectively, at the level of stationary means. Then, the intermediate species $X$ encodes the quantity $a+b$, and the final output species $Y$ encodes $(a+b) \times c$ at the stationary mean level.}
\end{figure}
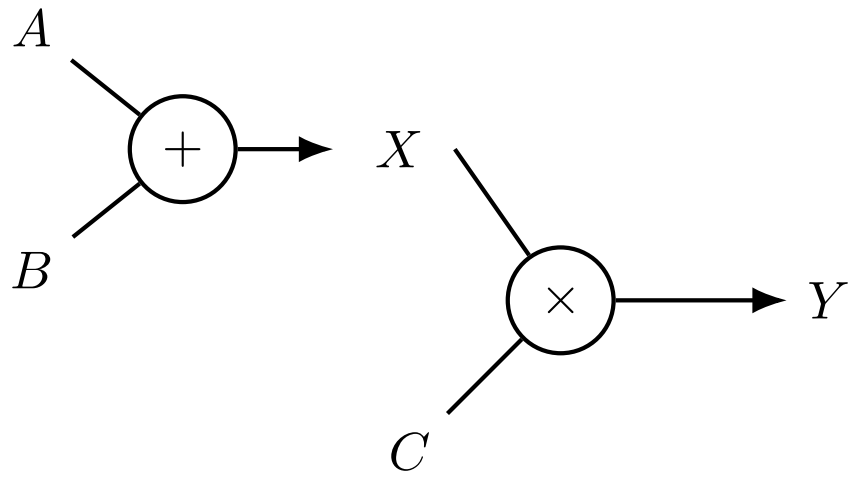
An important point is that the same arithmetic expression can be implemented by various reaction-network architectures. For example,  $(a+b) \times c$ can also be computed using the following architecture.
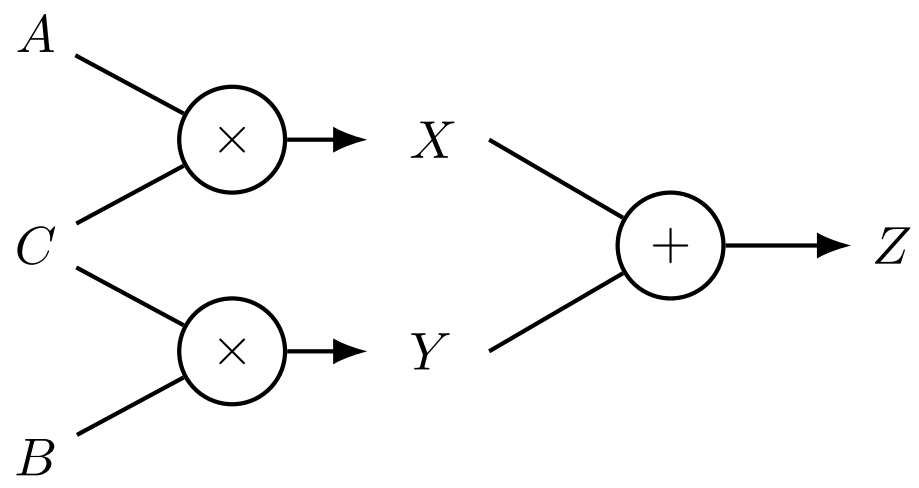
\begin{figure}[H]
\centering
\begin{tikzpicture}[
  thick, >=Latex,
  module/.style={draw,circle,minimum size=7mm,inner sep=0pt},
  every node/.style={font=\normalsize}
]

  \node (A) at (-1, 2.9) {$A$};
  \node (C) at (-1, 1.5) {$C$};
  \node (B) at (-1, 0.1) {$B$};

  \node[module] (m1) at (0.3, 2.2) {$\times$};
  \node[module] (m2) at (0.3, 0.8) {$\times$};
  \node[module] (p1) at (3.2, 1.5) {$+$};

  \coordinate (Xtip) at (1.2, 2.2);
  \coordinate (Ytip) at (1.2, 0.8);

  \coordinate (Xstart) at (2.0, 2.2);
  \coordinate (Ystart) at (2.0, 0.8);

  \draw (A) -- (m1);
  \draw (C) -- (m1);
  \draw (C) -- (m2);
  \draw (B) -- (m2);

  \draw[->] (m1) -- (Xtip)
    node[pos=1, anchor=west, xshift=4pt] {$X$};
  \draw (Xstart) -- (p1);

  \draw[->] (m2) -- (Ytip)
    node[pos=1, anchor=west, xshift=4pt] {$Y$};
  \draw (Ystart) -- (p1);

  \draw[->] (p1) -- ++(1.2,0) node[right] {$Z$};

\end{tikzpicture}
\caption{Schematic representation of the computation of $(a + b) \times c$. Species $A$, $B$, and $C$ encode the values $a$, $b$, and $c$, respectively, at the level of stationary means. The intermediate species $X$ and $Y$ encode $a\times c$ and $b\times c$, respectively, and the output species $Z$ encodes $(a + b) \times c$ at the stationary-mean level.}
\end{figure}

\begin{figure}[H]
\centering
\begin{tikzpicture}[
  thick, >=Latex,
  sum/.style={draw,circle,minimum size=7mm,inner sep=0pt},
  every node/.style={font=\normalsize}
]

  \node (A1) at (-2, 3) {$\emptyset \xrightleftharpoons[1]{a_1 } A_{0,1}$};
  \node (A2) at (-2, 1) {$\emptyset \xrightleftharpoons[1]{a_2} A_{0,2}$};
  \node (A3) at (-2, 0) {$\vdots$};
  \node (A4) at (-2, -1) {$\emptyset \xrightleftharpoons[1]{a_{2^{m} - 1} } A_{0,2^{m} -1}$};
  \node (A5) at (-2, -3) {$\emptyset \xrightleftharpoons[1]{a_{2^{m}} } A_{0,2^{m}}$};

  \node (B1) at (1, 2) {$A_{1,1}$};
  \node (B2) at (1, 0) {$\vdots$};
  \node (B3) at (1, -2) {$A_{1, 2^{m -1}}$};

  \node (C1) at (2, 1.5) {$\ddots$};
  \node (C2) at (2, -1.5) {$\iddots$};

  \node (D1) at (3, 1) {$A_{m -1, 1}$};
  \node (D2) at (3, -1) {$A_{m -1 , 2}$};
  \node (E1) at (6, 0) {$A_{m , 1}$};

  \node[sum] (s1) at (-1,2) {Op};
  \node[sum] (s2) at (-1,-2) {Op};
  \node[sum] (s3) at (4.5,0) {Op};

  \coordinate (Xtip) at (2.0,2.0);
  \coordinate (Ytip) at (2.0,-2.0);

  \coordinate (Xstart) at (2.8,2.0);
  \coordinate (Ystart) at (2.8,0.0);

  \draw (A1) -- (s1);
  \draw (A2) -- (s1);
  \draw (A4) -- (s2);
  \draw (A5) -- (s2);

  \draw[->] (s1) -- (B1);
  \draw (D1) -- (s3);

  \draw[->] (s2) -- (B3);
  \draw (D2) -- (s3);

  \draw[->] (s3) -- (E1);
\end{tikzpicture}
\caption{Binary-tree architecture for iterated operations. Here `Op' denotes either addition ($+$) or  multiplication ($\times$).}
\label{graph: addition-binary-tree}
\end{figure}
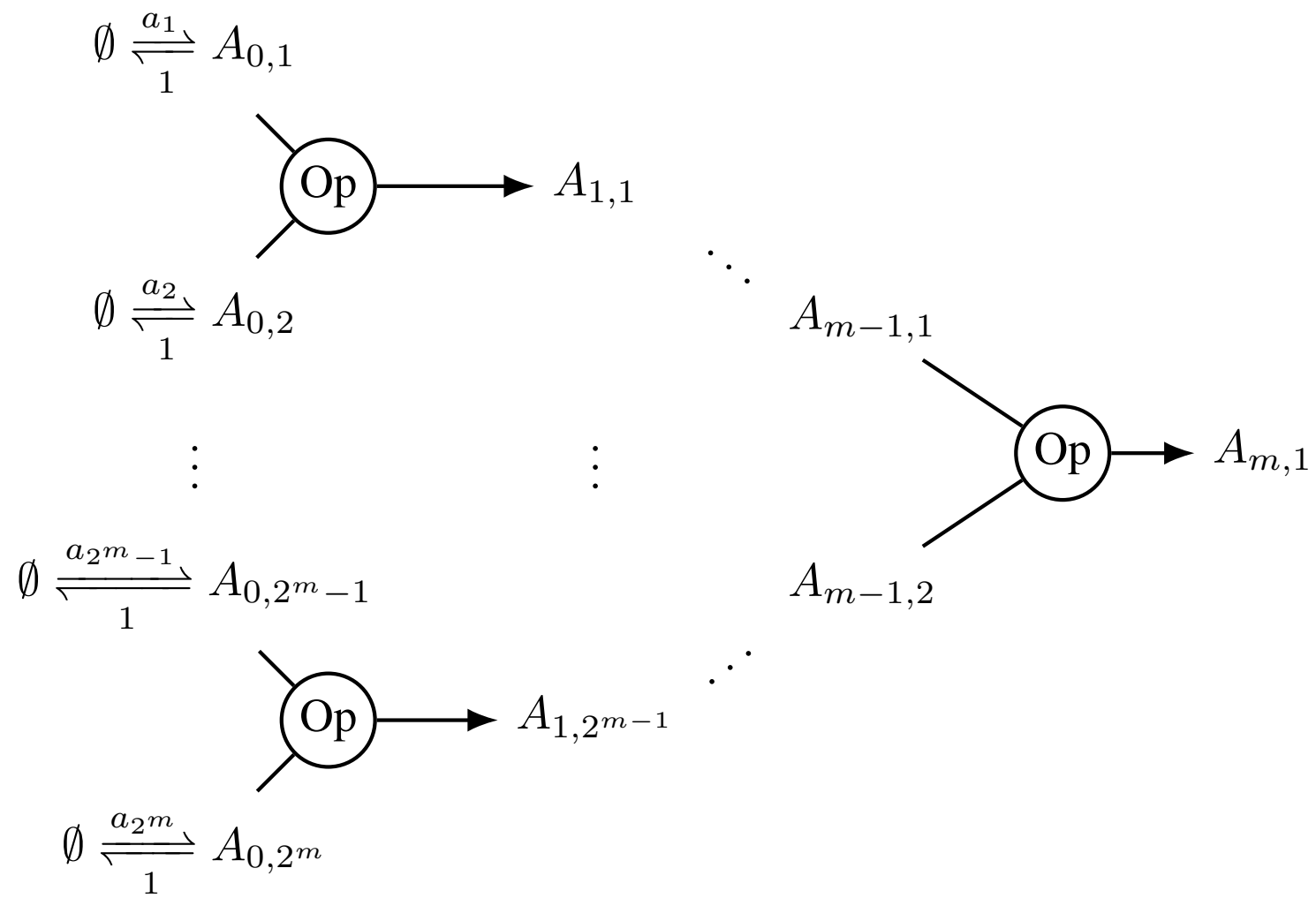

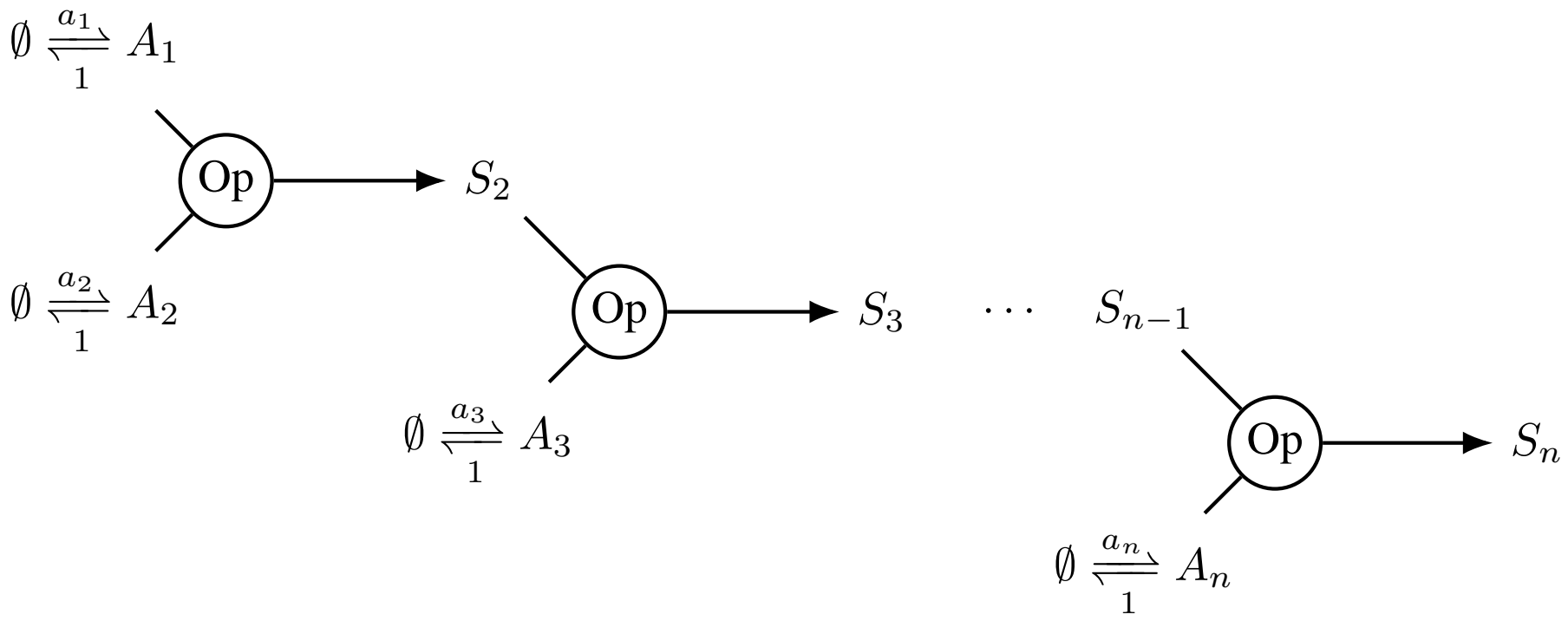
\begin{figure}[H]
\centering
\begin{tikzpicture}[
  thick, >=Latex,
  sum/.style={draw,circle,minimum size=7mm,inner sep=0pt},
  every node/.style={font=\normalsize}
]

  \node (A1) at (-2, 3) {$\emptyset \xrightleftharpoons[1]{a_1 } A_{1}$};
  \node (A2) at (-2, 1) {$\emptyset \xrightleftharpoons[1]{a_2 } A_{2}$};
  \node (A3) at (1, 0) {$\emptyset \xrightleftharpoons[1]{a_3 } A_{3}$};
  \node (A4) at (6, -1) {$\emptyset \xrightleftharpoons[1]{a_n } A_{n}$};

  \node (B1) at (1, 2) {$S_{2}$};
  \node (B2) at (4,1) {$S_3$};
  \node (B3) at (5,1) {$\cdots$};
  \node (B4) at (6,1) {$S_{n-1}$};
  \node (B5) at (9,0) {$S_{n}$};

  \node[sum] (s1) at (-1,2) {Op};
  \node[sum] (s2) at (2,1) {Op};
  \node[sum] (s3) at (7,0) {Op};

  \coordinate (Xtip) at (2.0,2.0);
  \coordinate (Ytip) at (2.0,-2.0);

  \coordinate (Xstart) at (2.8,2.0);
  \coordinate (Ystart) at (2.8,0.0);

  \draw (A1) -- (s1);
  \draw (A2) -- (s1);
  \draw (B1) -- (s2);
  \draw (A3) -- (s2);

  \draw[->] (s1) -- (B1);
  \draw[->] (s2) -- (B2);
  \draw[->] (s3) -- (B5);
  \draw (B4) -- (s3);
  \draw (A4) -- (s3);

\end{tikzpicture}
\caption{Linear-chain architecture for iterated operations. Here `Op' denotes either addition ($+$) or multiplication ($\times$).}
\label{graph: addition-linear-chain}
\end{figure}

Architectural choices become particularly important for larger composite calculations. To study architectural effects in a simple and systematic setting, we consider iterated addition and multiplication of $n$ input values. Two natural choices are a balanced binary tree, which combines inputs pairwise at each level (Fig.~\ref{graph: addition-binary-tree}), and a linear chain, which incorporates one input at a time (Fig.~\ref{graph: addition-linear-chain}). For the binary-tree architecture, we assume $n=2^m$. For multiplication, we take the input processes to be mutually independent. We compare these architectures in Section~\ref{sec:architecture}.

The main issue is not merely whether an architecture computes the desired quantity in stationary mean, but how the architecture used to realize the computation affects its performance. In particular, different architectures may agree at the level of stationary means while exhibiting substantially different mean convergence times, variances, and autocovariances, and hence different readout costs.

\section{Analysis of a general random-input immigration-death system}
\label{sec:theoretical-tools}

Before analyzing the individual modules, we take a brief, but necessary, analytical detour to develop a collection of results that will be used repeatedly throughout the paper. Note that all three modules share the same underlying  structure: an upstream process determines the birth intensity of a new species, while that species undergoes linear degradation. Thus, each module can be viewed as a random-input immigration-death process, or equivalently, in the jargon of queuing systems, as a Cox-driven infinite-server queue. Rather than analyze this structure separately for each module, we study it once in a general setting. To ensure that the resulting properties propagate through repeated feedforward compositions, we will also introduce a notion of admissibility that guarantees the existence of a stationary distribution and the moment bounds needed later in the paper.

\begin{defn}
    Let $Z=(Z(t))_{t\ge0}$ be a $d$-dimensional continuous-time Markov chain. We say that $Z$ is \textit{admissible} if
    \begin{enumerate}
        \item $Z$ is non-explosive,
        \item $Z$ admits a stationary distribution $\pi$, and
        \item $\E_\pi[(1+\|Z\|_1)^p]<\infty$ for every integer $p\ge 1$.
    \end{enumerate}
\end{defn}

We consider systems of the form $(Z,X)$, where $Z$ is an admissible CTMC with state space $E\subseteq\Z_{\ge0}^d$ and $X$ is governed by
\begin{align}\label{eq:generalized-version}
    \emptyset \xrightarrow{\lambda(Z(t))}X, \quad X \xrightarrow{\gamma} \emptyset,
\end{align}
for some $\gamma>0$. We assume that $Z$ evolves autonomously (that is, its dynamics do not depend on $X$), and that $\lambda$ satisfies the following condition.

\begin{cond}\label{condition:lambda}
    The function $\lambda:E\to\mathbb R_{\ge0}$ satisfies this condition if there exist a constant $C>0$ and an integer $m\ge0$ such that
    \begin{align*}
        \lambda(z) \le C (1+||z||_1) ^m, \quad z \in E.
    \end{align*} 
\end{cond}

Taking $\lambda(z)=z_A$, $\lambda(z)=z_A+z_B$, or $\lambda(z)=z_Az_B$ gives the birth intensities used in the identification, addition, and multiplication modules, respectively, and each choice satisfies Condition~\ref{condition:lambda} (with $m=1,1,$ and $2$, respectively).

For such Cox-driven infinite-server systems, the conditional-Poisson representations used below are classical, and related formulations for queues in random environments can be found in
~\cite{jansen2016large, koops2017networks, boxma2019infinite, dean2020functional, dean2022noise}. Since the formulations in these works are slightly different from the one required here, we state and prove the corresponding results.

\subsection{Finite-time conditional Poisson representation and moments}

We begin with the finite-time representation of the random-input immigration--death process. The result is standard and follows from binomial survival of the initial population together with Poisson thinning of the subsequent arrivals.

Recall that we say a  stochastic intensity $\lambda(Z(\cdot ))$ is \textit{locally integrable} if, almost surely,
\[
\int_0^T \lambda(Z(s))\,ds<\infty
\]
for every finite $T>0$.

\begin{prop}[Finite-time conditional Poisson representation]\label{thm:finite-time-cox-process-representation}
Let $Z$ be a stochastic process and suppose that $\lambda(Z(\cdot))$ is non-negative and locally integrable.
Consider the process $X$ governed by
\[
\emptyset \xrightarrow{\lambda(Z(t))} X,
\qquad
X \xrightarrow{\gamma} \emptyset,
\]
for some $\gamma>0$. Assume that $Z$ evolves autonomously. Then, for every $t\ge0$, 
\[
X(t)=X_{\mathrm{old}}(t)+X_{\mathrm{new}}(t),
\]
where, conditional on $X(0)$ and the path $\{\lambda(Z(s)):0\le s\le t\}$,
\[
X_{\mathrm{old}}(t)\sim \operatorname{Binomial}\!\left(X(0),e^{-\gamma t}\right),
\]
and
\[
X_{\mathrm{new}}(t)\sim \operatorname{Poisson}\!\left(M_t^{(0)}\right),
\qquad
M_t^{(0)}:=\int_0^t e^{-\gamma(t-s)}\lambda(Z(s))\,ds,
\]
and moreover $X_{\mathrm{old}}(t)$ and $X_{\mathrm{new}}(t)$ are conditionally independent.
\end{prop}

Note that if $Z$ is a non-explosive continuous-time Markov chain and $\lambda$ is finite-valued on the state space, then $\lambda(Z(\cdot))$ is locally integrable almost surely. This follows because  on every bounded time interval $Z$ visits only finitely many states, and hence $\lambda(Z(t))$ is bounded there. Thus, admissibility of $Z$ together with Condition~\ref{condition:lambda} is more than sufficient for the local integrability required in the proposition above.

\begin{proof}
Decompose the molecules present at time $t$ into those already present at time $0$, denoted $X_{\text{old}}(t)$, and those born during $(0,t]$, denoted $X_{\text{new}}(t)$, so that
\[
X(t)=X_{\text{old}}(t)+X_{\text{new}}(t).
\]

Given $X(0)$, each of the $X(0)$ initial molecules survives independently until time $t$ with probability $e^{-\gamma t}$. Hence
\[
X_{\mathrm{old}}(t)\sim \operatorname{Binomial}\!\left(X(0),e^{-\gamma t}\right).
\]

Conditional on the path $\{\lambda(Z(s)):0\le s\le t\}$, births during $(0,t]$ form an inhomogeneous Poisson process with intensity $s\mapsto\lambda(Z(s))$. A molecule born at time $s$ survives until time $t$ with probability $e^{-\gamma(t-s)}$, independently of the other molecules and of the birth process. Thus, thinning the birth process by the time-dependent retention probability $e^{-\gamma(t-s)}$ gives
\[
X_{\text{new}}(t) \sim
 \operatorname{Poisson} \left(
\int_0^t e^{-\gamma(t-s)} \lambda(Z(s)) ds
\right)
=
\operatorname{Poisson} \left(M_t^{(0)}\right).
\]

Finally, the remaining lifetimes of the molecules present at time $0$ are independent of the birth process after time zero, and of the lifetimes of molecules born after time $0$. Hence $X_{\text{old}}(t)$ and $X_{\text{new}}(t)$ are conditionally independent given $X(0)$ and $\{\lambda(Z(s)):0\le s\le t\}$.
\end{proof}

\begin{cor}[Finite-time mean, variance, and autocovariance] \label{cor:finite-time-moments}

Under the assumptions and notation of Proposition~\ref{thm:finite-time-cox-process-representation}, suppose additionally that
\[
\E[X(0)^2]<\infty \qquad \text{and} \qquad \E[(M_T^{(0)})^2] < \infty
\]
for every finite $T>0$. Then, for $0\le s\le t$,
    \begin{align*}
        \E [X(t)] &= e^{- \gamma t} \E [X(0)] + \E [M_t ^{(0)}],\\
        \var(X(t) ) &= e^{- \gamma t} (1 - e^{- \gamma t}) \E [X(0)] + \E [M_t ^{(0)}] + e^{-2 \gamma t} \var(X(0)) + \var(M_t ^{(0)}) + 2 e^{- \gamma t} \cov (X(0), M_t ^{(0)}),\\
        \cov(X(s), X(t)) &= e^{- \gamma (t-s)} \left(e^{- \gamma s} (1 - e^{- \gamma s}) \E [X(0)] + \E [M_s ^{(0)}] \right) + \cov (e^{- \gamma s} X(0)+M_s ^{(0)}, e^{-\gamma t} X(0) + M_t ^{(0)}).
    \end{align*}
    In particular, if $X(0)$ is independent of the stochastic process $\{\lambda(Z(s)):s\ge 0\}$, and hence of $M_t^{(0)}$ for every $t\ge 0$, then
\begin{align*}
\var(X(t))
&=
e^{-\gamma t}(1-e^{-\gamma t})\E[X(0)]
+\E[M_t^{(0)}]
+e^{-2\gamma t}\var(X(0))
+\var(M_t^{(0)}),
\\
\cov(X(s),X(t))
&=
e^{-\gamma(t-s)}
\left(
e^{-\gamma s}(1-e^{-\gamma s})\E[X(0)]
+\E[M_s^{(0)}]
\right) + e^{-\gamma(s+t)}\var(X(0))
+\cov(M_s^{(0)},M_t^{(0)}).
\end{align*}
\end{cor}

\begin{proof}
    For $t\ge 0$, let
\[
\F_t:=\sigma\left(X(0),\{\lambda(Z(u)):0\le u\le t\}\right)
\]
denote the $\sigma$-algebra generated by $X(0)$ and the intensity path up to time $t$. Also, note that all quantities below are finite under the second-moment assumptions in the statement of the corollary.

    From Proposition \ref{thm:finite-time-cox-process-representation}, conditioning on $X(0)$ and the intensity path gives
\[
\E\left[X(t) \ | \ \F_t \right] = e^{-\gamma t}X(0)+M_t^{(0)}.
\]
Taking expectations therefore yields
\[
\E[X(t)]
=
e^{-\gamma t}\E[X(0)]+\E[M_t^{(0)}].
\]

For the variance, Proposition \ref{thm:finite-time-cox-process-representation} and the conditional independence of $X_{\text{old}}(t)$ and $X_{\text{new}}(t)$ give
\[
\begin{aligned}
\var\left(X(t) \ | \  \F_t \right) &= \var\left(X_{\text{old}}(t) \ | \  \F_t \right)  + \var\left(X_{\text{new}}(t) \ | \  \F_t \right) \\
&= X(0)e^{-\gamma t}(1-e^{-\gamma t})+M_t^{(0)}.
\end{aligned}
\]
Hence, by the law of total variance,
\begin{align*}
\var(X(t)) &= \E\left[\var\left(X(t) \ | \ \F_t \right)\right]  + \var\left(\E\left[X(t) \ | \  \F_t \right]\right)\\
&=
e^{-\gamma t}(1-e^{-\gamma t})\E[X(0)]
+\E[M_t^{(0)}]
+\var\left(e^{-\gamma t}X(0)+M_t^{(0)}\right)\\
&=
e^{-\gamma t}(1-e^{-\gamma t})\E[X(0)]
+\E[M_t^{(0)}]
+e^{-2\gamma t}\var(X(0))
+\var(M_t^{(0)})
+2e^{-\gamma t}\cov(X(0),M_t^{(0)}).
\end{align*}

We turn now to the covariance.  First, since the evolution of $X$ up to time $s$ depends on the intensity path only through its restriction to $[0,s]$, the same formulas  as above (for $X(t)$) hold when conditioning on $\mathcal F_t$:
\[
\E[X(s)\mid\mathcal F_t]
=
e^{-\gamma s}X(0)+M_s^{(0)},
\]
and
\[
\var(X(s)\mid\mathcal F_t)
=
X(0)e^{-\gamma s}(1-e^{-\gamma s})+M_s^{(0)}.
\]

We write
\[
J_{s,t}:=\int_s^t e^{-\gamma(t-u)}\lambda(Z(u))\,du.
\]
The same survival and thinning argument, applied from time $s$ to time $t$, gives
\begin{align}
\label{eq:conditional_X}
\E\left[X(t) \ | \ X(s),\mathcal F_t\right]
=
e^{-\gamma(t-s)}X(s)+J_{s,t}.
\end{align}
Since $X(s)$ is included in the conditioning and $J_{s,t}$ is
$\mathcal F_t$-measurable, the tower property gives
\begin{align*}
\E\left[X(s)X(t) \ | \ \mathcal F_t\right]
&=
\E\left[
X(s)\E\left[X(t) \ | \ X(s),\mathcal F_t\right]
\ | \ \mathcal F_t
\right]\\
&=
e^{-\gamma(t-s)}\E\left[X(s)^2 \ | \ \mathcal F_t\right]
+
J_{s,t}\E\left[X(s) \ | \ \mathcal F_t\right].
\end{align*}
Also, taking conditional expectations of the formula \eqref{eq:conditional_X}  given $\mathcal F_t$,
and again using the tower property, yields
\[
\E\left[X(t) \ | \ \mathcal F_t\right]
=
e^{-\gamma(t-s)}\E\left[X(s) \ | \ \mathcal F_t\right]+J_{s,t}.
\]
Subtracting the product of the conditional means therefore gives
\begin{align*}
\cov\left(X(s),X(t) \ | \ \mathcal F_t\right)
&=
\E\left[X(s)X(t) \ | \ \mathcal F_t\right]
-
\E\left[X(s) \ | \ \mathcal F_t\right]
\E\left[X(t) \ | \ \mathcal F_t\right]\\
&=
e^{-\gamma(t-s)}
\left(
\E\left[X(s)^2 \ | \ \mathcal F_t\right]
-
\E\left[X(s) \ | \ \mathcal F_t\right]^2
\right)\\
&=
e^{-\gamma(t-s)}
\var\left(X(s) \ | \ \mathcal F_t\right)\\
&=
e^{-\gamma(t-s)}
\left(
X(0)e^{-\gamma s}(1-e^{-\gamma s})+M_s^{(0)}
\right).
\end{align*}
Finally, by the law of total covariance and the conditional mean formulas,
\begin{align*}
\cov(X(s),X(t))
&=
\E\left[\cov\left(X(s),X(t) \ | \ \mathcal F_t\right)\right] + 
\cov\left(
\E\left[X(s) \ | \ \mathcal F_t\right],
\E\left[X(t) \ | \ \mathcal F_t\right]
\right)\\
&=
e^{-\gamma(t-s)}
\left(
e^{-\gamma s}(1-e^{-\gamma s})\E[X(0)]
+\E[M_s^{(0)}]
\right)\\
&\qquad+
\cov\left(
e^{-\gamma s}X(0)+M_s^{(0)},
e^{-\gamma t}X(0)+M_t^{(0)}
\right).
\end{align*}
Under the additional independence assumption, $\cov(X(0),M_r^{(0)})=0$ for every $r\ge0$. The final two formulas therefore follow immediately.
\end{proof}

\subsection{Stationary conditional Poisson representation and moments}

We next construct a stationary version of the random-input immigration--death system by allowing births to occur throughout the infinite past. The resulting conditional Poisson representation yields stationary mean, variance, and autocovariance formulas that will be used to analyze the computational modules.

\begin{prop}[Stationary conditional Poisson representation]\label{thm:stationary-poisson-representation}
    Let $Z$ be an admissible CTMC, let $\pi_Z$ be a stationary distribution satisfying the moment condition in the definition of admissibility, and let $\{Z(t)\}_{t\in\mathbb R}$ denote a two-sided stationary version of $Z$ with marginal distribution $\pi_Z$. Let $\lambda$ be a function satisfying Condition~\ref{condition:lambda}. Consider the immigration-death process $X$ given in \eqref{eq:generalized-version}. Then, there exists a stationary joint process $(Z,X)$ with stationary distribution $\pi$ whose $Z$-marginal is $\pi_Z$, such that, conditioned on the entire path $\{\lambda(Z(s)) :s \in \mathbb{R} \}$, for every $t\in \R$,
    \begin{align*}
        X(t) \sim \operatorname{Poisson} (M_t), \quad M_t := \int_{- \infty} ^{t} e^{- \gamma (t-s)} \lambda (Z(s)) ds = \int_0 ^{\infty} e^{- \gamma u} \lambda(Z(t-u)) du,
    \end{align*}
    and $M_t$ is almost surely finite.
\end{prop}

\begin{proof}
By Tonelli's theorem, stationarity, and Condition~\ref{condition:lambda}, there exist a constant $C>0$ and an integer $m \ge 0$ for which
\begin{align*}
    \E[M_t] &= \int_0^\infty e^{-\gamma u}\E[\lambda(Z(t-u))] du = \frac{\E_{\pi_Z}[\lambda(Z)]}{\gamma} \le \frac{C}{\gamma}\E_{\pi_Z}[(1+\|Z\|_1)^m] <\infty.
\end{align*}
 $M_t$ is a non-negative random variable, so $M_t < \infty$ almost surely. 
 
 Conditional on the path $\{Z(s):s\in\mathbb R\}$, let $\{T_k\}$ be the points of an inhomogeneous Poisson process on $\mathbb R$ with intensity $s\mapsto\lambda(Z(s))$. Give each molecule born at $T_k$ an independent exponential lifetime $L_k$ with rate $\gamma$, with these lifetimes independent of $Z$ and the birth process. We define $X(t)$ to be the number of molecules present at time $t$, counting all arrivals from the infinite past that have survived until time $t$:
    \begin{align*}
        X(t) := \sum_k \textbf{1}_{\{T_k \le t < T_k + L_k\} }.
    \end{align*}
    Then, for fixed $t$, Poisson thinning gives
    \begin{align*}
        X(t) \mid \{\lambda(Z(s)) : s \in \mathbb{R} \} \sim \operatorname{Poisson} \left(\int_{- \infty} ^t e^{- \gamma (t-s)} \lambda(Z(s)) ds \right) =: \operatorname{Poisson}(M_t).
    \end{align*}
Since $M_t<\infty$ almost surely, the conditional Poisson representation implies that $X(t)<\infty$ almost surely. By construction, births occur with intensity $\lambda(Z(t))$, while the exponential lifetimes give total death intensity $\gamma X(t)$, as required. Hence, $X$ satisfies \eqref{eq:generalized-version}.
 
 Now we prove stationarity. Define $Z^{(h)} (s) := Z(s+h)$ and $T_k ^{(h)} := T_k - h$. Conditional on $Z^{(h)}$, the shifted birth times $\{T_k^{(h)}\}$ form an inhomogeneous Poisson process with intensity $s\mapsto\lambda(Z^{(h)}(s))$, while the lifetimes $L_k$ remain independent $\operatorname{Exp}(\gamma)$ random variables.
 Under the shifted environment, 
    \begin{align*}
        X^{(h)} (t) = \sum_k \textbf{1}_{\{ T_k ^{(h)} \le t < T_k ^{(h)} + L_k \}} = \sum_{k} \textbf{1}_{ \{ T_k \le t+h < T_k + L_k\} } = X(t+h).
    \end{align*}
    Since $Z^{(h)}$ has the same law as $Z$, and the conditional law of the shifted birth times and their lifetimes given $Z^{(h)}$ agrees with that of the original birth times and lifetimes given $Z$, the joint construction is invariant under time shifts. Consequently,
    \begin{align*}
        (Z(t_j + h), X(t_j + h))_{\{ j=1, \dots, r\} } \overset{d}{=} (Z(t_j), X(t_j))_{\{ j =1, \dots, r\} } ,
    \end{align*}
    for any $t_1, \dots, t_r$ and $h$. Therefore, the finite-dimensional distributions of $(Z,X)$ are invariant under time shifts, so $(Z,X)$ is a stationary process. 
\end{proof}

\begin{cor}[Stationary mean, variance, and autocovariance]\label{lem : cox-process-variance}
For the stationary joint process constructed in Proposition~\ref{thm:stationary-poisson-representation}, define
\[
C_\lambda(t):=\cov_\pi(\lambda(Z(0)),\lambda(Z(t))),
\qquad t\ge0.
\]
Then, 
    \begin{align*}
        \E_\pi [X] &= \frac{\E_\pi[\lambda(Z)]}{\gamma},\\
        \var_\pi (X) &= \frac{\E_\pi[\lambda(Z)]}{\gamma} + \frac{1}{\gamma} \int_0 ^{\infty} e^{- \gamma u} C_\lambda(u) du = \E_\pi [X] + \frac{1}{\gamma} \cov_\pi (\lambda(Z), X),\\
        \cov_{\pi} (X(0), X(t)) &= e^{- \gamma t} \frac{\E_\pi[\lambda(Z)]}{\gamma} + \frac{1}{2\gamma} \int_0 ^{\infty} e^{-\gamma u} (C_\lambda (t+u) + C_\lambda(|t-u|)) du.
    \end{align*}
\end{cor}

\begin{proof}
We first verify that the second moments used below are finite.
By Condition~\ref{condition:lambda} and the admissibility of $Z$,
\[
    \E_{\pi_Z}[\lambda(Z)^2] < \infty.
\]
Moreover, by the Cauchy--Schwarz inequality,
\begin{align*}
    M_t^2 &= \left(\int_0^\infty e^{-\gamma u}\lambda(Z(t-u))\,du \right)^2 \le\left(\int_0^\infty e^{-\gamma u}\,du\right)\left(\int_0^\infty e^{-\gamma u}\lambda(Z(t-u))^2\,du\right).
\end{align*}
Hence, by Tonelli's theorem and stationarity,
\[
    \E[M_t^2] \le \frac{1}{\gamma^2} \E_{\pi_Z}[\lambda(Z)^2] < \infty.
\]
Let
\[
\F:=\sigma\left(\{\lambda(Z(u)):u\in\mathbb R\}\right)
\]
be the $\sigma$-algebra generated by the entire stationary intensity path.  Since the process is stationary and $X(t)\mid\mathcal F\sim\operatorname{Poisson}(M_t)$ from Proposition \ref{thm:stationary-poisson-representation},
\[
    \E_\pi[X^2]
    =
    \E[M_t]+\E[M_t^2]
    <\infty.
\]
Thus, all second moments appearing in the variance and covariance calculations below are finite.  By Proposition~\ref{thm:stationary-poisson-representation},
\begin{equation}
    \label{eq:stationary-conditional-moments}
\E[X(t) \ | \ \F]=M_t,
\qquad
\var(X(t) \ | \ \F)=M_t.
\end{equation}
Hence, by stationarity, the tower property, and the first equality in \eqref{eq:stationary-conditional-moments},
\begin{align*}
    \E_\pi [X] = \E [\E [X(0) \mid \mathcal{F}]] = \E [M_0] = \int_0 ^\infty e^{- \gamma u} \E_\pi [\lambda(Z(-u))] du.
\end{align*}
By stationarity, $\E_\pi [\lambda(Z(-u))] = \E_\pi [\lambda(Z(0))]$, and hence
\begin{align*}
    \E_\pi[X] = \E_\pi [\lambda(Z(0))] \int_0 ^\infty e^{-\gamma u} du = \frac{\E_\pi [\lambda(Z)]}{\gamma}.
\end{align*}
We turn to the variance.  By stationarity, the second equality in \eqref{eq:stationary-conditional-moments}, and the law of total variance, we have
\begin{align*}
    \var_\pi (X) = \E_\pi [\var (X(0) \mid \mathcal{F})] + \var_\pi (\E [X(0) \mid \mathcal{F}])= \E_\pi [M_0] + \var_\pi (M_0)= \frac{\E_\pi [\lambda(Z)]}{\gamma} + \var_\pi (M_0).
\end{align*}
Now,
\begin{align*}
    \var_\pi (M_0) &= \cov_\pi (M_0 , M_0)= \int_0 ^\infty \int_0 ^\infty e^{- \gamma (u+v)} \cov_\pi (\lambda(Z(-u)), \lambda(Z(-v))) \, du \, dv.
\end{align*}
Stationarity implies $\cov_\pi (\lambda(Z(-u)), \lambda(Z(-v))) = C_\lambda(|u-v|)$. Therefore,
\begin{align*}
    \var_\pi (M_0) = 2 \int_{u \ge v} e^{- \gamma (u+v)} C_\lambda(u-v) \, du\, dv= \frac{1}{\gamma} \int_0 ^\infty e^{- \gamma u} C_\lambda(u) \, du.
\end{align*}
Also, since $\lambda(Z(0))$ is $\F$-measurable,
\begin{align*}
    \cov_\pi (\lambda(Z(0)), X(0)) &= \cov_\pi (\lambda(Z(0)), \E [X(0) \mid \mathcal{F}])= \cov_\pi (\lambda(Z(0)), M_0)\\
    &= \int_0 ^\infty e^{- \gamma u} \cov_\pi (\lambda(Z(0)), \lambda(Z(-u)))\,du\\
    &= \int_0 ^\infty e^{- \gamma u} C_\lambda(u)\, du,
\end{align*}
which completes the proof of the variance formula. Finally, for the autocovariance, for $t\ge0$, define
\[ 
J_{0,t} := \int_0^t e^{-\gamma(t-u)} \lambda(Z(u))\,du.
\]
Conditioning on $X(0)$ and $\F$,
the same survival and thinning argument as in the finite-time representation gives
\[
\E[X(t)\mid X(0),\mathcal F] = e^{-\gamma t}X(0)+J_{0,t}.
\]
Since $J_{0,t}$ is $\mathcal F$-measurable,
\begin{align*}
\cov(X(0),X(t)\mid\mathcal F) &= \cov\!\left( X(0), \E[X(t)\mid X(0),\mathcal F] \,\middle|\, \mathcal F \right)\\
&=e^{-\gamma t} \var(X(0)\mid\mathcal F)\\
&=e^{-\gamma t}M_0.
\end{align*}
Hence, by the law of total covariance,
\begin{align*}
    \cov_\pi(X(0),X(t)) &= \E_\pi\!\left[ \cov(X(0),X(t)\mid\mathcal F) \right] + \cov_\pi\!\left( \E[X(0)\mid\mathcal F], \E[X(t)\mid\mathcal F] \right)\\
    &= e^{-\gamma t}\E_\pi[M_0] + \cov_\pi(M_0,M_t) \\
    &= e^{-\gamma t} \frac{\E_\pi[\lambda(Z)]}{\gamma} + \cov_\pi(M_0,M_t).
\end{align*}
As in the variance calculation, we obtain
\begin{align*}
    \cov_\pi(M_0,M_t) &= \int_0^\infty \int_0^\infty e^{-\gamma(u+v)} \cov_\pi\!\left( \lambda(Z(-u)), \lambda(Z(t-v)) \right) \,du\,dv \\
    &= \int_0^\infty \int_0^\infty e^{-\gamma(u+v)} C_\lambda(|t+u-v|) \,du\,dv\\
    &= \int_{u \ge v}e^{-\gamma(u+v)}C_\lambda(|t+u-v|)dudv + \int_{v > u}e^{-\gamma(u+v)}C_\lambda(|t+u-v|)dudv\\
    &= \frac{1}{2\gamma} \int_0^\infty e^{-\gamma u} \left( C_\lambda(t+u) + C_\lambda(|t-u|) \right)\,du.
\end{align*}
Substituting this expression into the preceding formula for $\cov_\pi(X(0),X(t))$ gives the result.
\end{proof}

\begin{rem}
    Note that
    \begin{align*}
        \var_\pi (X) = \E_\pi [X] +\frac{1}{\gamma }\int_0 ^{\infty} e^{-\gamma u} C_{\lambda} (u) du = \E_\pi [X]+\var (M_0 ) \geq \E_\pi [X].
    \end{align*}
    Thus, the stationary variance of the output species is always at least its stationary mean.
\end{rem}

\subsection{Closure property of admissibility}

We now show that adjoining one of our computational modules preserves admissibility. Together with the admissibility of the birth--death input systems, this allows the stationary construction and moment bounds of the previous subsection to be applied repeatedly through a finite feedforward composition of modules.

\begin{thm}\label{thm:admissibility-of-module}
    Let $Z$ be an admissible CTMC on a state space $E \subseteq \mathbb{Z}_{\ge 0} ^d$ with a stationary distribution $\pi_Z$ satisfying the moment condition in the definition of admissibility.  Let $\lambda : E \rightarrow \mathbb{R}_{\ge 0}$ be a real-valued function satisfying Condition~\ref{condition:lambda}. Now, suppose that $X$ is governed by the reactions in \eqref{eq:generalized-version}, with $Z$ evolving autonomously. 
    Then, the process $(Z,X)$ is also admissible.
\end{thm}

\begin{proof} 
We verify the three conditions for admissibility in turn.
\vspace{.1in}

\noindent 
    \textbf{(1) Non-explosivity.}

    Fix an initial state $(z,x) \in E \times \mathbb{Z}_{\ge 0}$, and a finite $T > 0$. Since $Z$ is non-explosive, it has only finitely many jumps on $[0,T]$ almost surely. Since $\lambda$ is finite-valued on $E$,
    \begin{align*}
        \sup_{0 \le s \le T} \lambda(Z(s)) < \infty \quad \text{almost surely}.
     \end{align*}
     Therefore, conditioned on the upstream path, the number of births of $X$ in $[0,T]$ is finite almost surely. Also, the number of deaths of $X$ during $[0,T]$ is at most $X(0) + (\text{number of births during } [0,T])$, and is therefore almost surely finite. Thus, $(Z,X)$ is non-explosive.

     \vspace{.1in}
     
    \noindent \textbf{(2) Existence of a stationary distribution.}

    Since $Z$ admits a stationary distribution $\pi_Z$, by Proposition~\ref{thm:stationary-poisson-representation}, $(Z,X)$ admits a stationary distribution $\pi$.

\vspace{.1in}

\noindent
    \textbf{(3) Finiteness of polynomial moments.}

Now, we work with the stationary joint process $(Z,X)$ with distribution $\pi$ constructed from (2). Fix an integer $p \ge 1$. By Jensen's inequality,
\begin{align*}
    M_t ^p = \left( \int_0 ^\infty e^{-\gamma u} \lambda(Z(t-u)) du \right)^p \le \frac{1}{\gamma^p} \int_0 ^\infty \lambda(Z(t-u))^p \gamma e^{-\gamma u} du.
\end{align*}
Therefore, by stationarity of $Z$,
\begin{align*}
    \E [M_t ^p] \le \frac{1}{\gamma^p} \int_0 ^\infty \gamma e^{-\gamma u} \E_{\pi_Z}[\lambda(Z)^p] du = \frac{1}{\gamma^p} \E_{\pi_Z}[\lambda(Z)^p] < \infty,
\end{align*}
where the last inequality follows from the polynomial upper bound of $\lambda$ and the admissibility of $Z$. Now, once again let
\begin{align*}
    \F := \sigma (\{\lambda(Z(s)) \mid s \in \R \}).
\end{align*}
Then, by Proposition \ref{thm:stationary-poisson-representation}, $X(t) \mid \F \sim \operatorname{Poisson}(M_t)$. Hence,
\begin{align*}
    \E_\pi [X^p] = \E_\pi [\E[X(t)^p \mid \F]] = \E_\pi \left[\sum_{k=1} ^p S(p,k) M_t ^k \right],
\end{align*}
where $S(p,k)$ are the Stirling numbers of the second kind. Since the $p$-th moment of a Poisson random variable is a polynomial of degree $p$ , there exists some constant $C_p > 0$ such that
\begin{align*}
    \sum_{k=1} ^p S(p,k)m^k \le C_p (1+m^p), \quad m \ge 0.
\end{align*}
Thus,
\begin{align*}
    \E_\pi [X^p] \le C_p (1+\E[M_t ^p]) < \infty.
\end{align*}
Finally, since $(1+||Z||_1 + X) ^p \le 2^{p-1} ((1+||Z||_1)^p + X^p)$,  the admissibility of $Z$, together with the previous bounds, gives
\begin{align*}
    \E_\pi [(1+ ||Z||_1 + X)^p] < \infty.
\end{align*}
Since $p \ge 1$ was arbitrary, all polynomial moments under $\pi$ are finite.
\end{proof}

\begin{lem}\label{lem:admissibility-birth-death}
    The birth-death process
    \begin{align*}
        \emptyset \xrightarrow{a} A,\quad A \xrightarrow{1} \emptyset,
    \end{align*}
    is an admissible process.
\end{lem}

\begin{proof}
    This process is the standard immigration-death process (or $M/M/\infty$ queue), and is non-explosive. Indeed, births form a Poisson process with rate $a$, so there are almost surely finitely many births on every finite time interval, while the number of deaths cannot exceed $A(0) + (\text{number of birth events of }A)$, hence almost surely finite. Also, it is known that this process has a stationary distribution $\operatorname{Poisson}(a)$, which has finite polynomial moments.
\end{proof}

\begin{cor}
    The continuous-time Markov chains associated with the structures in Figs.~\ref{graph: addition-binary-tree} and \ref{graph: addition-linear-chain} are admissible.
\end{cor}

\begin{proof}
    This follows from Lemma~\ref{lem:admissibility-birth-death}
    and repeated application of Theorem~\ref{thm:admissibility-of-module} in feedforward order.
\end{proof}

\section{Analysis of single computation modules}
\label{sec:convergence-mean}

We now apply the results of Section \ref{sec:theoretical-tools} to the three computational modules introduced in Section \ref{sec:modules}. We first establish, in Section \ref{sec:mean_level}, their correctness at the level of stationary means and derive formulas for the stationary variance and autocovariance of their outputs. We then study, in Section \ref{sec:convergence}, the rates at which the output means and autocovariances converge to their stationary values.

\subsection{Mean-level correctness and stationary moments}
\label{sec:mean_level}

Throughout this subsection, we assume that the upstream process \(Z\) is admissible and work with the stationary joint distribution \(\pi\) constructed in Proposition \ref{thm:stationary-poisson-representation}. For a stationary scalar process \(Y\), we write \(C_Y(t)=\operatorname{Cov}_\pi(Y(0),Y(t))\), for \(t\geq 0\).

\begin{lem}[Mean-level correctness of computation modules]
    The identification and addition modules \eqref{eq : identification-module} and \eqref{eq : single-addition-module} compute their respective operations correctly at the stationary mean level. Also, if the input species $A,B$ are independent under the stationary distribution $\pi$, the multiplication module \eqref{eq : single-multiplication-module} computes multiplication correctly at the stationary mean level.
\end{lem}

\begin{proof}
    The intensity functions $\lambda(z) = z_A, z_A+ z_B$, and $z_A z_B$ all satisfy Condition~\ref{condition:lambda}. Therefore, by Corollary~\ref{lem : cox-process-variance}, the stationary output means of the three modules are $\E_\pi [A]$, $\E_\pi [A] + \E_\pi [B]$, and $\E_\pi [AB]$, respectively. For the multiplication module, independence under $\pi$ implies $\E_\pi [AB] = \E_\pi [A] \E_\pi [B]$, which completes the proof.
\end{proof}

\begin{lem}[Stationary variance and autocovariance of identification module]
    For a single identification module \eqref{eq : identification-module}, for every $t \ge 0$,
    \begin{align*}
    \begin{split}
        \cov_\pi (X(0), X(t)) &= e^{-t} \E_\pi [A] + \frac{1}{2} \int_0 ^\infty e^{-u} \left( C_A(t+u)+C_A(|t-u|) \right)du,\\
        \var_{\pi} (X) &=   \E _\pi [A] + \cov_\pi (A, X).
    \end{split}
    \end{align*}
\end{lem}

\begin{proof}
    Applying Corollary~\ref{lem : cox-process-variance} with $\gamma = 1$ and $\lambda(Z(t)) = A(t)$ gives the result.
\end{proof}

\begin{lem}[Stationary variance and autocovariance of addition module]\label{lem:5.3}
    For a single addition module \eqref{eq : single-addition-module}, for every $t \ge 0$,
    \begin{align*}
    \begin{split}
        \cov_\pi (X(0), X(t))&= e^{-t} (\E_\pi [A] + \E_\pi [B]) + \frac{1}{2} \int_0 ^\infty e^{-u} (C_{A+B} (t+u) + C_{A+B}(|t-u|)) du,\\
        \var_{\pi} (X) &= \E _\pi [A] + \E_\pi [B] + \cov_\pi (A, X) + \cov_\pi (B, X).
    \end{split}
    \end{align*}
\end{lem}

\begin{proof}
    Applying Corollary~\ref{lem : cox-process-variance} with $\gamma = 1$ and $\lambda(Z(t)) = A(t) + B(t)$, and using the fact that $\cov_\pi (A+B, X) = \cov_\pi (A,X) + \cov_\pi (B,X)$ gives the result.
\end{proof}

\begin{lem}[Stationary variance and autocovariance of multiplication module]\label{lem:5.4}
    For a single multiplication module \eqref{eq : single-multiplication-module}, if the input species $A$ and $B$ are independent under the stationary distribution $\pi$, then for every $t \ge 0$,
    \begin{align*}
    \begin{split}
        \cov_\pi (X(0), X(t)) &= e^{-t} (\E_\pi [A] \cdot \E_\pi [B]) + \frac{1}{2} \int_0 ^\infty e^{-u} (C_{AB} (t+u) + C_{AB} (|t-u|)) du,\\
        \var_{\pi} (X) &= \E _\pi [A] \cdot \E_\pi [B] + \cov_\pi (AB, X).
    \end{split}
    \end{align*}
\end{lem}

\begin{proof}
    Applying Corollary~\ref{lem : cox-process-variance} with $\gamma = 1$ and $\lambda(Z(t)) = A(t)  B(t)$ gives the result.
\end{proof}

\subsection{Convergence of the mean and autocovariance}
\label{sec:convergence}

In this subsection, we study the convergence of the output mean and autocovariance to their stationary values. These convergence properties determine how quickly the bias and the time-dependent sampling error approach their stationary counterparts. We first establish convergence rate estimates and convergence time bounds, and then apply these results to the three individual modules. We adopt the notion of mean convergence rate introduced in~\cite{anderson2025chemical, anderson2026computing}.

\begin{defn}
    Let $\{X(t)\}_{t\ge 0}$ be a family of integrable random variables, and let $x \in \mathbb{R}_{\ge 0}$. Assume that
    \begin{align*}
        |\E[X(t)]-x | \rightarrow 0 \quad \text{as } t\rightarrow\infty.
    \end{align*} 
    Then, the \emph{rate of convergence} of $\E[X(t)]$ to $x$ is defined by
    \begin{align*}
        \rho = -\limsup_{t\rightarrow\infty} \frac{\log |\E[X(t)]- x|}{t}.
    \end{align*}
\end{defn}    

\begin{defn}
    Let $\{X(t)\}_{t\ge 0}$ be a family of square-integrable random variables. Let $C_X: [0,\infty) \rightarrow \R$ be a continuous function. Suppose that
    \begin{align*}
        \sup_{\tau \ge 0}|\cov(X(t), X(t+\tau)) - C_X(\tau)| \rightarrow 0 \quad \text{as } t \rightarrow \infty.
    \end{align*}
    We then define the \emph{rate of convergence} of the autocovariance of $X$ to $C_X$ by
    \begin{align*}
        \sigma = - \limsup_{t \rightarrow \infty} \frac{1}{t}\log  \left(\sup_{\tau\ge0}|\cov(X(t), X(t+\tau)) - C_X(\tau)| \right).
    \end{align*}
    We say that the autocovariance of $X$ converges uniformly to $C_X (\tau)$ at rate $\sigma$.
\end{defn}

\begin{lem}[\cite{anderson2025chemical}, Lemma 3.2]\label{lem:rate-of-convergence}
    Suppose there exist $C>0$, $\rho>0$, and a non-negative polynomial $p(t)$ such that
    \begin{align*}
        |\E[X(t)] - x| \le C\,p(t)e^{-\rho t},
    \end{align*}
    for all sufficiently large $t$. Then $\E [X(t)]$ converges to $x$ at rate at least $\rho$. Conversely, if $\E[X(t)]$ converges to $x$ at a rate greater than $\rho$, then there exist $c>0$ and $t_0\ge 0$ such that
    \begin{align*}
        |\E[X(t)]-x | \le c e^{-\rho t} \quad \text{for all } t\ge t_0.
    \end{align*}
    If $t \mapsto |\E  [ X(t)] - x |$ is bounded on every finite interval, then the bound can be extended to all $t \ge 0$ after increasing $c$ if necessary.
\end{lem}

See \cite{anderson2025chemical} for the proof of the above lemma.

\begin{lem}[Autocovariance convergence rates]\label{lem:covariance-rate-of-convergence}
Suppose there exist constants $C>0$, $\sigma>0$ and a non-negative polynomial $p(t)$ such that
\begin{align*}
    \sup_{\tau \ge 0}|\cov (X(t), X(t+\tau)) - C_X(\tau)| \le C p(t)  e^{- \sigma t} 
\end{align*}
for all sufficiently large $t$. Then, the autocovariance of $X$ converges uniformly to $C_X (\tau)$ at a rate at least $\sigma$. Conversely, if the autocovariance of $X$ converges uniformly to $C_X (\tau)$ at a rate greater than $\sigma$, then there exist $c > 0$ and $t_0 \ge 0$ such that
\begin{align*}
    \sup_{\tau \ge 0}|\cov(X(t), X(t+\tau)) - C_X (\tau)| \le c e^{-\sigma t} \quad \text{for all } t \ge t_0.
\end{align*}
If, for every finite $T>0$,
\begin{align*}
    \sup_{0 \le t \le T} \sup_{\tau \ge 0} |\cov(X(t), X(t+\tau))- C_X (\tau)| < \infty,
\end{align*}
then the inequality holds for every $t \ge 0$, after increasing $c$ if necessary.
\end{lem}

\begin{proof}
    By assumption, we have
    \begin{align*}
        \sup_{\tau\ge0} |\cov(X(t),X(t+\tau))-C_X(\tau)| \le C p(t)e^{-\sigma t} \le C(1+p(t))e^{-\sigma t},
    \end{align*}
    for all sufficiently large $t$. Consequently,
    \begin{align*}
        \limsup_{t\to\infty}\frac{1}{t} \log\left(\sup_{\tau\ge0} |\cov(X(t),X(t+\tau))-C_X(\tau)| \right) \le\limsup_{t\to\infty} \left( \frac{\log C}{t}+\frac{\log(1+p(t))}{t} -\sigma \right) =-\sigma,
    \end{align*}
    since $\log(1+p(t))/t\to0$ as $t\to\infty$. Thus the autocovariance converges uniformly to $C_X$ at a rate at least $\sigma$. Conversely, suppose that the covariance convergence rate is strictly greater than $\sigma$. Choose $\delta>0$ such that
    \begin{align*}
        \limsup_{t\to\infty}\frac{1}{t} \log\left( \sup_{\tau\ge0} |\cov(X(t),X(t+\tau))-C_X(\tau)| \right) <-(\sigma+\delta).
    \end{align*}
    By the definition of the $\limsup$, there exists $t_0>0$ such that
    \begin{align*}
        \sup_{\tau\ge0} |\cov(X(t),X(t+\tau))-C_X(\tau)| \le e^{-(\sigma+\delta)t} \le c e^{-\sigma t} \quad\text{for all }t\ge t_0,
    \end{align*}
    for some $c>0$. Finally, under the additional local boundedness assumption, let
    \begin{align*}
        M:= \sup_{0\le t\le t_0}\sup_{\tau\ge0} |\cov(X(t),X(t+\tau))-C_X(\tau)| <\infty.
    \end{align*}
    Taking $c'=\max\{c,Me^{\sigma t_0}\}$ completes the proof.
\end{proof}

\begin{lem}[Transfer of mean and autocovariance convergence rates]
\label{lem:general-convergence-rate}

Let $X$ be a process governed by the reactions \eqref{eq:generalized-version}, with $Z$ admissible and $\lambda$ satisfying Condition~\ref{condition:lambda}. Fix a stationary distribution $\pi_Z$ of $Z$ satisfying the moment condition in the definition of admissibility, and let $\pi$ be the corresponding joint stationary distribution of $(Z,X)$ constructed in Proposition~\ref{thm:stationary-poisson-representation}. Write
\[
    \bar\lambda:=\E_{\pi_Z}[\lambda(Z)], \quad C_\lambda(\tau):= \cov_{\pi_Z}(\lambda(Z(0)),\lambda(Z(\tau))).
\]
Suppose that $\E[X(0)]<\infty$ and that $\E[\lambda(Z(t))]$ converges to $\bar\lambda$ at a rate greater than $\rho>0$. Assume also that
\[
    \sup_{0\le t\le T}\E[\lambda(Z(t))]<\infty,
\]
for every finite $T>0$. Then $\E[X(t)]$ converges to $\E_\pi[X]$ at a rate at least $\min\{\rho,\gamma\}$. Suppose additionally that $\E[X(0)^2]<\infty$, that
\[
    \sup_{0\le t\le T} \E[\lambda(Z(t))^2]<\infty,
\]
for every finite $T>0$, and that the autocovariance of $\lambda(Z(t))$ converges uniformly to $C_\lambda$ at a rate greater than $\sigma>0$. Then the autocovariance of $X$ converges uniformly to $C_X$ at a rate at least $\min\{\rho,\sigma,\gamma\}$.
\end{lem}

\begin{proof}
See Appendix~\ref{sec:ProofOfLemmaTransfer} for the proof.
\end{proof}

The convergence rates above describe the asymptotic decay of the errors. We now define the convergence times as the times after which the errors remain below a prescribed tolerance level, with the autocovariance measured uniformly over all non-negative time lags.

\begin{defn}
    Let $\pi$ be the stationary distribution of the joint process $(Z,X)$ specified above, and write $x:=\mathbb E_\pi[X]$ and $C_X(\tau):=\operatorname{Cov}_\pi(X(0),X(\tau))$. For each tolerance $\varepsilon>0$, we define the \emph{mean convergence time} $t_{\mathrm{mean}}(\varepsilon)$ as
\begin{align*}
    t_{\mathrm{mean}}(\varepsilon) := \inf\left\{ t\ge 0: \left|\E[X(s)]-x\right|\le \varepsilon \text{ for all } s\ge t \right\}.
\end{align*}
Similarly, we define the \emph{covariance convergence time} $t_{\mathrm{cov}}(\varepsilon)$ as
\begin{align*}
    t_{\mathrm{cov}}(\varepsilon) := \inf\left\{ t\ge 0: \sup_{\tau \ge 0}\left|\cov (X(s) , X(s+\tau))-C_X (\tau)\right|\le \varepsilon \text{ for all } s\ge t \right\}.
\end{align*}
We adopt the convention $\inf\emptyset=+\infty$, so that either convergence time is infinite if no such time exists.
\end{defn}

\begin{lem}[Mean and autocovariance convergence times]
\label{lem:general-mean-convergence-time}
Let $X$ be a process governed by reactions
\eqref{eq:generalized-version}.
Under the assumptions of the first part of
Lemma~\ref{lem:general-convergence-rate}, there exist
constants $C_1>0$ and $\varepsilon_0\in(0,1)$,
independent of $\varepsilon$, such that, for every
$0<\varepsilon\le\varepsilon_0$,
\[
    t_{\mathrm{mean}}(\varepsilon)
    \le
    \frac{C_1}{\min\{\rho,\gamma\}}
    \log_+\left(\frac1\varepsilon\right).
\]

Under the additional assumptions of the second part of Lemma~\ref{lem:general-convergence-rate}, there exists a constant $C_2>0$, independent of $\varepsilon$, such that, for every $0<\varepsilon\le\varepsilon_0$,
\[
    t_{\mathrm{cov}}(\varepsilon)
    \le
    \frac{C_2}{\min\{\rho,\sigma,\gamma\}}
    \log_+\left(\frac1\varepsilon\right),
\]
where $\log_+(x):=\max\{\log x,0\}$.
\end{lem}

\begin{proof}
Set $r_{\mathrm{mean}}:=\min\{\rho,\gamma\}$ and $r_{\mathrm{cov}}:=\min\{\rho,\sigma,\gamma\}$.
Under the respective assumptions, the preceding lemma and the definition of convergence rate imply that the mean and autocovariance errors are bounded by $D_1e^{-r_{\mathrm{mean}}t/2}$ and
$D_2e^{-r_{\mathrm{cov}}t/2}$, respectively, for suitable constants $D_1,D_2>0$ and all $t\ge0$.
Solving these inequalities for $t$ gives bounds involving $\log_+(D_i/\varepsilon)$. For sufficiently small $\varepsilon>0$, these are bounded by constant multiples of $\log_+(1/\varepsilon)$, which gives the claimed bounds.
\end{proof}

We now apply the general results to the three elementary modules.
Throughout the following three lemmas, assume that $Z$ is admissible,
$\E[X(0)]<\infty$, and each input species has locally bounded first
moments. For the autocovariance conclusions, additionally assume
$\E[X(0)^2]<\infty$ and locally bounded second moments of the inputs.

\begin{lem}[Mean and autocovariance convergence of identification]
\label{lem:4.3}
    For a single identification module
    \eqref{eq : identification-module}, suppose that $\E[A(t)]$
    converges to $a$ at a rate greater than $\rho_A$, and that the autocovariance of $A$ converges uniformly to $C_A (\tau)$ at a rate greater than $\sigma_A$. Then $\E[X(t)]$ converges to $a$ at a rate at least $\min \{\rho_A , 1\}$, and the autocovariance of $X$ converges
    uniformly to $C_X$ at a rate at least $\min \{\rho_A , \sigma_A , 1\}$.
    Moreover, there are constants $C_1, C_2>0$ and $\varepsilon_0>0$, all independent of $\varepsilon$, such that for all $0<\varepsilon \le \varepsilon_0$,
    \begin{align*}
        t_{\mathrm{mean}}(\varepsilon)
        \le  
            \frac{C_1}{\min \{\rho_A , 1\}}
            \log\frac{1}{\varepsilon}
        ,\qquad
        t_{\mathrm{cov}}(\varepsilon) \le 
            \frac{C_2}{\min \{\rho_A , \sigma_A ,1\}}
            \log\frac{1}{\varepsilon}.
    \end{align*}
\end{lem}

\begin{proof}
    For every finite $T>0$, the assumptions give
\[
    \sup_{0\le t\le T} \E[\lambda(Z(t))^2] = \sup_{0\le t\le T}\E[A(t)^2] <\infty.
\]
Thus, the required local moment conditions of Lemma~\ref{lem:general-convergence-rate} are satisfied. Since $\lambda(Z(t))=A(t)$, the input intensity has mean and uniform autocovariance convergence rates greater than $\rho_A$ and $\sigma_A$, respectively. Since $\gamma=1$ and $\E_\pi[X]=\E_\pi[A]=a$, the result follows from Lemmas~\ref{lem:general-convergence-rate}
    and~\ref{lem:general-mean-convergence-time}.
\end{proof}

\begin{lem}[Mean and autocovariance convergence of addition]
\label{lem:4.4}
    For a single addition module \eqref{eq : single-addition-module}, suppose that $\E[A(t)]$ and $\E[B(t)]$ converge to $a$ and $b$ at rates greater than $\rho_A$ and $\rho_B$, respectively. Then $\E[X(t)]$ converges to $a+b$ at a rate at least
    $\min \{\rho_A, \rho_B , 1\}$, and
    \[
        t_{\mathrm{mean}}(\varepsilon)
        \le
            \frac{C_1}{\min \{\rho_A , \rho_B, 1\}}
            \log\frac{1}{\varepsilon},
    \]
    for some constant $C_1 > 0$ independent of $\varepsilon$. Suppose additionally that $A$ and $B$ are independent as stochastic processes and that their autocovariances converge uniformly to $C_A (\tau)$ and $C_B (\tau)$ at rates greater than $\sigma_A$ and $\sigma_B$, respectively. Then the autocovariance of $X$ converges uniformly to $C_X (\tau)$ at a rate at least $\min \{\rho_A,\rho_B,\sigma_A,\sigma_B,1 \}$, and
    \[
        t_{\mathrm{cov}}(\varepsilon)
        \le 
            \frac{C_2}{\min \{\rho_A,\rho_B,\sigma_A,\sigma_B,1 \}}
            \log\frac{1}{\varepsilon},
    \]
    for some constant $C_2 > 0$ independent of $\varepsilon$. The convergence-time bounds hold for all $0<\varepsilon\le\varepsilon_0$, for some $\varepsilon_0\in(0,1)$ independent of $\varepsilon$.
\end{lem}

\begin{proof}
    Here $\lambda(Z(t))=A(t)+B(t)$. For every finite $T>0$,
\begin{align*}
\sup_{0\le t\le T} \E[\lambda(Z(t))^2]=\sup_{0\le t\le T}\E[(A(t)+B(t))^2] \le 2\sup_{0\le t\le T}\E[A(t)^2] +2\sup_{0\le t\le T}\E[B(t)^2] <\infty.
\end{align*}
Thus, the required local second-moment condition holds without assuming independence of $A$ and $B$. By the mean convergence assumptions, there exist $c_A,c_B>0$ such that, for all sufficiently large $t$,
    \begin{align*}
        |\E[\lambda(Z(t))]-(a+b)|
        \le |\E[A(t)]-a|+|\E[B(t)]-b| \le c_Ae^{-\rho_A t}+c_Be^{-\rho_B t} \le
        (c_A+c_B)e^{-\min\{\rho_A,\rho_B\}t}.
    \end{align*}
    
Since the input mean convergence rates are strictly greater than $\rho_A$ and $\rho_B$, the exponential bounds above remain valid with slightly larger exponents $\rho_A+\delta$ and $\rho_B+\delta$, for some $\delta>0$. Consequently, the mean convergence rate of $\lambda(Z(t))$ is strictly greater than $\min\{\rho_A,\rho_B\}$.

The mean convergence assertions therefore follow from the first parts of Lemmas~\ref{lem:general-convergence-rate} and~\ref{lem:general-mean-convergence-time}. This argument does not require independence. For the covariance assertions, however, independence of the input processes yields
    \begin{align*}
        \cov(\lambda(Z(t)),\lambda(Z(t+\tau)))=
        \cov(A(t),A(t+\tau))
        +\cov(B(t),B(t+\tau)),
    \end{align*}
    and
    \[
        C_\lambda(\tau)=C_A(\tau)+C_B(\tau).
    \]
    By the uniform covariance convergence assumptions, there exist $d_A,d_B>0$ such that, for all sufficiently large $t$ and every $\tau\ge0$,
    \begin{align*}
        &|\cov(\lambda(Z(t)),\lambda(Z(t+\tau)))
            -C_\lambda(\tau)|\\
        &\quad\le
        |\cov(A(t),A(t+\tau))-C_A(\tau)|
        +|\cov(B(t),B(t+\tau))-C_B(\tau)|\\
        &\quad\le
        d_Ae^{-\sigma_A t}+d_Be^{-\sigma_B t} \le
        (d_A+d_B)e^{-\min\{\sigma_A,\sigma_B\}t}.
    \end{align*}
Since the input autocovariance convergence rates are strictly greater than $\sigma_A$ and $\sigma_B$, the exponential bounds above remain valid with slightly larger exponents. Therefore, the uniform autocovariance convergence rate of $\lambda(Z(t))$ is strictly greater than $\min\{\sigma_A,\sigma_B\}$. Applying the second parts of Lemmas~\ref{lem:general-convergence-rate} and~\ref{lem:general-mean-convergence-time}, with $\gamma=1$, proves the autocovariance assertions.
\end{proof}

\begin{lem}[Mean and autocovariance convergence of multiplication]
    Consider a single multiplication module \eqref{eq : single-multiplication-module} with input processes $A$ and $B$ that are independent. Suppose that $\E[A(t)]$ and $\E[B(t)]$ converge to $a$ and $b$ at rates greater than $\rho_A$ and $\rho_B$, respectively, and that their autocovariances converge uniformly to $C_A (\tau)$ and $C_B (\tau)$ at rates greater than $\sigma_A$ and $\sigma_B$, respectively. Then, $\E[X(t)]$ converges to $ab$ at a rate at least $\min \{\rho_A , \rho_B ,1\}$, and the autocovariance of $X$ converges uniformly to $C_X (\tau)$ at a rate at least $\min \{\rho_A , \rho_B , \sigma_A, \sigma_B , 1\}$. Moreover,
    \begin{align*}
        t_{\mathrm{mean}}(\varepsilon) \le
            \frac{C_1}{\min \{\rho_A , \rho_B ,1\}}
            \log\frac{1}{\varepsilon}
        , \quad
        t_{\mathrm{cov}} (\varepsilon)\le
            \frac{C_2}{\min \{\rho_A , \rho_B , \sigma_A, \sigma_B , 1\}}
            \log\frac{1}{\varepsilon},
    \end{align*}
    for some constants $C_1 , C_2>0$ independent of $\varepsilon$. The convergence-time bounds hold for all $0<\varepsilon\le\varepsilon_0$, for some $\varepsilon_0\in(0,1)$ independent of $\varepsilon$.
\end{lem}

\begin{proof}
    Here $\lambda(Z(t))=A(t)B(t)$.  Since $A$ and $B$ are independent processes, for every finite $T>0$,
\begin{align*}
\sup_{0\le t\le T}\E[\lambda(Z(t))^2] =\sup_{0\le t\le T}
\E[A(t)^2]\E[B(t)^2] \le \left(\sup_{0\le t\le T}\E[A(t)^2]\right)
\left(\sup_{0\le t\le T}\E[B(t)^2]\right)
<\infty.
\end{align*}
Therefore, the required local second-moment condition is satisfied. By the input convergence assumptions, there exist $c_A,c_B,d_A,d_B>0$ such that, for all sufficiently large $t$ and every $\tau\ge0$,
    \begin{align*}
        |\E[A(t)]-a|\le c_Ae^{-\rho_A t},& \quad|\E[B(t)]-b|\le c_Be^{-\rho_B t},\\
        |\cov(A(t),A(t+\tau))-C_A(\tau)|
        \le d_Ae^{-\sigma_A t}, & \quad|\cov(B(t),B(t+\tau))-C_B(\tau)|\le d_Be^{-\sigma_B t}.
    \end{align*}
    By the independence of $A$ and $B$,
    \begin{align}
    \label{eq:multiplication-lambda-upper-bound}
    \begin{split}
        |\E[\lambda(Z(t))]-ab|
        &=
        |\E[A(t)]\E[B(t)]-ab|\\
        &\le
        |\E[A(t)]-a|\,|\E[B(t)]-b|
        +a|\E[B(t)]-b|
        +b|\E[A(t)]-a|\\
        &\le
        (c_Ac_B+a c_B+b c_A)
        e^{-\min\{\rho_A,\rho_B\}t}.
    \end{split}
    \end{align}
    For the covariance, independence of the two input processes gives
    \begin{align*}
    \cov(A(t)B(t),A(t+\tau)B(t+\tau)) &=\cov(A(t),A(t+\tau))
        \cov(B(t),B(t+\tau))\\
        &\hspace{.2in} +
        \cov(A(t),A(t+\tau))
        \E[B(t)]\E[B(t+\tau)]\\
        &\hspace{.2in} +
        \cov(B(t),B(t+\tau))
        \E[A(t)]\E[A(t+\tau)].
    \end{align*}
    The corresponding stationary covariance is
    \[
        C_\lambda(\tau)
        =
        C_A(\tau)C_B(\tau)
        +b^2C_A(\tau)
        +a^2C_B(\tau).
    \]
    The mean bounds imply, uniformly over $\tau\ge0$,
    \begin{align*}
        |\E[A(t)]\E[A(t+\tau)]-a^2|
        \le
        (2a c_A+c_A^2)e^{-\rho_A t},\quad
        |\E[B(t)]\E[B(t+\tau)]-b^2|
        \le
        (2b c_B+c_B^2)e^{-\rho_B t}.
    \end{align*}
    They also imply
    \[
        |\E[A(t)]\E[A(t+\tau)]|\le(a+c_A)^2,
        \quad
        |\E[B(t)]\E[B(t+\tau)]|\le(b+c_B)^2.
    \]
    Moreover, stationarity and Cauchy--Schwarz give
    \[
        |C_A(\tau)|\le\var_\pi(A),
        \quad
        |C_B(\tau)|\le\var_\pi(B).
    \]
    Subtracting the stationary covariance from the finite-time
    covariance and applying these bounds yields
    \begin{align*}
        |\cov(\lambda(Z(t)),\lambda(Z(t+\tau)))
            -C_\lambda(\tau)|
        &\le
        d_A\bigl(\var_\pi(B)+(b+c_B)^2\bigr)
        e^{-\sigma_A t}\\
        &+
        d_B\bigl(\var_\pi(A)+(a+c_A)^2\bigr)
        e^{-\sigma_B t}
        +d_A d_B e^{-(\sigma_A+\sigma_B)t}\\
        &+\var_\pi(B)(2a c_A+c_A^2)e^{-\rho_A t}+\var_\pi(A)(2b c_B+c_B^2)e^{-\rho_B t}.
    \end{align*}
    Consequently, there exists a constant $d_\lambda>0$, independent of $t$ and $\tau$, such that
    \[
        \sup_{\tau\ge0}
        |\cov(\lambda(Z(t)),\lambda(Z(t+\tau)))
            -C_\lambda(\tau)|
        \le d_\lambda e^{-\min\{\rho_A,\rho_B,\sigma_A,\sigma_B\} t}
    \]
    for all sufficiently large $t$. Since all input rate assumptions are strict, the four exponents may be increased by a sufficiently small positive amount in the preceding argument. Thus the mean and uniform covariance convergence rates of $\lambda(Z(t))$ are strictly greater than $\min\{\rho_A,\rho_B\}$ and $\min\{\rho_A,\rho_B,\sigma_A,\sigma_B\}$, respectively. Applying Lemmas~\ref{lem:general-convergence-rate} and~\ref{lem:general-mean-convergence-time} with $\gamma=1$ and $\E_\pi[X]=ab$ proves the result.
\end{proof}

\begin{rem}
The local moment assumptions above are preserved under finite feedforward compositions of the modules considered in the next section, provided that all species have finite initial second moments and that distinct inputs to multiplication modules are independent as processes. Indeed, local second-moment bounds for the input intensity imply local second-moment bounds for the output through the finite-time conditional Poisson representation. The claim therefore follows inductively in feedforward order, starting from the birth--death input species.
\end{rem}

\section{Architectures and readout performance}
\label{sec:architecture}

In this section, we study how selected feedforward architectures affect convergence times and stationary fluctuations. We first analyze identification chains, then compare binary-tree and linear-chain implementations of iterated addition and multiplication, using lumping for the addition analysis. For each architecture, we derive convergence-time bounds and study stationary variance and autocovariance.

\subsection{Identification}
\label{sec:6}

We first consider a sequence of $n$ identification modules with unit rate constants, driven by an input species $A_0$ with stationary mean $a$. The final output species $A_n$ and every intermediate species in the chain have the same stationary mean $a$. We analyze this architecture first because identification chains arise in the analysis of addition architectures through lumping. 

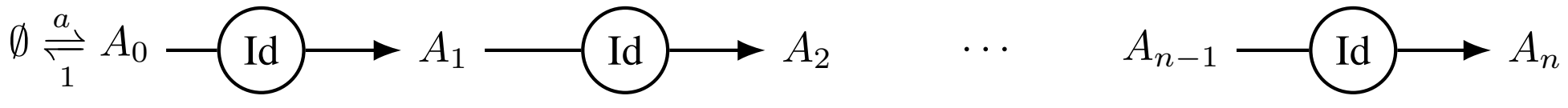
\begin{figure}[H]
\centering
\begin{tikzpicture}[
  thick, >=Latex,
  sum/.style={draw,circle,minimum size=7mm,inner sep=0pt},
  every node/.style={font=\normalsize}
]
  \node (A_0) at (-1.5, 2) {$\emptyset \xrightleftharpoons[1]{a} A_0$};
  \node (A_1) at (1.5, 2) {$A_1$};
  \node (A_2) at (4.5, 2) {$A_2$};
  \node (A_3) at (7.5, 2) {$A_{n-1}$};
  \node (A_4) at (10.5, 2) {$A_n$};

  \node[sum] (s1) at (0,2) {Id};
  \node[sum] (s2) at (3,2) {Id};
  \node (s3) at (6,2) {$\cdots$};
  \node[sum] (s4) at (9,2) {Id};

  \draw (A_0) -- (s1);
  \draw[->] (s1) -- (A_1);
  \draw (A_1) -- (s2);
  \draw[->] (s2) -- (A_2);
  \draw (A_3) -- (s4);
  \draw[->] (s4) -- (A_4);
\end{tikzpicture}
\caption{A chain of $n$ identification modules with source mean $a$.}
\label{eq : identification-sequence}
\end{figure}

\begin{thm}[Mean and autocovariance convergence times for an identification chain]
\label{thm:6.1}
Fix $a>0$. For each $n\ge 1$, consider a sequence of $n$
identification modules as in Fig.~\ref{eq : identification-sequence}.
Suppose there exists a constant $D<\infty$, independent of $n$, such that
\[
    \sup_{0\le i\le n}|\E[A_i(0)]-a|\le D.
\]
Then there exists a constant $C_{\mathrm{mean}}<\infty$,
independent of $n$ and $\varepsilon$, such that
\[
    t_{\mathrm{mean}}(\varepsilon)
    \le C_{\mathrm{mean}}\cdot \left(n+\log\frac{1}{\varepsilon}\right),
\]
for all $n\ge1$ and $0<\varepsilon\le1$.

Additionally, suppose there exists a constant $V<\infty$,
independent of $n$, such that
\[
    \sup_{0\le i\le n}\var(A_i(0))\le V.
\]
Then there exists a constant $C_{\mathrm{cov}}<\infty$,
independent of $n$ and $\varepsilon$, such that
\[
    t_{\mathrm{cov}}(\varepsilon)
    \le C_{\mathrm{cov}}\cdot \left(n+\log\frac{1}{\varepsilon}\right),
\]
for all $n\ge1$ and $0<\varepsilon\le1$.
\end{thm}

\begin{proof}
    Taking expectations in the conditional-Poisson representation of Proposition~\ref{thm:finite-time-cox-process-representation}, and defining $e_i (t) := \E [A_i (t)] - a$, we obtain
    \begin{align*}
        e_i (t) = e^{-t} e_i (0) + \int_0 ^t e^{-(t-s)} e_{i-1} (s) ds, \quad 1 \le i \le n,
    \end{align*}
    where $e_0 (t) = e_0 (0) e^{-t}$. This recursive relation has the solution
    \begin{align*}
        e_n (t) = e^{-t} \sum_{k=0} ^n \frac{t^k}{k!} e_{n-k} (0).
    \end{align*}
    Hence,
    \begin{align*}
        |e_n (t)| \le D e^{-t} \sum_{k=0} ^n \frac{t^k}{k!} = D \mathbb P(\operatorname{Poisson}(t) \le n).
    \end{align*}
    If $D=0$ or $\varepsilon\ge D$, then $t_{\mathrm{mean}}(\varepsilon)=0$. Thus, we may assume $0<\varepsilon<D$. For $t\ge n$, the Chernoff bound gives
\begin{align*}
    D\,\mathbb P(\operatorname{Poisson}(t)\le n) \le D\exp\left(-\frac{(t-n)^2}{2t}\right).
\end{align*}
Set $L:=\log(D/\varepsilon)>0$. Taking $t=2(n+L)$, we have
\[
    \frac{(t-n)^2}{2t} =L+\frac{n^2}{4(n+L)}\ge L.
\]
Since the Chernoff upper bound is nonincreasing for $t\ge n$, we conclude that
\begin{align}\label{eq:t-mean-epsilon-bound}
    t_{\mathrm{mean}}(\varepsilon)
    \le 2n+2\log(D/\varepsilon).
\end{align}
Since $n\ge1$ and $0<\varepsilon\le1$, the claimed uniform bound follows with $C_{\mathrm{mean}}=2(1+\log\max\{1,D\})$.

For the covariance, let us first write $C_i (\tau) := \cov_\pi (A_i (0), A_i (\tau))$. We prove that
    \begin{align}\label{eq:cov-sup-upper-bound-stuff}
        \sup_{\tau \ge 0} |\cov(A_i (t), A_i (t+\tau)) - C_i (\tau)| \le C 6^i e^{-t/2} \quad \text{for all } t \ge 0,
    \end{align}
    where $C :=4(1+a+D+V)$. If $i=0$, applying Corollary \ref{cor:finite-time-moments} gives
    \begin{align*}
        \cov(A_0 (t), A_0 (t+\tau)) - C_0 (\tau) =e^{- \tau} \{ (\E[A_0 (0)] - a) e^{-t} + (\var (A_0 (0)) - \E [A_0 (0)]) e^{-2t} \}.
    \end{align*}
    Therefore,
    \begin{align*}
        \sup_{\tau \ge 0} |\cov(A_0 (t), A_0 (t+\tau)) - C_0 (\tau)| \le D e^{-t} + (V+a+D)e^{-2t} \le C e^{-t/2}.
    \end{align*}
    Now, suppose that the induction bound holds at level $i-1$. We can obtain an analogous bound from \eqref{eq:covariance-convergence-upper-bound} by replacing $X$ with $A_{i}$ and $\lambda(Z)$ with $A_{i-1}$. The mean-error and initial variance terms (first and second terms) are bounded by
    \begin{align*}
        |\E [A_i (t)] - a|, \quad |\var (A_i (0)) - \E [A_i (0)]|e^{-2t}
    \end{align*}
    after taking the supremum over $\tau \ge 0$. For the two initial-input covariance terms (third and fourth terms), the induction hypothesis gives
    \begin{align*}
        \var (A_{i-1} (s)) &\le \var_\pi (A_{i-1}) + C 6^{i-1},\\
        |\cov (A_i (0), A_{i-1} (s))| &\le \sqrt{\var (A_i (0)) (\var_\pi (A_{i-1}) + C 6^{i-1})}.
    \end{align*}
    Thus, these two covariance terms are bounded by
    \begin{align*}
        &\sqrt{\var (A_i (0)) (\var_\pi (A_{i-1}) + C 6^{i-1})} \{ e^{-t} \int_0 ^{t+\tau} e^{-(t+\tau - v)} dv + e^{-(t+\tau)} \int_0 ^t e^{-(t-u)}du\}\\
        &=\sqrt{\var (A_i (0)) (\var_\pi (A_{i-1}) + C 6^{i-1})} \{e^{-t} (1-e^{-(t+\tau)}) + e^{-(t+\tau)} (1-e^{-t}) \}\\
        & \le 2 \sqrt{\var (A_i (0)) (\var_\pi (A_{i-1}) + C 6^{i-1})} e^{-t},
    \end{align*}
    after taking the supremum for $\tau$. For the two stationary-covariance tail integrals (sixth and seventh terms), since $|C_{i-1} (\tau)| \le \var_\pi (A_{i-1})$, we obtain the following bound:
    \begin{align*}
        &\var_\pi (A_{i-1})\cdot \left( \int_{-\infty} ^0 e^{-(t-u)} du \int_{-\infty} ^{t+\tau} e^{-(t+\tau - v) } dv + \int_0 ^t e^{-(t-u)} du \int_{-\infty} ^0 e^{-(t+\tau - v)} dv\right)\\
        &= \var_\pi (A_{i-1}) \cdot ( e^{-t} + (1-e^{-t}) e^{-(t+\tau)} ) \le 2 e^{-t} \var_\pi (A_{i-1}) .
    \end{align*}
    Finally, applying the induction hypothesis to the covariance-difference double integral (fifth term) gives
    \begin{align*}
        &C 6^{i-1} \int_0 ^t \int_0 ^{t+\tau} e^{-(t-u)} e^{-(t+\tau - v)} e^{-\min \{u,v\}/2} dv du 
        \\
        &\le C 6^{i-1} (1+e^{-\tau}) \int_0 ^t e^{-(t-s)} e^{-s/2} ds \le 2 C 6^{i-1} \int_0 ^t e^{-(t-s)} e^{-s/2} ds.
    \end{align*}
    Summing up the upper bounds, we get
    \begin{align*}
        \sup_{\tau\ge 0} |\cov (A_i (t), A_i (t+\tau)) - C_i (\tau)| &\le |\E [A_i (t)] - a| + |\var (A_i (0)) - \E [A_i (0)]|e^{-2t}\\
        & + 2 \sqrt{\var (A_i (0)) (\var_\pi (A_{i-1}) + C 6^{i-1})} e^{-t} + 2 \var_\pi (A_{i-1})e^{-t}\\
        &+2 C 6^{i-1} \int_0 ^t e^{-(t-s)} e^{-s/2} ds.
    \end{align*}
    Note that
    \begin{align*}
        \var_\pi (A_i) = a + \int_0 ^{\infty} e^{-s} C_{i-1} (s) ds \le a + \var_\pi (A_{i-1}),
    \end{align*}
    so we get $\var_\pi (A_i) \le a(i+1)$. We also have the following bound on the mean discrepancy:
\begin{align*}
    |\E[A_i(t)]-a| \le D e^{-t}\sum_{k=0}^{i}\frac{t^k}{k!} \le D2^i e^{-t} \sum_{k=0}^{i}\frac{(t/2)^k}{k!} \le D2^i e^{-t/2}.
\end{align*}
 Using the known bounds
    \begin{align*}
        |\E [A_i (t)] - a| & \le D 2^i e^{-t/2}, \quad \var (A_i (0))  \le V, \quad \E [A_i (0)] \le a+D, \quad \var_\pi (A_{i-1}) \le ai,
        \end{align*}
        we get
        \begin{align*}
            \sup_{\tau\ge 0} |\cov (A_i (t), A_i (t+\tau)) - C_i (\tau)| &\le \{ D 2^i + a + D + V + 2\sqrt{V(ai + C 6^{i-1})} + 2ai + 4C 6^{i-1}\} e^{-t/2}\\
            & \le\{  5 C 6^{i-1} + D 2^i + 3ai + a + D + 2V\} e^{-t/2}\\
            & \le \{5C 6^{i-1} + (4a + 3D + 2V) 6^{i-1} \} e^{-t/2}\\
            &\le C 6^i e^{-t/2}.
        \end{align*}
        This proves the claim at level $i$. Taking $i=n$ in \eqref{eq:cov-sup-upper-bound-stuff}, we obtain, for $0<\varepsilon\le1$,
\[
    t_{\mathrm{cov}}(\varepsilon) \le 2n\log 6+2\log\frac{C}{\varepsilon}
    \le C_{\mathrm{cov}} \left(n+\log\frac1\varepsilon\right),
\]
where $C_{\mathrm{cov}}:=2\log(6C)$ depends only on $a,D,V$ and is independent of $n$ and $\varepsilon$.

\end{proof}

\begin{thm}[Autocovariance of the sequence of identification modules]\label{thm:7.1}
For a sequence of identification modules in Fig.~\ref{eq : identification-sequence} and every $t \ge 0$,
    \begin{align*}
    \begin{split}
        \cov_{\pi} (A_n (0), A_n (t)) &= \E_\pi [A_n ] e^{-t} \sum_{j=0} ^n \sum_{k=0} ^j \frac{\binom{j+k}{k}}{2^{j+k}  (j-k)!} t^{j-k},\\
        \var_\pi (A_n ) &= \E_\pi [A_n ] (2n+1) \frac{\binom{2n}{n}}{4^n} \sim \frac{2}{\sqrt{\pi}} \sqrt{n} ~\E_\pi [A_n],
    \end{split} 
    \end{align*}
    where the asymptotic relation holds as $n \to \infty$.
\end{thm}

\begin{proof}
For $t\geq 0$, define $C_n(t):=\cov_\pi(A_n(0),A_n(t))$. By Corollary~\ref{lem : cox-process-variance},
\begin{align*}
    C_{i+1}(t) = a e^{-t} + \frac12\int_0^\infty e^{-u} \left( C_i(t+u)+C_i(|t-u|) \right)du .
\end{align*}
Define the operator $\mathcal K$ by
\begin{align*}
    (\mathcal K f)(t) := \frac12\int_0^\infty e^{-u} \left( f(t+u)+f(|t-u|) \right)du .
\end{align*}
Then $C_{i+1}(t)=a e^{-t}+(\mathcal K C_i)(t)$. Now set
\begin{align*}
    G_j(t):= e^{-t} \sum_{k=0}^{j} \frac{(j+k)!}{2^{j+k}j!k!(j-k)!} t^{j-k}.
\end{align*}
We claim that $C_n (t) = a \sum_{j=0} ^n G_j (t)$. Note that $G_0(t)=e^{-t}$. Indeed, for every integer \(r\ge 0\),
\[
\mathcal K(e^{-t}t^r) = e^{-t}\left[
\frac{t^{r+1}}{2(r+1)} + \sum_{q=0}^r \binom rq \frac{q!}{2^{q+2}}t^{r-q} + \frac{r!}{2^{r+2}} \right].
\]
Therefore, after substituting the polynomial expression for \(G_j\), a direct coefficient comparison gives, for each \(0\le \ell\le j+1\),
\[
\text{Coefficient of }e^{-t}t^{j+1-\ell}\,\text{ in } (\mathcal KG_j)(t)
=
\frac{(j+1+\ell)!}
{2^{j+1+\ell}(j+1)!\ell!(j+1-\ell)!}.
\]
This is precisely the coefficient of \(t^{j+1-\ell}\) in \(e^tG_{j+1}(t)\). Thus \(\mathcal KG_j=G_{j+1}\). We now prove the formula by induction. For $n=0$, $C_0(t)=a e^{-t}=a G_0(t).$
Suppose the formula holds at level $i$; that is, $C_i(t)=a \sum_{j=0}^{i}G_j(t).$
Then
\begin{align*}
    C_{i+1}(t)= a G_0(t)+(\mathcal K C_i)(t) = a G_0(t)+a \sum_{j=0}^{i}(\mathcal K G_j)(t) = a G_0(t)+ a \sum_{j=0}^{i}G_{j+1}(t) =a \sum_{j=0}^{i+1}G_j(t).
\end{align*}
This completes the induction. Setting $t=0$ gives
\begin{align*}
    \var_\pi(A_n) =a\sum_{j=0}^{n}G_j(0)=a\sum_{j=0}^{n} \frac{\binom{2j}{j}}{4^j}.
\end{align*}
Using the identity $\sum_{j=0}^{n}\frac{\binom{2j}{j}}{4^j} =(2n+1)\frac{\binom{2n}{n}}{4^n}$ gives the stated variance formula. The asymptotic equivalence follows from Stirling's formula.
\end{proof}

\subsection{Addition}

In this subsection, we study compositions of addition modules using the architectures illustrated in Figs.~\ref{graph: addition-binary-tree} and \ref{graph: addition-linear-chain}. The key observation is that compositions of addition modules can be analyzed in terms of identification chains through lumping. For the binary-tree architecture, the sums of the species at each level form a single identification chain. For the linear-chain architecture, the addition network is obtained by lumping a collection of identification chains.

\begin{defn}
Let $X(t)$ be a non-explosive, time-homogeneous continuous-time Markov chain on a countable state space $E$, and let $\mathcal{E}=\{E_i\}_{i\in I}$ be a partition of $E$. Define $h:E\to I$ by $h(x)=i$ whenever $x\in E_i$, and let $Y(t)=h(X(t))$.

We say that $X$ is \emph{strongly lumpable} with respect to $\mathcal{E}$ if, for every initial distribution of $X(0)$, the process $Y$ is a time-homogeneous continuous-time Markov chain with the same transition law, independent of the choice of the initial distribution.
\end{defn}

\begin{lem}\label{lem: lumping}
    (\cite{ball1993lumpability}, Theorem 2.4) Let $\mathcal{E} = \{ E_1 , E_2 , \cdots \}$ be a partition of a countable state space and let $Q=(q_{kl})$ denote the generator of the CTMC from the preceding definition. The CTMC is strongly lumpable with respect to $\mathcal{E}$ if and only if, for every pair $i\neq j$, there exists a constant $\lambda_{i,j} \in \mathbb{R}$ such that
    \begin{align*}
        Q_{ij}\textbf{1} = [q_{kl}] _{k 
        \in E_i, l \in E_j} \textbf{1} = \lambda_{ij} \textbf{1},
    \end{align*}
    where $Q_{ij}$ is a block matrix corresponding to $E_i$ and $E_j$, and $\textbf{1}$ is a column vector of ones of appropriate dimension.
\end{lem}

\subsubsection{Binary-tree architecture}

We first consider the binary-tree addition architecture with $n=2^m$ input species. At each level, we consider the sum of the species copy numbers, from the $n$ inputs at level zero to the final output $A_{m,1}$ at level $m$. We show that the joint process of these level sums forms an identification chain through lumping. This reduction allows us to derive convergence-time bounds and explicit formulas for the stationary variance and autocovariance of the final output.

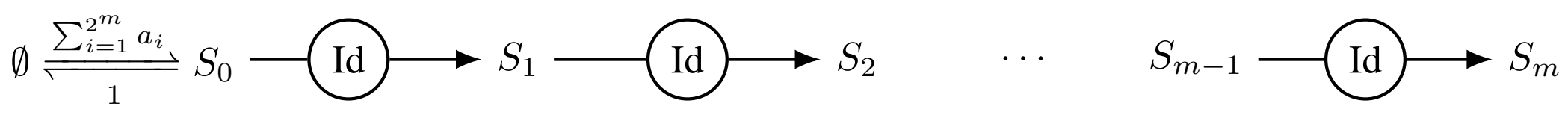
\begin{figure}[htbp]
\centering
\begin{tikzpicture}[
  thick, >=Latex,
  sum/.style={draw,circle,minimum size=7mm,inner sep=0pt},
  every node/.style={font=\normalsize}
]

  \node (A_0) at (-2, 2) {$\emptyset \xrightleftharpoons[1]{\sum_{i=1} ^{2^{m}} a_i } S_0$};
  \node (A_1) at (1.5, 2) {$S_1$};
  \node (A_2) at (4.5, 2) {$S_2$};
  \node (A_3) at (7.5, 2) {$S_{m -1}$};
  \node (A_4) at (10.5, 2) {$S_{m}$};

  \node[sum] (s1) at (0,2) {Id};
  \node[sum] (s2) at (3,2) {Id};
  \node (s3) at (6,2) {$\cdots$};
  \node[sum] (s4) at (9,2) {Id};

  \draw (A_0) -- (s1);
  \draw[->] (s1) -- (A_1);
  \draw (A_1) -- (s2);
  \draw[->] (s2) -- (A_2);
  \draw (A_3) -- (s4);
  \draw[->] (s4) -- (A_4);
\end{tikzpicture}
\caption{Identification-chain representation of binary-tree addition.}
\label{eq : addition-translation}
\end{figure}

\begin{thm}[Lumping of binary-tree addition]\label{thm : addition-binary-tree-translation}
    For the binary-tree addition architecture in Fig.~\ref{graph: addition-binary-tree}, let $x = (x_{i,j})$, where $x_{i,j}:= n_{A_{i,j}}$, and define the level sums
    \begin{align*}
        S_i (x) = \sum_{j=1} ^{2^{m -i}} x_{i,j}, \quad 0 \leq i \leq m .  
    \end{align*}
    Let $\Phi(x) = (S_0 (x) , \cdots , S_{m} (x))$, and let 
    \begin{align*}
        E_s = \{ x: \Phi(x) = s\}, \quad s \in \mathbb{Z}_{\geq 0} ^{m +1}.
    \end{align*}
    Then, the CTMC associated with the binary-tree addition network in Fig.~\ref{graph: addition-binary-tree} is strongly lumpable with respect to the partition $\{E_s  \}$. Moreover, the lumped chain is the CTMC associated with the reaction network shown in Fig.~\ref{eq : addition-translation}.    
\end{thm}

\begin{proof} 
Here and below, $e_0, \ldots, e_m$ denote the standard basis vectors of $\R^{m+1}$, corresponding to the level-sum coordinates $S_0 , \ldots, S_m$, respectively. Death transitions are understood only when the corresponding coordinate is positive; otherwise their rates are zero.

Fix $s\in\mathbb Z_{\geq 0} ^{m+1}$ and $x\in E_s$. By Lemma~\ref{lem: lumping}, it suffices to show that the total transition rate from $x$ to any other lump depends only on $s$. There are three types of reactions.

\smallskip
\noindent \textbf{(i) Births into level $i=0$.} 

Each reaction $0\rightarrow A_{0,j}$ changes $\Phi(x) = s$ to $s+e_0$, and its rate is $a _j$.
Therefore the total rate from $x$ to $E_{s+e_0}$ is $\sum_{j=1}^{2^{m}}a_j $, which depends only on $s$.

\vspace{.1in}
\noindent
\textbf{(ii) Production from level $i-1$ to level $i$, $1\le i\le m$.}

The reactions producing $A_{i,j}$ are
\begin{align*}
    A_{i-1,2j-1}\rightarrow A_{i-1,2j-1}+A_{i,j}, \quad A_{i-1,2j}\rightarrow A_{i-1,2j}+A_{i,j}.
\end{align*}
Each such reaction increases $S_i$ by one and changes no other $S_k$. Hence these reactions induce transitions from $E_s$ to $E_{s+e_i}$. Their total rate is
\begin{align*}
\sum_{j=1}^{2^{\,m-i}}\bigl(x_{i-1,2j-1}+x_{i-1,2j}\bigr) =\sum_{k=1}^{2^{\,m-(i-1)}} x_{i-1,k} = S_{i-1}(x)=s_{i-1},
\end{align*}
again depending only on $s$. 

\vspace{.1in}
\noindent
\textbf{(iii) Degradation at level $i$, $0\le i\le m $.}

Each reaction $A_{i,j}\rightarrow \emptyset$ changes $s$ to $s-e_i$.
The total rate from $x$ to $E_{s-e_i}$ is
\begin{align*}
\sum_{j=1}^{2^{\,m-i}} x_{i,j}=S_i(x)=s_i,
\end{align*}
which depends only on $s$.

\smallskip
No other reaction changes $\Phi(x)$. So, for every pair of distinct lumps $E_s$ and $E_t$, the total transition rate from any $x\in E_s$ to $E_t$ is a constant depending only on $s,t$. Hence, by Lemma~\ref{lem: lumping}, the chain is strongly lumpable with respect to $\{E_s\}$. Moreover, the transition rates computed above coincide with those of the CTMC associated with the reaction network in Fig.~\ref{eq : addition-translation}. Hence, the two CTMCs have the same generator.
\end{proof}

\begin{thm}[Mean and autocovariance convergence times for binary-tree addition]
\label{thm:6.5}
For each $n=2^m$, $m\ge1$, consider the binary-tree addition
architecture in Fig.~\ref{graph: addition-binary-tree}.
Suppose there exists a constant $D<\infty$, independent of $n$,
such that
\[
    \max_{i,j}
    |\E[A_{i,j}(0)]-\E_\pi[A_{i,j}]|\le D.
\]
Then there exists a constant $C_{\mathrm{mean}}<\infty$,
independent of $n$ and $\varepsilon$, such that
\[
    t_{\mathrm{mean}}(\varepsilon)
    \le C_{\mathrm{mean}} \cdot
    \left(\log n+\log\frac1\varepsilon\right)
\]
for all $n=2^m\ge2$ and $0<\varepsilon\le1$.

Additionally, suppose there exist constants $V,A<\infty$,
independent of $n$, such that
\[
    \max_{i,j}\var(A_{i,j}(0))\le V,
    \qquad
    \max_{1\le i\le n}a_i\le A.
\]
Then there exists a constant $C_{\mathrm{cov}}<\infty$,
independent of $n$ and $\varepsilon$, such that
\[
    t_{\mathrm{cov}}(\varepsilon)
    \le C_{\mathrm{cov}} \cdot
    \left(\log n+\log\frac1\varepsilon\right)
\]
for all $n=2^m\ge2$ and $0<\varepsilon\le1$.
\end{thm}

\begin{proof}
    By Theorem~\ref{thm : addition-binary-tree-translation}, the binary-tree addition architecture is strongly lumpable, and the resulting lumped chain is the sequence of identification modules shown in Fig.~\ref{eq : addition-translation}. We use the explicit estimates established in the proof of Theorem~\ref{thm:6.1}. For each level $i$, the stationary mean of $S_i$ is $s := \sum_{j=1} ^n \E_\pi  [A_{0, j}] = \sum_{j=1} ^n a_j$. Define
\begin{align*}
    D_m:=\max_{0\leq i\leq m }\left|\E[S_i(0)]-s\right|.
\end{align*}
Define $a_{i,j} := \E_\pi [A_{i,j}]$. Since $S_i(0)=\sum_{j=1}^{2^{m -i}}A_{i,j}(0)$,  and $s=\sum_{j=1}^{2^{m -i}}a_{i,j}$,
we have
\begin{align*}
    \left|\E[S_i(0)]-s\right| \leq \sum_{j=1}^{2^{m -i}} \left|\E[A_{i,j}(0)]-a_{i,j}\right| \leq D 2^{m -i},
\end{align*}
by the assumption on the initial discrepancies. Hence $D_{m} \le nD$. By \eqref{eq:t-mean-epsilon-bound} and $D_m\le nD$,
\[
t_{\mathrm{mean}}(\varepsilon) \le 2m+2\log\frac{n\max\{1,D\}}{\varepsilon}.
\]
Since $m=\log_2 n$, this is bounded by $C_{\mathrm{mean}}(\log n+\log(1/\varepsilon))$, where $C_{\mathrm{mean}}$ depends only on $D$. For the covariance, the Cauchy-Schwarz inequality gives
    \begin{align*}
        \var (S_i (0)) \le \left( \sum_{j=1} ^{2^{m-i}} \sqrt{\var (A_{i,j} (0))} \right)^2 \le V 2^{2(m-i)} \le V n^2.
    \end{align*}
    Then, using \eqref{eq:cov-sup-upper-bound-stuff},
\begin{align*}
    \sup_{\tau\ge0} |\cov(A_{m,1}(t),A_{m,1}(t+\tau)) - C_{A_{m,1}}(\tau)|  \le 4(1+s+D_m+Vn^2)6^m e^{-t/2}.
\end{align*}
Since $s\le nA$, $D_m\le nD$, and $6^m=n^{\log_2 6}$, the right-hand side is bounded by
\[
    4(1+A+D+V)n^{2+\log_2 6}e^{-t/2}.
\]
Consequently,
\[
    t_{\mathrm{cov}}(\varepsilon)
    \le 2\log\frac{4(1+A+D+V)}{\varepsilon}
    +(4+2\log_2 6)\log n.
\]
This establishes the claimed bound with a constant $C_{\mathrm{cov}}$ depending only on $A,D,V$.
\end{proof}

\begin{thm}[Autocovariance for binary-tree addition]
\label{thm:binary-tree-addition-variance}
    Consider the binary-tree addition architecture with $n = 2^{m}$ inputs. Then, for every $t \ge 0$,
    \begin{align*}
    \begin{split}
        \cov_{\pi} (A_{m,1} (0), A_{m,1} (t)) &= \E_\pi [A_{m,1}] e^{-t} \sum_{j=0} ^{m} \sum_{k=0} ^j \frac{ 
        \binom{j+k}{j}}{2^{j+k} (j-k)!} t^{j-k},\\
        \var_\pi (A_{m,1}) &= \E_\pi [A_{m,1}] (2m + 1)\frac{\binom{2m }{m }}{4^{m}} \sim \frac{2}{\sqrt{\pi}} \sqrt{\log_2 (n)} \E_\pi [A_{m, 1}],
    \end{split} 
    \end{align*}
    where the asymptotic relation holds as $n = 2^m \to \infty$.
\end{thm}

\begin{proof}
By Theorem~\ref{thm : addition-binary-tree-translation}, the binary tree addition structure is strongly lumpable to the sequence of identification modules in Fig.~\ref{eq : addition-translation}. Therefore, the output autocovariance of the original binary-tree architecture is exactly equal to the autocovariance of the final species in an identification chain of length $m$ with stationary input mean $s:= \sum_{j=1} ^n a_j$.

By Theorem~\ref{thm:7.1}, since $S_m=A_{m,1}$ and $\E_\pi[S_m]=s$, we obtain
\begin{align*}
    \cov_\pi(A_{m ,1}(0),A_{m ,1}(t)) = \E_\pi[A_{m ,1}]e^{-t} \sum_{j=0}^{m } \sum_{k=0}^{j} \frac{(j+k)!}{2^{j+k}j!k!(j-k)!} t^{j-k}.
\end{align*}
Setting $t=0$ and using the identity
\begin{align*}
    \sum_{j=0}^{m}\frac{\binom{2j}{j}}{4^j} = (2m+1)\frac{\binom{2m}{m}}{4^m},
\end{align*}
we obtain the variance formula. The asymptotic equivalence follows from Stirling's
formula and $m=\log_2 n$.
\end{proof}

\subsubsection{Linear-chain architecture}

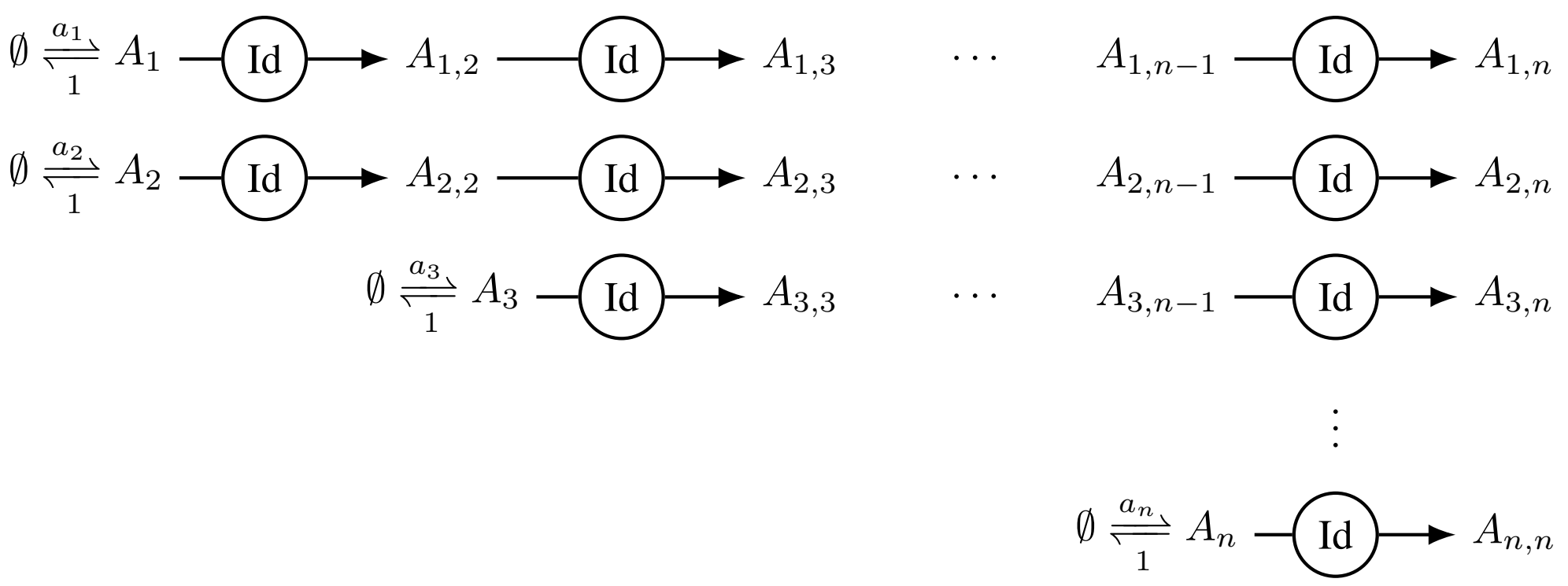
\begin{figure}[htbp]
\centering
\begin{tikzpicture}[
  thick, >=Latex,
  sum/.style={draw,circle,minimum size=7mm,inner sep=0pt},
  every node/.style={font=\normalsize}
]

  \node (A_0) at (-1.5, 2) {$\emptyset \xrightleftharpoons[1]{a_1 } A_1$};
  \node (A_1) at (1.5, 2) {$A_{1,2}$};
  \node (A_2) at (4.5, 2) {$A_{1,3}$};
  \node (A_3) at (7.5, 2) {$A_{1,n-1}$};
  \node (A_4) at (10.5, 2) {$A_{1,n}$};
  \node (B_0) at (-1.5, 1) {$\emptyset \xrightleftharpoons[1]{a_2 } A_2$};
  \node (B_1) at (1.5, 1) {$A_{2,2}$};
  \node (B_2) at (4.5, 1) {$A_{2,3}$};
  \node (B_3) at (7.5, 1) {$A_{2,n-1}$};
  \node (B_4) at (10.5, 1) {$A_{2,n}$};
  \node (C_0) at (1.5, 0) {$\emptyset \xrightleftharpoons[1]{a_3 }A_3$};
  \node (C_1) at (4.5, 0) {$A_{3,3}$};
  \node (C_2) at (7.5, 0) {$A_{3,n-1}$};
  \node (C_3) at (10.5, 0) {$A_{3,n}$};
  \node (C_4) at (9, -1) {$\vdots$};
  \node (D_0) at (7.5, -2) {$\emptyset \xrightleftharpoons[1]{a_n }A_{n}$};
  \node (D_1) at (10.5, -2) {$A_{n,n}$};

  \node[sum] (s1) at (0,2) {Id};
  \node[sum] (s2) at (3,2) {Id};
  \node (s3) at (6,2) {$\cdots$};
  \node[sum] (s4) at (9,2) {Id};

  \node[sum] (s5) at (0,1) {Id};
  \node[sum] (s6) at (3,1) {Id};
  \node (s7) at (6,1) {$\cdots$};
  \node[sum] (s8) at (9,1) {Id};
  \node[sum] (s9) at (3,0) {Id};
  \node (s10) at (6,0) {$\cdots$};
  \node[sum] (s11) at (9,0) {Id};
  \node[sum] (s12) at (9,-2) {Id};

  \draw (A_0) -- (s1);
  \draw[->] (s1) -- (A_1);
  \draw (A_1) -- (s2);
  \draw[->] (s2) -- (A_2);
  \draw (A_3) -- (s4);
  \draw[->] (s4) -- (A_4);
  \draw (B_0) -- (s5);
  \draw[->] (s5) -- (B_1);
  \draw (B_1) -- (s6);
  \draw[->] (s6) -- (B_2);
  \draw (B_3) -- (s8);
  \draw[->] (s8) -- (B_4);
  \draw (C_0) -- (s9);
  \draw[->] (s9) -- (C_1);
  \draw (C_2) -- (s11);
  \draw[->] (s11) -- (C_3);
  \draw (D_0) -- (s12);
  \draw[->] (s12) -- (D_1);
\end{tikzpicture}
\caption{Decomposed representation of the linear-chain addition network}
\label{eq : summation-translation-linear}
\end{figure}

We now consider the linear-chain addition architecture. To analyze this architecture, we introduce an auxiliary network in which each input species generates its own identification chain, as illustrated in Fig.~\ref{eq : summation-translation-linear}. The sums of the species in each column correspond to the partial sums in the original linear-chain addition network. We show that lumping the auxiliary network through these column sums yields the CTMC associated with the original linear-chain addition network.

For each row $i$, let $\ell_1=n-1$ and $\ell_i=n-i+1$ for $2\le i\le n$. We use the aliases
\[
A_{1,1}:=A_1,\qquad A_{i,i-1}:=A_i,\quad 2\le i\le n,
\]
so that the species in row $i$ are indexed by $n-\ell_i\le j\le n$. Each species in this row has stationary mean $a_i$.

\begin{thm}[Lumping of linear-chain addition]\label{thm: addition-linear-translation}

For the decomposition network in Fig.~\ref{eq : summation-translation-linear}, define
\begin{align*}
S_k:=\sum_{i=1}^k A_{i,k},\quad 2\le k\le n,
\end{align*}
and define the projection
\begin{align*}
    \Phi(x)=\bigl(x_1,\dots,x_n,s_2,\dots,s_n\bigr),
\end{align*}
where $x_i = n_{A_i}$ and $s_k=\sum_{i=1}^k n_{A_{i,k}}$. Then the decomposition network is strongly lumpable with respect to the partition induced by $\Phi$. Moreover, the lumped CTMC is the process associated with the linear-chain addition network in Fig.~\ref{graph: addition-linear-chain}.
\end{thm}

\begin{proof}
As in the proof of Theorem~\ref{thm : addition-binary-tree-translation}, we apply Lemma~\ref{lem: lumping} by checking the transition rates between lumps.

For the projected state $(x_1,\ldots,x_n,s_2,\ldots,s_n)$, the birth and death rates of each input coordinate $x_i$ are $a_i$ and $x_i$, respectively. The rates of increase and decrease of $s_2$ are $x_1+x_2$ and $s_2$, respectively. For $3\le k\le n$, the corresponding rates for $s_k$ are $s_{k-1}+x_k$ and $s_k$.

These rates depend only on the projected state, so the decomposition network is strongly lumpable. Moreover, they are exactly the transition rates of the linear-chain addition network in Fig.~\ref{graph: addition-linear-chain}.
\end{proof}

\begin{thm}[Mean and autocovariance convergence times for linear-chain addition]
\label{thm:6.8}
For each $n\ge2$, consider the linear-chain addition architecture in Fig.~\ref{graph: addition-linear-chain}. Let $\Phi$ be the projection defined in Theorem~\ref{thm: addition-linear-translation}. Suppose that the initial distribution of the original network is the pushforward under $\Phi$ of an initial distribution on the auxiliary decomposition network in Fig.~\ref{eq : summation-translation-linear}.

Suppose there exists a constant $D<\infty$, independent of $n$, such that
\[
\max_{1\le i\le n}\max_{n-\ell_i\le j\le n} |\E[A_{i,j}(0)]-a_i|\le D.
\]
Then there exists a constant $C_{\mathrm{mean}}<\infty$, independent of $n$ and $\varepsilon$, such that
\[
    t_{\mathrm{mean}}(\varepsilon) \le C_{\mathrm{mean}}  \cdot \left(n+\log\frac1\varepsilon\right)
\]
for all $n\ge2$ and $0<\varepsilon\le1$.

Additionally, suppose there exist constants $V,A<\infty$, independent of $n$, such that
\[
\max_{1\le i\le n}\max_{n-\ell_i\le j\le n}
\var(A_{i,j}(0))\le V,
\qquad
\max_{1\le i\le n}a_i\le A.
\]
Then there exists a constant $C_{\mathrm{cov}}<\infty$,
independent of $n$ and $\varepsilon$, such that
\[
    t_{\mathrm{cov}}(\varepsilon) \le C_{\mathrm{cov}} \cdot \left(n+\log\frac1\varepsilon\right)
\]
for all $n\ge2$ and $0<\varepsilon\le1$.
\end{thm}

\begin{proof}
    See Appendix~\ref{sec:proof-thm-6.8} for the proof. The argument applies the identification-chain estimates from Theorem~\ref{thm:6.1} to each row of the auxiliary network and controls the cross-row covariances that may arise from correlated initial conditions.
\end{proof}

\begin{thm}[Autocovariance for linear-chain addition]
For the linear chain addition architecture in Fig.~\ref{graph: addition-linear-chain} and every $t\geq 0$,
\begin{align*}
    \cov_\pi(S_n(0),S_n(t)) = \sum_{i=1}^n a_i e^{-t} \sum_{j=0}^{\ell_i} \sum_{k=0}^{j} \frac{(j+k)!}{2^{j+k}j!k!(j-k)!} t^{j-k}.
\end{align*}
In particular,
\begin{align*} 
    \var_\pi(S_n) = \sum_{i=1}^n a_i (2\ell_i+1)\frac{\binom{2\ell_i}{\ell_i}}{4^{\ell_i}},
\end{align*}
where $\ell_1=n-1$ and $\ell_i:=n-i+1$ for $ 2\leq i\leq n$. If, in addition, there exist constants $0<a_{\min}\le a_{\max}<\infty$, independent of $n$, such that $a_{\min}\le a_i\le a_{\max}$ for all $1\le i\le n$, then
\begin{align*}
    \frac{2a_{\min}}{3a_{\max}}\sqrt n\,\E_\pi[S_n]
    \le \var_\pi(S_n)
    \le \frac{\sqrt2\,a_{\max}}{a_{\min}}
        \sqrt n\,\E_\pi[S_n].
\end{align*}
\label{thm:6.9}
\end{thm}

\begin{proof}
Construct the identification chains in the auxiliary decomposition
network independently in stationarity. By
Theorem~\ref{thm: addition-linear-translation}, the product
stationary law projects to a stationary law $\pi$ for the
original linear-chain addition network.

Since $S_n(t)=\sum_{i=1}^n A_{i,n}(t)$ under this
representation and the rows are independent, the stationary
autocovariance of $S_n$ is the sum of the autocovariances
of the individual rows. Each row is an identification chain
of length $\ell_i$ with input mean $a_i$, so the stated
autocovariance formula follows from Theorem~\ref{thm:7.1}.
Setting $t=0$ gives the variance formula.

Finally, for every integer $\ell\ge0$, the elementary bounds
\[
    \sqrt{\ell+1}
    \le (2\ell+1)\frac{\binom{2\ell}{\ell}}{4^\ell}
    \le \sqrt{2\ell+1}
    \le \sqrt2\,\sqrt{\ell+1}
\]
hold, as can be verified by induction.
Moreover,
\[
    \frac23n^{3/2}
    \le\sum_{i=1}^n\sqrt{\ell_i+1}
    =\sqrt n+\sum_{k=2}^n\sqrt k \le n^{3/2}.
\]
The uniform bounds on the $a_i$ therefore imply
\[
    \frac23a_{\min}n^{3/2}
    \le \var_\pi(S_n)
    \le \sqrt2\,a_{\max}n^{3/2}.
\]
Also, $a_{\min}n \le \E_\pi[S_n]=\sum_{i=1}^n a_i \le a_{\max}n$. Consequently,
\[
    \frac{2a_{\min}}{3a_{\max}}\sqrt n
    \le \frac{\var_\pi(S_n)}{\E_\pi[S_n]}
    \le \frac{\sqrt2\,a_{\max}}{a_{\min}}\sqrt n.
\]
\end{proof}

All the preceding results can be summarized succinctly using standard asymptotic notation. Recall that $O(g)$ denotes an upper bound by a constant multiple of $g$, while $\Theta(g)$ denotes both upper and lower bounds by positive constant multiples of $g$, in the relevant asymptotic regime. Throughout, the implicit constants are independent of $n$ and $\varepsilon$, under the stated assumptions.

In this notation, Theorems~\ref{thm:6.5} and~\ref{thm:6.8} show that binary-tree addition has mean and autocovariance convergence-time upper bounds of $O(\log n+\log(1/\varepsilon))$, compared with $O(n+\log(1/\varepsilon))$ for linear chains. Furthermore, if the source means are uniformly bounded above and bounded away from zero, independently of $n$, Theorems~\ref{thm:binary-tree-addition-variance} and~\ref{thm:6.9} show that the stationary Fano factor grows as $\Theta(\sqrt{\log n})$ for binary trees, compared with $\Theta(\sqrt n)$ for linear chains. These comparisons are summarized in Table~\ref{table:addition-summary}.

\begin{table}[htbp]
\centering
\renewcommand{\arraystretch}{1.4} 
\begin{tabular}{|c|c|c|}
\hline
    & Binary tree 
    & Linear chain\\ 
\hline    
    Mean and autocovariance convergence times 
    & $O(\log n + \log \frac{1}{\varepsilon})$
    & $O(n + \log \frac{1}{\varepsilon})$\\ 
\hline
    Stationary Fano factor ($\var_\pi / \E_\pi$)
    & $\Theta(\sqrt{\log n})$
    & $\Theta(\sqrt{n})$\\ 
\hline    
\end{tabular}

\caption{Comparison of binary-tree and linear-chain
addition architectures. The stationary Fano-factor
asymptotics assume $0<a_{\min}\le a_i\le a_{\max}<\infty$
uniformly in $n$. The convergence-time bounds hold under
the uniform initial-moment assumptions of the corresponding
theorems, with an additional uniform upper bound on the
source means for autocovariance convergence. The convergence-time bounds hold uniformly in $n$ and
$0<\varepsilon\le1$, under the stated assumptions.}
\label{table:addition-summary}
\end{table}

\subsection{Multiplication}

In this subsection, for simplicity, we consider the computation of $a^n$ for a fixed $a>0$. Specifically, we assume that all source species in Figs.~\ref{graph: addition-binary-tree} and \ref{graph: addition-linear-chain} have the same stationary mean $a_j = a$ for $1 \le j \le n$. We compare the binary-tree and linear-chain architectures in terms of convergence time bounds and the stationary variance and autocovariance of their outputs. Throughout this subsection, we consider mutually independent source processes rather than repeated use of a single random input.

\subsubsection{Binary-tree architecture}

We consider a balanced binary-tree multiplication architecture of depth $m$, with $n=2^m$ input species, each having stationary mean $a$. At level $i$, each species $A_{i,j}$ has stationary mean $a^{2^i}$, so the final output $A_{m,1}$ has stationary mean $a^n$. We first derive bounds on the mean and autocovariance convergence times of the final output, and then study its stationary variance and autocovariance.

\begin{thm}[Mean and autocovariance convergence times for binary-tree multiplication]
\label{lem:6.11}
Fix $a>0$. For each $n=2^m$, $m\ge1$, consider the balanced binary-tree multiplication architecture in Fig.~\ref{graph: addition-binary-tree}. Assume that the initial copy numbers of all distinct species in the architecture are mutually independent.

Suppose there exists a constant $D<\infty$, independent of $n$,
such that
\[
    \max_{0\le i\le m}\max_{1\le j\le 2^{m-i}}
    \left|\E[A_{i,j}(0)]-a^{2^i}\right|\le D.
\]
Then there exists a constant $C_{\mathrm{mean}}<\infty$,
independent of $n$ and $\varepsilon$, such that
\[
    t_{\mathrm{mean}}(\varepsilon)
    \le C_{\mathrm{mean}} \cdot
    \left(n+\log\frac1\varepsilon\right)
\]
for all $n=2^m\ge2$ and $0<\varepsilon\le1$.

Additionally, suppose there exists a constant $V<\infty$,
independent of $n$, such that
\[
    \max_{0\le i\le m}\max_{1\le j\le 2^{m-i}}
    \var(A_{i,j}(0))\le V.
\]
Then there exists a constant $C_{\mathrm{cov}}<\infty$,
independent of $n$ and $\varepsilon$, such that
\[
    t_{\mathrm{cov}}(\varepsilon) 
    \le C_{\mathrm{cov}} \cdot
    \left(n+\log\frac1\varepsilon\right)
\]
for all $n=2^m\ge2$ and $0<\varepsilon\le1$.
\end{thm}

\begin{proof}
See Appendix~\ref{sec:appendix-c} for the proof. The proof uses a recursion for the mean error and establishes simultaneous bounds on second moments and autocovariance errors, with constants growing at most exponentially in $n$.  
\end{proof}

To analyze stationary autocovariances for the multiplication architectures, we first establish bounds relating them to stationary variances. These bounds show that, for each fixed time $t$, the autocovariance and variance have the same asymptotic order as the number of inputs grows. Thus, it suffices to study stationary variances to obtain the autocovariance asymptotics below.

For the binary-tree architecture, we write $A_i$ for a representative species at level $i$, since all species at that level have the same stationary process law. When referring to two sibling inputs, we retain the distinct labels $A_{i,1}$ and $A_{i,2}$.

\begin{lem}[Comparison of stationary autocovariance and variance]\label{lem:6.12}
    For the species at each level of the binary-tree multiplication architecture in Fig.~\ref{graph: addition-binary-tree},
    \begin{align*}
        e^{-t} \var_{\pi} (A_i ) \leq \cov_{\pi} (A_i (0), A_i (t)) \le \var_\pi (A_i  ),
    \end{align*}
    for all $0 \le i \le m$ and $t \ge 0$. The same upper and lower bounds hold for the linear-chain architecture in Fig.~\ref{graph: addition-linear-chain}, with $S_i$ in place of $A_i$ and for $1 \le i \le n$.
\end{lem}

\begin{proof}
The upper bound follows from the Cauchy--Schwarz inequality and stationarity:
\begin{align*}
    \cov_\pi (A_i (0), A_i (t)) \leq \sqrt{\var_\pi (A_i (0)) \var_\pi (A_i (t))} = \var_\pi (A_i).
\end{align*}
We prove the lower bound by induction on $i$. For $i=0$, 
    \begin{align*}
        e^{-t} \var_{\pi} (A_0 ) = a e^{-t} = \cov_{\pi} (A_0 (0) , A_0 (t)).
    \end{align*}
    Suppose that the bound holds at level $i \ge 0$. By Corollary \ref{lem : cox-process-variance},
    \begin{align*}
        \cov_{\pi} (A_{i+1}(0), A_{i+1}(t)) = e^{-t} \E_\pi[\lambda(Z)] + \frac{1}{2} \int_0 ^{\infty} e^{-u} (C_\lambda (t+u) + C_\lambda(|t-u|)) du,
    \end{align*}
    and
    \begin{align*}
        \var_\pi (A_{i+1}) = \E_\pi [\lambda(Z)] + \int_0 ^{\infty} e^{-u} C_\lambda(u) du.
    \end{align*}
    Therefore,
    \begin{align*}
        \cov_{\pi} (A_{i+1}(0), A_{i+1}(t)) - e^{-t} \var_\pi (A_{i+1}) = \frac{e^{-t}}{2} \int_0 ^{t} (e^{u} - e^{-u}) C_\lambda(u) du + \frac{e^{-t}}{2} (e^{2t} - 1) \int_t ^{\infty} e^{-u} C_\lambda(u) du.
    \end{align*}
    For the binary tree multiplication module,
    \begin{align*}
    \begin{split}
        C_\lambda(s) &= \cov_\pi  (\lambda(Z(0)), \lambda(Z(s)))\\ &= \cov_\pi (A_{i,1} (0) A_{i,2}(0), A_{i,1}(s)A_{i,2}(s)) \\
        &= \cov_\pi (A_i (0), A_i (s)) ^2 + 2 \E_\pi [A_i ]^2 \cov_\pi (A_i (0), A_i (s))  \geq 0,
    \end{split}
    \end{align*}
    since $\cov_\pi (A_i (0), A_i (s)) \ge e^{-s} \var_\pi (A_i) \ge 0$ by induction. Therefore,
    \begin{align*}
        \cov_{\pi} (A_{i+1}(0), A_{i+1}(t)) - e^{-t} \var_{\pi} (A_{i+1}) \ge 0.
    \end{align*}
    The upper bound also holds for the linear-chain architecture. For the lower bound, we proceed by induction on $i$. The base case is $S_1 = A_1$, for which 
    \[
        \cov_\pi (S_1 (0), S_1 (t)) = a e^{-t} = e^{-t} \var_\pi (S_1).
    \]
    Suppose the lower bound holds for $S_i$, where $1 \le i < n$. For the module producing $S_{i+1}$, we have
    \begin{align*}
    \begin{split}
        C_\lambda(s) &= \cov _\pi (\lambda(Z(0)), \lambda(Z(s)))\\ &= \cov_\pi (S_{i} (0) A_{i+1}(0), S_{i}(s)A_{i+1}(s)) \\
        &= a e^{-s} \cov_\pi (S_i (0), S_i (s)) + a^2 \cov_\pi (S_i (0), S_i (s)) + a (\E_\pi [S_i]) ^2 e^{-s} \geq 0,
    \end{split}
    \end{align*}
    so we get
    \begin{align*}
          \cov_\pi (S_{i+1} (0), S_{i+1} (t)) - e^{-t} \var_\pi (S_{i+1} ) \ge 0.
    \end{align*}
\end{proof}

We now use the preceding comparison between stationary variance and autocovariance to derive asymptotic estimates for the output of the binary-tree multiplication architecture.

\begin{thm}[Stationary autocovariance asymptotics for binary-tree multiplication]
Consider the balanced binary-tree multiplication structure in Fig.~\ref{graph: addition-binary-tree} with $n=2^{m}$ mutually independent source processes, each having stationary mean $a >0$. Let $A_{m,1}$ be the final output species and define
\begin{align*}
    C_m(t):=\cov_\pi(A_{m,1}(0),A_{m,1}(t)).
\end{align*}
Then, for each fixed $t\geq 0$, as $n=2^{m}\to\infty$, there exist constants $0< c_1 \le c_2$, depending only on $a$ and $t$, such that
\begin{alignat*}{2}
a^{2n-1}e^{c_1 n}
&\le C_m(t)\le a^{2n-1}e^{c_2 n},
&\qquad& a>1,\\
e^{c_1 n}
&\le C_m(t)\le e^{c_2 n},
&\qquad& a=1,\\
a^ne^{c_1 n}
&\le C_m(t)\le a^ne^{c_2 n},
&\qquad& 0<a<1.
\end{alignat*}
\label{thm:multiplication-binary-tree}
\end{thm}

\begin{proof}
See Appendix~\ref{sec:appendix-a-2} for the proof. The argument combines Lemmas~\ref{lem:6.12} and \ref{lem:5.4} to derive recurrence bounds for the stationary variance, which then yield the autocovariance asymptotics.
\end{proof}

\begin{rem}
    Theorem~\ref{thm:multiplication-binary-tree} shows that the stationary variance and autocovariance of binary-tree multiplication exhibit exponential amplification relative to the output mean when $0<a<1$. This does not necessarily imply growth of the absolute variance or autocovariance, since the output mean $a^n$ decreases exponentially in this regime. For $a\ge1$, the theorem gives exponential growth of the absolute stationary variance and autocovariance. In contrast, for binary-tree addition with identical input means, the stationary variance grows as $\Theta(n\sqrt{\log n})$, which is subexponential in $n$. 
    
    This difference arises from how input covariances propagate through the modules. For independent stationary input processes $A$ and $B$, the addition module has 
    \[
    C_{A+B}(s)=C_A(s)+C_B(s),
\]
    whereas the product input to the multiplication module satisfies
\[
    C_{AB}(s) =C_A(s)C_B(s) +\E_\pi[A]^2 C_B(s) +\E_\pi[B]^2 C_A(s).
\]
Thus, the multiplication module involves both a nonlinear covariance term and terms weighted by the squared stationary input means. These terms can produce exponential amplification of fluctuations relative to the output mean.
\end{rem}

\subsubsection{Linear-chain architecture}

We now consider the linear-chain multiplication architecture shown in Fig.~\ref{graph: addition-linear-chain}. Let $S_i$ denote the output species after multiplying the first $i$ inputs, with $S_1 = A_1$. The stationary mean of $S_i$ is $a^i$, so the final output $S_n$ has stationary mean $a^n$. We first derive bounds on the mean and autocovariance convergence times, and then study the stationary autocovariance bounds for comparison with the binary-tree architecture.

\begin{thm}[Mean and autocovariance convergence times for linear-chain multiplication]
\label{thm:mean-convergence-of-linear-chain-multiplication}
Fix $a>0$. For each $n\ge2$, consider the linear-chain multiplication architecture in Fig.~\ref{graph: addition-linear-chain} computing $a^n$. Let $S_i$ denote the output species after multiplying the first $i$ inputs, with $S_1=A_1$. Assume that the initial copy numbers of all distinct species in the architecture are mutually independent.

Suppose there exist constants $D,C<\infty$, independent of $n$, such that
\[
    \max_{1\le i\le n}|\E[S_i(0)]-a^i|\le D,
    \qquad
    \max_{1\le i\le n}|\E[A_i(0)]-a|\le C.
\]
Then there exists a constant $C_{\mathrm{mean}}<\infty$, independent of $n$ and $\varepsilon$, such that
\[
    t_{\mathrm{mean}}(\varepsilon)
    \le C_{\mathrm{mean}} \cdot
    \left(n+\log\frac1\varepsilon\right)
\]
for all $n\ge2$ and $0<\varepsilon\le1$. Additionally, suppose there exists a constant $V<\infty$, independent of $n$, such that
\[
    \max_{1\le i\le n}
    \max\{\var(S_i(0)),\var(A_i(0))\}\le V.
\]
Then there exists a constant $C_{\mathrm{cov}}<\infty$, independent of $n$ and $\varepsilon$, such that
\[
    t_{\mathrm{cov}}(\varepsilon)
    \le C_{\mathrm{cov}} \cdot
    \left(n+\log\frac1\varepsilon\right)
\]
for all $n\ge2$ and $0<\varepsilon\le1$.
\end{thm}

\begin{proof}
    See Appendix \ref{sec:appendix-d} for the proof. The argument bounds the mean-error recursion and propagates second-moment and autocovariance-error bounds, with constants growing at most exponentially in $n$.
\end{proof}

We next turn to the stationary variance and autocovariance of the linear-chain output. Unlike binary-tree multiplication, where the stationary variance recursion involves a quadratic term in the variance of the preceding level, the corresponding recursion for linear-chain multiplication is linear. This leads to different bounds depending on the value of $a$. As in the binary-tree case, Lemma~\ref{lem:6.12} allows us to obtain autocovariance bounds from stationary variance bounds.

\begin{thm}[Stationary autocovariance bounds for linear-chain multiplication]
\label{thm:fixed-time-aucotovariance-bounds-for-linear-chain-multiplication}
Consider the linear-chain multiplication structure in
Fig.~\ref{graph: addition-linear-chain} computing $a^n$.
Let $S_n$ be the output species after multiplying the first
$n$ inputs, and define
$C_n(t):=\cov_\pi(S_n(0),S_n(t))$.

Then, for each fixed $t\geq 0$, as $n\to\infty$, there exist constants $c_1 , c_2 > 0$, depending only on $a$ and $t$, for which the following bounds for $C_n(t)$ hold:
\begin{alignat*}{2}
c_1a^{2n-1}
&\le C_n(t)\le c_2a^{2n-1}\left(1+\frac{1}{2a}\right)^n,
&\qquad& a>1,\\
c_1
&\le C_n(t)\le c_2\left(\frac32\right)^n,
&\qquad& a=1,\\
c_1a^n
&\le C_n(t)\le c_2a^n\left(a+\frac12\right)^n,
&\qquad& \frac12<a<1,\\
c_1 2^{-n}
&\le C_n(t)\le c_2 n2^{-n},
&\qquad& a=\frac12,\\
c_1a^n
&\le C_n(t)\le c_2a^n,
&\qquad& 0<a<\frac12.
\end{alignat*}
In particular, when $0<a<\frac12$, the bounds match,
giving $C_n(t)=\Theta (a^n)$.
\end{thm}

\begin{proof}
See Appendix~\ref{sec:appendix-e} for the proof. The argument combines Lemmas~\ref{lem:6.12} and \ref{lem:5.4} to derive recursive upper and lower bounds for the stationary variance, which then yield the autocovariance bounds.
\end{proof}

The most notable threshold in Theorem~ \ref{thm:fixed-time-aucotovariance-bounds-for-linear-chain-multiplication} occurs at $a=1/2$. When $a<1/2$, the stationary variance of the linear-chain output is of the same order as its stationary mean, $a^n$. At $a=1/2$, our upper and lower bounds differ by at most a factor of order $n$. For larger values of $a$, we obtain upper and lower bounds that do not generally match. The distinction between $a<1$ and $a>1$ also reflects the different behavior of the output mean $a^n$ in these two regimes.

Table~\ref{table:multiplication-summary} compares the stationary variance results for the two architectures, together with their convergence-time bounds. The variance results follow by taking $t=0$ in the preceding autocovariance theorems.

\begin{table}[h!]
\centering
\renewcommand{\arraystretch}{1.2}
\begin{tabular}{llcc}
\hline
\textbf{Quantity} & \textbf{Regime} & \textbf{Binary tree} & \textbf{Linear-chain bounds} \\
\hline
\multirow{5}{*}{Variance}
& $a>1$
& $a^{2n-1}\exp(\Theta(n))$
& $O\!\left(a^{2n-1}\left(1+\frac{1}{2a}\right)^n\right)$ \\

& $a=1$
& $\exp(\Theta(n))$
& $O\!\left(\left(\frac{3}{2}\right)^n\right)$ \\

& $\frac12<a<1$
& $a^n\exp(\Theta (n))$
& $O\!\left(a^n\left(a+\frac{1}{2}\right)^n\right)$ \\

& $a=\frac12$
& $2^{-n}\exp(\Theta (n))$
& $O\!\left(n 2^{-n}\right)$ \\

& $0<a<\frac12$
& $a^n\exp(\Theta(n))$
& $\Theta\!\left(a^n\right)$ \\
\hline
$ t_{\mathrm{mean}}(\varepsilon),t_{\mathrm{cov}}(\varepsilon)$
& 
& $O\!\left(n  + \log \frac{1}{\varepsilon}\right)$
& $O\!\left(n  + \log \frac{1}{\varepsilon}\right)$ \\
\hline
\end{tabular}
\caption{Summary of variance asymptotics and bounds for multiplication architectures}
\label{table:multiplication-summary}
\end{table}

The comparison reveals an important difference from addition. For $0<a\le 1/2$, the linear-chain architecture has asymptotically smaller stationary variance than the binary-tree architecture. This follows from the theoretical bounds: binary-tree multiplication exhibits exponential amplification relative to the output mean $a^n$, whereas the linear-chain variance is of order $a^n$ when $a<1/2$ and is at most of order $n2^{-n}$ when $a=1/2$. For $a>1/2$, the theoretical bounds do not establish a general ordering of the stationary variances. However, the simulations in Table~\ref{table:multiplication-simul} show larger stationary variances for binary-tree multiplication at each of the tested values of $a$, with $n=8$.

\begin{table}[h!]
\centering
\renewcommand{\arraystretch}{1.2}
\begin{tabular}{|c|c|c|c|c|c|c|}
\hline
$a$ & 1.0 & 1.2 & 1.4 & 1.6 & 1.8 & 2.0\\
\hline
Binary tree & 34.03 & 303.81 & 2408.34 & 12391.87 & 61136.38 & 241137.19\\
\hline
Linear chain & 20.19 & 187.59 & 1549.77 &  8757.58 & 40825.62 & 197498.47 \\
\hline
\end{tabular}
\caption{Estimated stationary variances for the computation of $a^8$ at different values of $a$, based on 10,000 independent stochastic simulation algorithm samples collected after burn-in.}
\label{table:multiplication-simul}
\end{table}

The convergence-time bounds provide another contrast with addition. Although the binary-tree multiplication architecture has logarithmic depth, both multiplication architectures have mean and autocovariance convergence-time upper bounds of order $n+\log(1/\varepsilon)$. Thus, our analysis does not establish a convergence-time advantage for binary-tree multiplication. For addition, by contrast, the binary-tree architecture gives better asymptotic bounds for both convergence times and stationary variance.

Taken together, these results suggest that the relative advantages of the architectures depend strongly on the operation being performed.

\section{Discussion}
\label{sec:discussion}

In this paper, we studied the computation of identification, addition, and multiplication using stochastic reaction networks, with values encoded by stationary means. A central feature of the stochastic setting is that establishing correctness and convergence of the output mean is not sufficient to characterize readout performance. One must also understand the variance and autocovariance of the output, and it is this analysis that reveals some of the more surprising effects of network architecture. For addition, the binary-tree architecture provides better asymptotic upper bounds on convergence times and a smaller stationary Fano factor than the linear-chain architecture, under the stated assumptions. For multiplication, however, the binary-tree architecture is not uniformly advantageous. In particular, multiplication can produce exponential amplification of stationary fluctuations relative to the output mean, despite the smaller depth of the binary tree. 
This contrast, particularly for multiplication, was unexpected and highlights the need for careful analysis of variance and autocovariance when evaluating stochastic computational architectures. In particular, architectural depth alone is not a reliable indicator of readout performance.

This work provides a starting point for the broader study of arithmetic computation using stochastic reaction networks. 
Our analysis of composite architectures has focused on
repeated applications of a single operation. An important
next step is to study computations that combine different
arithmetic operations, and to understand how fluctuations
propagate through such networks.
  Natural directions for future research also include constructing modules for additional arithmetic operations and extending the framework from positive numbers to signed real values, perhaps through dual-rail representations \cite{chen2023rate, chen2014deterministic}. An important challenge in these extensions will be not only to establish correctness at the stationary mean level, but also to understand how the resulting networks' convergence and fluctuation properties affect the cost of reading the computed values.

\section*{Acknowledgments}

DFA gratefully acknowledges financial support from the Office of the Vice Chancellor for Research at UW--Madison, with funding from the Wisconsin Alumni Research Foundation; from the Trustees of the William F. Vilas Estate; and from NSF grant DMS-2051498.

The authors used OpenAI’s ChatGPT as an interactive assistant for bounded tasks involving exposition, algebraic manipulation, and formalization. The authors formulated each task, examined the resulting output, and independently verified every calculation and argument before incorporating any material into the manuscript. The mathematical ideas, analysis, and conclusions, including any typos or errors, are entirely the authors' own.

\bibliographystyle{unsrt}
\bibliography{ref}

\clearpage
\newpage
\appendix
\appendixpage
\addtocontents{toc}{\protect\setcounter{tocdepth}{1}}

\section{Proof of Lemma~\ref{lem:general-convergence-rate}}
\label{sec:ProofOfLemmaTransfer}

\begin{proof}    
By Lemma~\ref{lem:rate-of-convergence} and the local boundedness assumption, there exists $c>0$ such that
\[
    |\E[\lambda(Z(t))]-\bar\lambda| \le ce^{-\rho t}, \quad t\ge0.
\]
 From Corollary~\ref{cor:finite-time-moments} and using that $\E_\pi [X] = \E_\pi[\lambda(Z)] / \gamma$,
    \begin{align*}
        \E [X(t)] - \E_\pi [X] &= e^{-\gamma t} \E[X(0)] + \int_0 ^t e^{-\gamma(t-s)} \E[\lambda(Z(s))] ds - \E_\pi [X]\\
        &= e^{- \gamma t} (\E[X(0)] - \E_\pi [X]) + \int_0 ^t e^{-\gamma(t-s)} (\E[\lambda(Z(s))] -\E_\pi [\lambda(Z)])ds.
        \end{align*}
Therefore,
\begin{align}\label{eq:mean-convergence-inequality}
\begin{split}
    |\E[X(t)] - \E_\pi [X]| &\le e^{-\gamma  t } |\E[X(0)] - \E_\pi [X]| + c \int_0 ^t e^{- \gamma (t-s)} e^{- \rho s} ds\\
    &=\begin{cases}
        |\E[X(0)] - \E_\pi [X]| e^{- \gamma t} + \frac{c}{|\gamma - \rho|} |e^{- \rho t} - e^{- \gamma t}|, &\text{if }\rho \neq \gamma,\\
        |\E[X(0)] - \E_\pi [X]| e^{- \gamma t} + c t e^{- \gamma t}, &\text{if }\rho = \gamma,
    \end{cases}\\
    & \le |\E[X(0)] - \E_\pi [X]| e^{- \gamma t}+ c t e^{- \min\{\rho, \gamma\} t}\\
    & \le (|\E[X(0)] - \E_\pi [X]| + c) (1+t)e^{- \min\{\rho, \gamma\} t} .
\end{split}
\end{align}
Therefore, by Lemma~\ref{lem:rate-of-convergence}, $\E [X(t)]$ converges at rate at least $\min \{\rho, \gamma\}$. 

We next assume the additional hypotheses in the second part of the lemma and establish uniform convergence of the autocovariance. The local second-moment bound, together with the assumed convergence of the input mean and autocovariance, implies
\[
    \sup_{t\ge0}\E[\lambda(Z(t))^2]<\infty.
\]
By Cauchy--Schwarz, this yields the local uniform boundedness required by Lemma~\ref{lem:covariance-rate-of-convergence}. Consequently, for some $c_\lambda>0$,
\[
\sup_{\tau\ge0}
\left|
\cov(\lambda(Z(t)),\lambda(Z(t+\tau)))
-C_\lambda(\tau)
\right|
\le c_\lambda e^{-\sigma t},
\qquad t\ge0.
\]
Write $C_\lambda(\tau) := \cov_\pi(\lambda(Z(0)),\lambda(Z(\tau)))$, $C_X(\tau) := \cov_\pi(X(0),X(\tau))$. By Corollaries~\ref{cor:finite-time-moments} and \ref{lem : cox-process-variance},
\begin{align*}
    \cov(X(t),X(t+\tau)) &= e^{-\gamma\tau} \left( e^{-\gamma t}(1-e^{-\gamma t})\E[X(0)] +\E[M_t^{(0)}] \right)+ \cov\left( e^{-\gamma t}X(0)+M_t^{(0)},e^{-\gamma(t+\tau)}X(0)+M_{t+\tau}^{(0)}\right),\\
    C_X(\tau)&=e^{-\gamma\tau}\E_\pi[X]+\int_{-\infty}^t\int_{-\infty}^{t+\tau}e^{-\gamma(t-u)}e^{-\gamma(t+\tau-v)}C_\lambda(|v-u|)\,dv\,du.
\end{align*}
Since $\E[X(t)]=e^{-\gamma t}\E[X(0)]+\E[M_t^{(0)}]$, the first term in the finite-time covariance formula satisfies
\begin{align*}
    e^{-\gamma\tau}\left(e^{-\gamma t}(1-e^{-\gamma t})\E[X(0)]+\E[M_t^{(0)}]\right)=e^{-\gamma\tau}\E[X(t)]-e^{-\gamma(2t+\tau)}\E[X(0)].
\end{align*}
Moreover,
\begin{align*}
    \cov\left(
        e^{-\gamma t}X(0),
        e^{-\gamma(t+\tau)}X(0)
    \right)
    =
    e^{-\gamma(2t+\tau)}\var(X(0)).
\end{align*}
Expanding the remaining covariance terms using the definition of $M_t ^{(0)}$, and splitting the stationary double integral into the common region $[0,t]\times[0,t+\tau]$ and its complement, the covariance error can be expressed as
\begin{align}\label{eq:covariance-convergence-upper-bound}
    \begin{split}
    &|\cov(X(t),X(t+\tau))-C_X(\tau)|\\
    &=
    \left|
    \begin{aligned}
        &e^{-\gamma\tau}
        \bigl(\E[X(t)]-\E_\pi[X]\bigr)\\
        &+e^{-\gamma(2t+\tau)}
        \bigl(\var(X(0))-\E[X(0)]\bigr)\\
        &+e^{-\gamma t}
        \int_0^{t+\tau}
        e^{-\gamma(t+\tau-v)}
        \cov(X(0),\lambda(Z(v)))\,dv\\
        &+e^{-\gamma(t+\tau)}
        \int_0^t
        e^{-\gamma(t-u)}
        \cov(X(0),\lambda(Z(u)))\,du\\
        &+\int_0^t\int_0^{t+\tau}
        e^{-\gamma(t-u)}e^{-\gamma(t+\tau-v)}
        \Bigl[
            \cov(\lambda(Z(u)),\lambda(Z(v)))
            -C_\lambda(|v-u|)
        \Bigr]\,dv\,du\\
        &-\int_{-\infty}^{0}\int_{-\infty}^{t+\tau}
        e^{-\gamma(t-u)}e^{-\gamma(t+\tau-v)}
        C_\lambda(|v-u|)\,dv\,du\\
        &-\int_0^t\int_{-\infty}^{0}
        e^{-\gamma(t-u)}e^{-\gamma(t+\tau-v)}
        C_\lambda(|v-u|)\,dv\,du
    \end{aligned}
    \right|.
     \end{split}
\end{align}
We estimate the terms on the right-hand side separately. In particular, for every $u,v\ge0$,
\begin{align}\label{eq:temp_proof}
    |\cov(\lambda(Z(u)),\lambda(Z(v))) -C_\lambda(|v-u|)| \le c_\lambda e^{-\sigma\min\{u,v\}} \quad\text{for all }u,v\ge0.
\end{align}
For the first term, \eqref{eq:mean-convergence-inequality} gives
\begin{align*}
    e^{-\gamma\tau}|\E[X(t)]-\E_\pi[X]|\le|\E[X(0)]-\E_\pi[X]|e^{-\gamma t}+ct e^{-\min\{\rho,\gamma\}t}.
\end{align*}
The second term satisfies
\begin{align*}
    e^{-\gamma(2t+\tau)}|\var(X(0))-\E[X(0)]|\le|\var(X(0))-\E[X(0)]|e^{-2\gamma t}.
\end{align*}
For the third and fourth terms, taking $u=v=s$ in \eqref{eq:temp_proof} yields
\begin{align*}
    \var(\lambda(Z(s)))\le\var_\pi (\lambda(Z)) + c_\lambda e^{-\sigma s} \le \var_\pi(\lambda(Z))+c_\lambda\quad\text{for all }s\ge0.
\end{align*}
Thus, by the Cauchy--Schwarz inequality,
\begin{align*}
    |\cov(X(0),\lambda(Z(s)))|\le \sqrt{\var(X(0)) \var(\lambda(Z(s)))} \le \sqrt{\var(X(0))\bigl(\var_\pi(\lambda(Z))+c_\lambda\bigr)}.
\end{align*}
Consequently,
\begin{align*}
    \left|e^{-\gamma t}\int_0^{t+\tau}
    e^{-\gamma(t+\tau-v)}
    \cov(X(0),\lambda(Z(v)))\,dv
    \right|\le
    \frac{e^{-\gamma t}}{\gamma}
    \sqrt{
        \var(X(0))
        \bigl(\var_\pi(\lambda(Z))+c_\lambda\bigr)
    },
\end{align*}
and
\begin{align*}
    \left|
    e^{-\gamma(t+\tau)}
    \int_0^t
    e^{-\gamma(t-u)}
    \cov(X(0),\lambda(Z(u)))\,du
    \right|\le
    \frac{e^{-\gamma t}}{\gamma}
    \sqrt{
        \var(X(0))
        \bigl(\var_\pi(\lambda(Z))+c_\lambda\bigr)
    }.
\end{align*}
To bound the fifth term, we use the assumed covariance estimate and split the integration domain into two regions;
\begin{align*}
    &\left|
    \int_0^t\int_0^{t+\tau}
    e^{-\gamma(t-u)}e^{-\gamma(t+\tau-v)}
    \Bigl[
        \cov(\lambda(Z(u)),\lambda(Z(v)))
        -C_\lambda(|v-u|)
    \Bigr]\,dv\,du
    \right|\\
    &\quad\le
    c_\lambda
    \int_0^t\int_0^{t+\tau}
    e^{-\gamma(t-u)}e^{-\gamma(t+\tau-v)}
    e^{-\sigma\min\{u,v\}}\,dv\,du\\
    &\quad=
    c_\lambda
    \int_0^t
    e^{-\gamma(t-u)}e^{-\sigma u}
    \left(
        \int_u^{t+\tau}e^{-\gamma(t+\tau-v)}\,dv
    \right)du+
    c_\lambda
    \int_0^t
    e^{-\gamma(t+\tau-v)}e^{-\sigma v}
    \left(
        \int_v^t e^{-\gamma(t-u)}\,du
    \right)dv\\
    &\quad\le
    \frac{c_\lambda}{\gamma}(1+e^{-\gamma\tau})
    \int_0^t e^{-\gamma(t-s)}e^{-\sigma s}\,ds \le
    \frac{2c_\lambda}{\gamma}\,
    t e^{-\min\{\sigma,\gamma\}t}.
\end{align*}
The last inequality follows from $\gamma(t-s)+\sigma s \ge \min\{\sigma,\gamma\}t$ for $0\le s\le t$. For the sixth and seventh terms, stationarity and the Cauchy--Schwarz inequality imply $|C_\lambda(h)|\le \var_\pi(\lambda(Z))$ for all $h\ge0$. Therefore,
\begin{align*}
    \left|
    \int_{-\infty}^{0}\int_{-\infty}^{t+\tau}
    e^{-\gamma(t-u)}e^{-\gamma(t+\tau-v)}
    C_\lambda(|v-u|)\,dv\,du
    \right| &\le
    \var_\pi(\lambda(Z))
    \left(\int_{-\infty}^{0}e^{-\gamma(t-u)}\,du\right)
    \left(
        \int_{-\infty}^{t+\tau}
        e^{-\gamma(t+\tau-v)}\,dv
    \right)\\
    &=
    \frac{\var_\pi(\lambda(Z))}{\gamma^2}
    e^{-\gamma t},
\end{align*}
and
\begin{align*}
    \left|
    \int_0^t\int_{-\infty}^{0}
    e^{-\gamma(t-u)}e^{-\gamma(t+\tau-v)}
    C_\lambda(|v-u|)\,dv\,du
    \right|&\le
    \var_\pi(\lambda(Z))
    \left(\int_0^t e^{-\gamma(t-u)}\,du\right)
    \left(
        \int_{-\infty}^{0}
        e^{-\gamma(t+\tau-v)}\,dv
    \right)\\
    &\quad=
    \frac{\var_\pi(\lambda(Z))}{\gamma^2}
    (1-e^{-\gamma t})e^{-\gamma(t+\tau)}\le
    \frac{\var_\pi(\lambda(Z))}{\gamma^2}
    e^{-\gamma t}.
\end{align*}
Combining these estimates using the triangle inequality yields
\begin{align*}
    |\cov(X(t),X(t+\tau))-C_X(\tau)| & \le
    |\E[X(0)]-\E_\pi[X]|e^{-\gamma t}
    +ct e^{-\min\{\rho,\gamma\}t}+
    |\var(X(0))-\E[X(0)]|e^{-2\gamma t}\\
    &\hspace{.2in} +\frac{2}{\gamma}
    \sqrt{
        \var(X(0))
        \bigl(\var_\pi(\lambda(Z))+c_\lambda\bigr)
    }\,
    e^{-\gamma t} +
    \frac{2c_\lambda}{\gamma}\,
    t e^{-\min\{\sigma,\gamma\}t}
    +\frac{2\var_\pi(\lambda(Z))}{\gamma^2}
    e^{-\gamma t}.
\end{align*}
Since the right-hand side does not depend on $\tau$, we may take the supremum over $\tau \ge 0$. Hence, setting
$r=\min\{\rho,\sigma,\gamma\}$, we obtain
\begin{align}\label{eq:cov_upper_bound_polynomial}
    \sup_{\tau\ge0}
    |\cov(X(t),X(t+\tau))-C_X(\tau)|
    \le D(1+t)e^{-rt},
\end{align}
where we may take
\begin{align*}
    D=|\E[X(0)]-\E_\pi[X]|+c
    +|\var(X(0))-\E[X(0)]|+
    \frac{2}{\gamma}
    \sqrt{
        \var(X(0))
        \bigl(\var_\pi(\lambda(Z))+c_\lambda\bigr)
    }+\frac{2c_\lambda}{\gamma}
    +\frac{2\var_\pi(\lambda(Z))}{\gamma^2}.
\end{align*}
This constant is finite by the moment assumptions. Therefore, by Lemma~\ref{lem:covariance-rate-of-convergence}, the autocovariance of $X$ converges uniformly to $C_X$ at rate at least $\min\{\rho,\sigma,\gamma\}$.
\end{proof}

\section{Proof of Theorem~\ref{thm:6.8}}
\label{sec:proof-thm-6.8}

\begin{proof}
We work with the auxiliary decomposition network initialized according to the distribution specified in the theorem. By Theorem~\ref{thm: addition-linear-translation}, its projection has the same law as the original linear-chain addition process under the corresponding initial distribution. The auxiliary network consists of disjoint identification chains, whose initial populations need not be independent.

Let $s:=\sum_{i=1}^n a_i$ and define
\[
g_\ell(t):=e^{-t}\sum_{p=0}^{\ell}\frac{t^p}{p!}
=\mathbb P(\operatorname{Poisson}(t)\le\ell).
\]
For $t\ge0$, we have
\[
g_\ell(t)\le1,
\qquad
g_\ell(t)\le2^\ell e^{-t/2},
\]
where the second inequality follows from Markov's inequality
applied to $2^{-\operatorname{Poisson}(t)}$.

Since $S_n(t)=\sum_{i=1}^n A_{i,n}(t)$, the mean estimate
from the proof of Theorem~\ref{thm:6.1} gives
\[
\begin{aligned}
|\E[S_n(t)]-s|
&\le D\sum_{i=1}^n g_{\ell_i}(t) \le nD\,2^{n-1}e^{-t/2}.
\end{aligned}
\]
Consequently,
\[
t_{\mathrm{mean}}(\varepsilon)
\le 2(n-1)\log2
+2\log\frac{n\max\{1,D\}}{\varepsilon}.
\]
This establishes the claimed mean convergence-time bound.

For the covariance, Theorem~\ref{thm:6.1} gives, for each row $i$,
\[
\sup_{\tau\ge0}
\left|
\cov(A_{i,n}(t),A_{i,n}(t+\tau))
-C_{A_{i,n}}(\tau)
\right|
\le 4(1+a_i+D+V)6^{\ell_i}e^{-t/2}.
\]
At stationarity, the different rows are independent, so
\[
C_{S_n}(\tau)=\sum_{i=1}^n C_{A_{i,n}}(\tau).
\]
However, cross-row covariances may be nonzero at finite times because the initial populations need not be independent. Let $\mathcal F$ be the sigma-field generated by all initial populations in the auxiliary network, and write
\[
M_i(u):=\E[A_{i,n}(u)\mid\mathcal F].
\]
Since distinct rows are driven by independent reaction randomness, they are conditionally independent given $\mathcal F$. Thus, for $i\neq j$, the law of total covariance gives
\[
\cov(A_{i,n}(t),A_{j,n}(t+\tau))
=\cov(M_i(t),M_j(t+\tau)).
\]
The conditional mean formula for an identification chain yields
\[
M_i(u)
=a_i+e^{-u}\sum_{p=0}^{\ell_i}
\frac{u^p}{p!}(A_{i,n-p}(0)-a_i).
\]
By Cauchy--Schwarz and the uniform variance bound,
\begin{align*}
    |\cov(M_i(t),M_j(t+\tau))| \le Vg_{\ell_i}(t)g_{\ell_j}(t+\tau) \le V2^{\ell_i}e^{-t/2} \le V2^{n-1}e^{-t/2}.
\end{align*}
Combining the within-row and cross-row bounds, we obtain
\[
\begin{aligned}
\sup_{\tau\ge0} |\cov(S_n(t),S_n(t+\tau))-C_{S_n}(\tau)| &\le \left\{ 4n(1+A+D+V)6^{n-1} +n(n-1)V2^{n-1} \right\}e^{-t/2}\\
&\quad\le \bigl(4(1+A+D+V)+V\bigr) n^2 6^{n-1}e^{-t/2}.
\end{aligned}
\]
Therefore,
\begin{align*}
    t_{\mathrm{cov}}(\varepsilon) \le 2(n-1)\log6+4\log n+2\log \frac{4(1+A+D+V)+V}{\varepsilon}.
\end{align*}
Since $\log n\le n$, this establishes the claimed bound with a constant $C_{\mathrm{cov}}$ depending only on $A,D,V$.
\end{proof}

\section{Proof of Theorem~\ref{lem:6.11}}
\label{sec:appendix-c}

\begin{proof}
\smallskip
\noindent
\textbf{Mean convergence.} Define $a_i := \E_\pi[A_{i,j}] = a^{2^i}$. For each level $i$, define the mean discrepancy
\begin{align*}
    d_i(t)
    := \sup_j \left|\E[A_{i,j}(t)]-a_i\right|,
    \qquad
    D_m := \max_{0\le i\le m} d_i(0).
\end{align*}
By assumption, $D_m\le D$, where $D$ is independent of $n$. At level $0$, each $A_{0,j}$ is an immigration-death process. Therefore,
\begin{align*}
    d_0(t) \le d_0(0)e^{-t} \le D_m e^{-t}.
\end{align*}

Now consider the node $A_{i+1,j}$. Its two input processes $A_{i,2j-1}$ and $A_{i,2j}$ belong to disjoint upstream subtrees. By the assumption that all distinct species have mutually independent initial copy numbers, together with the independence of the reaction-driving Poisson processes, these two subtrees evolve independently. Consequently, the two input processes are independent as stochastic processes. Taking expectations in Proposition~\ref{thm:finite-time-cox-process-representation} and using \eqref{eq:multiplication-lambda-upper-bound}, we obtain
\begin{align}
\label{eq:appendix-d-recursive}
    d_{i+1}(t)
    \le D_m e^{-t}
    + \int_0^t e^{-(t-s)}
        \bigl(d_i(s)^2+2a_i d_i(s)\bigr)\,ds.
\end{align}
Define $E_i(t):=\sum_{\ell=0}^i \frac{t^\ell}{\ell!}$, $M:=1\vee a$. We prove inductively that
\begin{align}\label{eq:appendix-d-upper-bound}
    d_i(t)
    \le e^{-t}M^{2^i-1}G_i E_i(t),
\end{align}
where $G_0:=1+D_m$, $G_{i+1}:=(D_m+2+H_i)G_i^2$, and
\begin{align*}
    H_i:=\int_0^\infty e^{-s}E_i(s)^2\,ds.
\end{align*}
Expanding $E_i(s)^2$ gives
\begin{align*}
    H_i
    =\sum_{j,k=0}^i
        \frac{1}{j!k!}
        \int_0^\infty e^{-s}s^{j+k}\,ds=\sum_{j,k=0}^i \binom{j+k}{j}
      \le \sum_{j,k=0}^i 2^{j+k}
      <4^{i+1}.
\end{align*}
For $i=0$, we have
\begin{align*}
    d_0(t)\le D_m e^{-t}\le G_0e^{-t}.
\end{align*}
Suppose \eqref{eq:appendix-d-upper-bound} holds at level $i$.
Since $a_i\le M^{2^i}$,
\eqref{eq:appendix-d-recursive} implies
\begin{align*}
    d_{i+1}(t)
    \le e^{-t}\Bigg[
        D_m
        +M^{2^{i+1}-2}G_i^2
            \int_0^t e^{-s}E_i(s)^2\,ds
        +2M^{2^{i+1}-1}G_i
            \int_0^t E_i(s)\,ds
    \Bigg].
\end{align*}
Since $M\ge1$, $G_i\ge1$, and
\begin{align*}
    \int_0^t E_i(s)\,ds
    =E_{i+1}(t)-1\le E_{i+1}(t),\quad \int_0^t e^{-s}E_i(s)^2\,ds
    \le H_i,
\end{align*}
we obtain
\begin{align*}
    d_{i+1}(t)
    \le e^{-t}M^{2^{i+1}-1}
        (D_m+2+H_i)G_i^2E_{i+1}(t)=e^{-t}M^{2^{i+1}-1}G_{i+1}E_{i+1}(t).
\end{align*}

This completes the induction. Since $D_m\le D$ and $H_i<4^{i+1}$, there exists a constant $c_D>0$, depending only on $D$, such that $\log G_0\le c_D$ and $\log(D_m+2+H_i)\le c_D(i+1)$. Iterating the recursion gives
\begin{align*}
    \log G_m =2^m\log G_0+ \sum_{i=0}^{m-1}2^{m-1-i}\log(D_m+2+H_i)\le c_D n+ c_D n\sum_{i=0}^{m-1}\frac{i+1}{2^{i+1}} \le 3c_D n.
\end{align*}
Set $\kappa:=\log M+3c_D$, which depends only on $a,D$. By \eqref{eq:appendix-d-upper-bound},
\begin{align*}
    d_m(t) \le e^{\kappa n} \mathbb P(\operatorname{Poisson}(t)\le m).
\end{align*}
Applying the Poisson-tail estimate used to obtain \eqref{eq:t-mean-epsilon-bound}, for $0<\varepsilon\le1$ we obtain
\begin{align*}
    t_{\mathrm{mean}}(\varepsilon)
    \le 2m+2\kappa n+2\log\frac1\varepsilon \le C_{\mathrm{mean}}
    \left(n+\log\frac1\varepsilon\right),
\end{align*}
where $C_{\mathrm{mean}}:=2(1+\kappa)$, since $m\le n$.

\medskip
\noindent
\textbf{Autocovariance convergence. }Set
\begin{align*}
    V_m
    :=\max_{0\le i\le m}\sup_j
        \var(A_{i,j}(0)).
\end{align*}
By the assumptions of the theorem, $D_m\le D$ and $V_m\le V$, where $D,V$ are independent of $n$. For each node, write
\begin{align*}
    C_{i,j}(\tau)
    :=\cov_\pi(A_{i,j}(0),A_{i,j}(\tau)).
\end{align*}
Define
\begin{align*}
    K_0:=4(1+a+D+V)^2,
    \qquad
    L:=64(1+D+V)^2,
    \qquad
    K_{i+1}:=LK_i^2.
\end{align*}
We prove simultaneously that, for every node $A_{i,j}$
and every $t\ge0$,
\begin{align}
    \max\left\{
        \sup_{s\ge0}\E[A_{i,j}(s)^2],
        \E_\pi[A_{i,j}^2]
    \right\}
    &\le K_i,
    \label{eq:appendix-d-second-moment-bound}\\
    d_i(t)
    &\le \sqrt{K_i}\,e^{-t/2},
    \label{eq:appendix-d-auxiliary-mean-bound}\\
    \sup_{\tau\ge0}
    \left|
        \cov(A_{i,j}(t),A_{i,j}(t+\tau))
        -C_{i,j}(\tau)
    \right|
    &\le K_i e^{-t/2}.
    \label{eq:appendix-d-covariance-bound}
\end{align}
For the base case, the immigration-death formulas give
\begin{align*}
    \E[A_{0,j}(t)]
    =a+\bigl(\E[A_{0,j}(0)]-a\bigr)e^{-t},\quad \var(A_{0,j}(t))
    =\E[A_{0,j}(t)]
      +e^{-2t}
        \bigl(\var(A_{0,j}(0))-\E[A_{0,j}(0)]\bigr).
\end{align*}
Hence
\begin{align*}
    \sup_{t\ge0}\E[A_{0,j}(t)^2]
    \le (a+D)^2+a+D+V
    \le K_0,
\end{align*}
and $\E_\pi[A_{0,j}^2]=a^2+a\le K_0$.
Moreover, $C_{0,j}(\tau)=ae^{-\tau}$ and
\begin{align*}
    \cov(A_{0,j}(t),A_{0,j}(t+\tau))
        -C_{0,j}(\tau)=
    e^{-\tau}\Big[
        \bigl(\E[A_{0,j}(0)]-a\bigr)e^{-t}
        +\bigl(\var(A_{0,j}(0))
               -\E[A_{0,j}(0)]\bigr)e^{-2t}
    \Big].
\end{align*}
Therefore,
\begin{align*}
    \sup_{\tau\ge0}
    \left|
        \cov(A_{0,j}(t),A_{0,j}(t+\tau))
        -C_{0,j}(\tau)
    \right|
    \le (a+2D+V)e^{-t}
    \le K_0e^{-t/2}.
\end{align*}
Together with $d_0(t)\le De^{-t}$, these estimates establish the three induction bounds at level $0$. Now, suppose the bounds hold at level $i$.
Fix a node $A_{i+1,j}$ and write
\begin{align*}
    A(t):=A_{i,2j-1}(t),\qquad
    B(t):=A_{i,2j}(t),\qquad
    \lambda(t):=A(t)B(t).
\end{align*}
The input processes $A$ and $B$ are independent.
Consequently,
\begin{align*}
    \sup_{t\ge0}\E[\lambda(t)^2] = \sup_{t\ge0} \E[A(t)^2]\E[B(t)^2] \le K_i^2,\quad \E_\pi[\lambda (0)^2] \le K_i^2.
\end{align*}
Also, $a_{i+1}=a_i^2\le K_i$, and
\begin{align*}
    \left|\E[\lambda(t)]-a_{i+1}\right| =
    \left|\E[A(t)]\E[B(t)]-a_i^2\right|\le
    \left|\E[A(t)]-a_i\right|\E[B(t)]
    +a_i\left|\E[B(t)]-a_i\right| \le 2K_i e^{-t/2}.
\end{align*}
To bound the input covariance error, independence gives
\begin{align*}
    \cov(\lambda(t),\lambda(t+\tau))
    &=
    \cov(A(t),A(t+\tau))
    \cov(B(t),B(t+\tau))\\
    &\quad+
    \cov(A(t),A(t+\tau))
    \E[B(t)]\E[B(t+\tau)]+
    \cov(B(t),B(t+\tau))
    \E[A(t)]\E[A(t+\tau)].
\end{align*}
For either $Y=A$ or $Y=B$, the induction hypothesis implies
\begin{align*}
    |\cov(Y(t),Y(t+\tau))|&\le K_i,\quad |C_Y(\tau)|\le K_i,\quad 
    \left|
        \E[Y(t)]\E[Y(t+\tau)]-a_i^2
    \right|
    \le 2K_i e^{-t/2},
\end{align*}
where $C_A=C_{i,2j-1}$ and $C_B=C_{i,2j}$.
Subtracting the stationary counterpart of the product
covariance identity therefore yields
\begin{align}\label{eq:appendix-d-input-covariance-bound}
    \sup_{\tau\ge0}
    \left|
        \cov(\lambda(t),\lambda(t+\tau))
        -C_\lambda(\tau)
    \right|
    \le (2+3+3)K_i^2e^{-t/2}
    =8K_i^2e^{-t/2},
\end{align}
where
$C_\lambda(\tau):=\cov_\pi(\lambda(0),\lambda(\tau))$. The finite-time mean formula now gives
\begin{align*}
    \left|\E[A_{i+1,j}(t)]-a_{i+1}\right|
    \le
    De^{-t}
    +2K_i\int_0^t e^{-(t-s)}e^{-s/2}\,ds\le (D+4)K_i e^{-t/2}\le \sqrt{K_{i+1}}\,e^{-t/2},
\end{align*}
since $L\ge(D+4)^2$.
This proves \eqref{eq:appendix-d-auxiliary-mean-bound} at level $i+1$. We next verify the second-moment bound. Define
\begin{align*}
    M_t^{(0)}
    :=\int_0^t e^{-(t-s)}\lambda(s)\,ds.
\end{align*}
The input second-moment estimate implies
\begin{align*}
    \sqrt{\E[(M_t^{(0)})^2]}
    \le
    \int_0^t e^{-(t-s)}
        \sqrt{\E[\lambda(s)^2]}\,ds
    \le (1-e^{-t})K_i.
\end{align*}
Moreover,
\begin{align*}
    \E[A_{i+1,j}(0)^2]
    \le V+(a_{i+1}+D)^2\le \bigl(V+(1+D)^2\bigr)K_i^2,
\end{align*}
and the mean formula gives
\begin{align*}
    \E[A_{i+1,j}(t)]
    \le e^{-t}(K_i+D)+(1-e^{-t})K_i
    \le K_i+D.
\end{align*}
By Corollary~\ref{cor:finite-time-moments}
and the Cauchy--Schwarz inequality,
\begin{align*}
    \E[A_{i+1,j}(t)^2]
    \le \E[A_{i+1,j}(t)]+
    \left(
        e^{-t}\sqrt{\E[A_{i+1,j}(0)^2]}
        +\sqrt{\E[(M_t^{(0)})^2]}
    \right)^2 \le
    \bigl(V+(1+D)^2+1+D\bigr)K_i^2\le K_{i+1}.
\end{align*}
The stationary moment formula similarly gives
\begin{align*}
    \E_\pi[A_{i+1,j}^2]
    \le a_{i+1}+\E_\pi[\lambda(0)^2]
    \le K_i+K_i^2
    \le K_{i+1}.
\end{align*}
Thus \eqref{eq:appendix-d-second-moment-bound}
also holds at level $i+1$. Finally, we bound the output covariance error. The input second-moment estimates imply
\begin{align*}
    \var_\pi(\lambda(0))&\le K_i^2,\quad
    \left|
        \cov(A_{i+1,j}(0),\lambda(s))
    \right|
    \le \sqrt{V}\,K_i
    \quad (s\ge0).
\end{align*}
Also, splitting the integration domain according to $u\le v$ and $v<u$ gives, for every $\tau\ge0$,
\begin{align*}
    \int_0^t\int_0^{t+\tau}
        e^{-(t-u)}e^{-(t+\tau-v)}
        e^{-\min\{u,v\}/2}\,dv\,du\le
    2\int_0^t e^{-(t-s)}e^{-s/2}\,ds.
\end{align*}
Subtracting the stationary covariance representation from the finite-time covariance formula, and applying these
estimates together with
\eqref{eq:appendix-d-input-covariance-bound}, yields
\begin{align*}
    &\sup_{\tau\ge0}
    \left|
        \cov(A_{i+1,j}(t),A_{i+1,j}(t+\tau))
        -C_{i+1,j}(\tau)
    \right|\\
    &\le
    (D+4)K_i e^{-t/2}
    +(V+K_i+D)e^{-2t}+
    2\sqrt{V}\,K_i e^{-t}
    +2K_i^2e^{-t}
    +16K_i^2
        \int_0^t e^{-(t-s)}e^{-s/2}\,ds\\
    &\le
    \bigl(2D+V+2\sqrt{V}+39\bigr)K_i^2e^{-t/2}\le
    K_{i+1}e^{-t/2}.
\end{align*}
Here the terms on the right respectively bound the mean
error, the initial variance correction, the two
initial-input covariance terms, the stationary past
contributions, and the input covariance error.
In the last two inequalities, we used
\begin{align*}
    \int_0^t e^{-(t-s)}e^{-s/2}\,ds
    =2(e^{-t/2}-e^{-t})
    \le 2e^{-t/2},
\end{align*}
$K_i\ge1$, and
$L\ge 2D+V+2\sqrt{V}+39$.
This proves \eqref{eq:appendix-d-covariance-bound} at level $i+1$ and completes the induction. Iterating $K_{i+1}=LK_i^2$ gives
\begin{align*}
    K_m
    =L^{2^m-1}K_0^{2^m}
    =L^{n-1}K_0^n.
\end{align*}
Therefore,
\begin{align*}
    \sup_{\tau\ge0}
    \left|
        \cov(A_{m,1}(t),A_{m,1}(t+\tau))
        -C_{m,1}(\tau)
    \right|
    \le L^{n-1}K_0^n e^{-t/2}.
\end{align*}
Since this upper bound decreases with $t$, for $0<\varepsilon\le1$ we obtain
\begin{align*}
    t_{\mathrm{cov}}(\varepsilon) \le 2(n-1)\log L+2n\log K_0+2\log\frac1\varepsilon \le C_{\mathrm{cov}}
    \left(n+\log\frac1\varepsilon\right),
\end{align*}
where $C_{\mathrm{cov}}:=2\max\{1,\log(LK_0)\}$ depends only on $a,D,V$ and is independent of $n$ and $\varepsilon$.

\end{proof}

\section{Proof of Theorem~\ref{thm:multiplication-binary-tree}}
\label{sec:appendix-a-2}

\begin{proof}
As above, $A_i$ denotes a representative species at level $i$ of the binary-tree architecture, with the common stationary process law of the species at that level. Let
\begin{align*}
    a_i:=\E_\pi[A_i]=a^{2^i}, \quad C_i(t):=\cov_\pi(A_i(0),A_i(t)),
\end{align*}
where $C_0(t)=a e^{-t}$. Now define $D_i:=C_i(0) = \var_\pi (A_i )$. Since the two sibling input processes at each multiplication module are independent, Corollary~\ref{lem : cox-process-variance} gives 
\begin{align*}
    D_{i+1} = a_i ^2 + \int_0 ^{\infty} e^{-u} (C_i (u)^2 + 2 a_i ^2 C_i (u) ) du.
\end{align*}
By Lemma~\ref{lem:6.12}, $e^{-t}D_i \leq C_i (t)\leq D_i$. Therefore, we obtain
\begin{align*}
    a_i^2+a_i^2 D_i+\frac13D_i^2 \leq D_{i+1} \leq a_i^2+2a_i^2D_i+D_i^2.
\end{align*}
Indeed,
\begin{align*}
    2a_i^2\int_0^\infty e^{-u}C_i(u)\,du \geq 2a_i^2D_i\int_0^\infty e^{-2u}\,du = a_i^2 D_i,
\end{align*}
and
\begin{align*}
    \int_0^\infty e^{-u}C_i(u)^2\,du \geq D_i^2\int_0^\infty e^{-3u}\,du = \frac13D_i^2.
\end{align*}

\smallskip
\noindent
\textbf{(i) $0 < a < 1$.}
Write $n=2^m$, so that the next level contains $2n$ source inputs. We normalize via $R_m:=\frac{D_m}{a^{n}}$. Dividing the above inequality by $a^{2n}$, we get
\begin{align*}
    1+\frac13R_m^2 \leq R_{m+1} \leq 1+2a^{n}R_m+R_m^2 \leq 1+2R_m+R_m^2.
\end{align*}
To justify the growth of $R_m$, observe that
\[
    1+\frac{r^2}{3}-r =\frac{(r-3/2)^2}{3}+\frac14>0.
\]
Thus, the lower recurrence implies that $R_m$ eventually exceeds $3$. Choose $m_0$ such that $R_{m_0}>3$. For every $m\ge m_0$,
\[
    \frac{R_{m+1}}{3} \ge\left(\frac{R_m}{3}\right)^2.
\]
Iterating gives
\[
    R_m \ge 3 \left(\frac{R_{m_0}}{3} \right)^{2^{m-m_0}}
\]
and hence
\[
    \log R_m \ge 2^{m - m_0} \log \left(\frac{R_{m_0}}{3} \right)
\]
On the other hand, since $R_m\ge1$, $R_{m+1}\le(1+R_m)^2\le4R_m^2$. Iterating this upper bound from $R_0 = 1$ yields $\log R_m \le (2^m - 1) \log 4$. Consequently, there exist constants $c_1 , c_2 > 0$, depending only on $a$, such that, for all sufficiently large $n=2^m$,
\[
c_1 n \le \log R_m \le c_2 n.
\]
Therefore,
\[
a^n e^{c_1 n} \le D_m \le a^n e^{c_2 n}.
\]

\smallskip
\noindent
\textbf{(ii) $a=1$.}
The same inequality becomes
\begin{align*}
    1+D_m+\frac13D_m^2 \leq D_{m+1} \leq 1+2D_m+D_m^2.
\end{align*}
Starting from $D_0=1$, the lower recurrence implies that $D_m$ is increasing and eventually exceeds $3$. Thereafter,
\[
    \frac{D_{m+1}}{3} \ge\left(\frac{D_m}{3}\right)^2.
\]
Together with the upper bound $D_{m+1}\le(1+D_m)^2\le4D_m^2$, same iteration as in part (i) therefore yields constants $c_1 , c_2 > 0$ such that, for all sufficiently large $n=2^m$, $c_1 n \le \log D_m \le c_2 n$. Thus,
\[
    e^{c_1 n} \le D_m \le e^{c_2 n}.
\]

\smallskip
\noindent
\textbf{(iii) $a>1$}. 
We normalize via $R_m:=D_m / a^{2n-1}$. Dividing the inequality by $a^{4n-1}$, we obtain
\begin{align*}
    a^{1-2n}+R_m+\frac{1}{3a}R_m^2 \leq R_{m+1} \leq a^{1-2n}+2R_m+\frac1aR_m^2.
\end{align*}
Since $R_0=1$, the lower recurrence
\[
R_{m+1}\ge R_m+\frac{R_m^2}{3a},
\]
implies that $R_m$ eventually exceeds $3a$. Choose $m_0$ such that $R_{m_0}>3a$. Then, for $m\ge m_0$,
\[
\frac{R_{m+1}}{3a} \ge\left(\frac{R_m}{3a}\right)^2,
\]
and hence,
\[
    \log R_m \ge 2^{m - m_0} \log \left( \frac{R_{m_0}}{3a} \right).
\]
Since $R_m\ge1$ and $a>1$, the upper recurrence gives
\[
    R_{m+1}\le(1+R_m)^2\le4R_m^2.
\]
Hence $\log R_m \le (2^m - 1) \log 4$. Consequently, there exist constants $c_1 , c_2 > 0$, depending only on $a$, such that for all sufficiently large $n=2^m$,
\[
c_1 n \le \log R_m \le c_2n.
\]
It follows that
\[
a^{2n-1} e^{c_1 n} \le D_m \le a^{2n-1} e^{c_2 n}
\]
Since for fixed $t\geq 0$,
\begin{align*}
    e^{-t}D_m\leq C_m(t)\leq D_m,
\end{align*}
the same asymptotic order holds for $C_m(t)$. Recalling that $n=2^m$, we conclude
\begin{alignat*}{2}
a^{2n-1}e^{c_1 n}
&\le C_m(t)\le a^{2n-1}e^{c_2 n},
&\qquad& a>1,\\
e^{c_1 n}
&\le C_m(t)\le e^{c_2 n},
&\qquad& a=1,\\
a^ne^{c_1 n}
&\le C_m(t)\le a^ne^{c_2 n},
&\qquad& 0<a<1.
\end{alignat*}
\end{proof}

\section{Proof of Theorem~\ref{thm:mean-convergence-of-linear-chain-multiplication}}
\label{sec:appendix-d}

\begin{proof}
\smallskip
\noindent
\textbf{Mean convergence.} Recall that $S_1=A_1$. Define
\begin{align*}
    e_i(t):=\left|\E[S_i(t)]-a^i\right|,
    \qquad
    c_i:=\left|\E[A_i(0)]-a\right|,
    \qquad
    d_i:=e_i(0).
\end{align*}

For each $1\le i<n$, the process $S_i$ depends only on the initial populations and reaction-driving randomness associated with the first $i$ source species and their intermediate multiplication modules. The source process $A_{i+1}$ is generated independently of this upstream network. Therefore, by the mutual independence of the initial copy numbers and the independent reaction-driving Poisson processes, $S_i$ and $A_{i+1}$ are independent as stochastic processes. Taking expectations in Proposition~\ref{thm:finite-time-cox-process-representation} and using \eqref{eq:multiplication-lambda-upper-bound}, we obtain
\begin{align}\label{eq:appendix-f-recursive}
    e_{i+1}(t)
    \le e^{-t}d_{i+1}
    +\int_0^t e^{-(t-s)}
        \left[
            (a+c_{i+1}e^{-s})e_i(s)
            +a^i c_{i+1}e^{-s}
        \right]\,ds.
\end{align}
We construct recursively a polynomial $Q_i$ with non-negative coefficients such that $e_i(t)\le e^{-t}Q_i(t)$. Since $S_1=A_1$ is an immigration-death process, we may take $Q_1(t)=d_1$. Suppose
\begin{align*}
    Q_i(t)
    =\sum_{\ell=0}^{i-1}b_{i,\ell}\frac{t^\ell}{\ell!},
    \qquad b_{i,\ell}\ge0,
\end{align*}
and define $R_i:=\sum_{\ell=0}^{i-1}b_{i,\ell}$. Then
\begin{align*}
    \int_0^t e^{-s}Q_i(s)\,ds
    \le
    \sum_{\ell=0}^{i-1}b_{i,\ell}
        \int_0^\infty e^{-s}\frac{s^\ell}{\ell!}\,ds=R_i.
\end{align*}
Substituting the bound for $e_i$ into
\eqref{eq:appendix-f-recursive} gives
\begin{align*}
    e_{i+1}(t)
    \le e^{-t}\left[
        d_{i+1}
        +c_{i+1}R_i
        +a\int_0^t Q_i(s)\,ds
        +a^i c_{i+1}t
    \right].
\end{align*}
We may therefore define
\begin{align*}
    Q_{i+1}(t)
    :=d_{i+1}
      +c_{i+1}R_i
      +a\int_0^t Q_i(s)\,ds
      +a^i c_{i+1}t.
\end{align*}
Integration maps $t^\ell/\ell!$ to
$t^{\ell+1}/(\ell+1)!$, so
\begin{align*}
    R_{i+1}
    =(a+c_{i+1})R_i+d_{i+1}+c_{i+1}a^i.
\end{align*}
Starting from $R_1=d_1$, iteration gives
\begin{align}\label{eq:appendix-f-coeff-def}
    R_i
    =
    d_1\prod_{q=2}^i(a+c_q)
    +\sum_{p=2}^i
        (d_p+c_p a^{p-1})
        \prod_{q=p+1}^i(a+c_q),
\end{align}
where an empty product equals one. By the assumptions of the theorem, $c_i\le C$
and $d_i\le D$, where $C,D$ are independent of $n$. Set $K:=\max\{1,a+C\}$. Using $a\le K$ in \eqref{eq:appendix-f-coeff-def}, we obtain
\begin{align*}
    R_i
    \le
    DK^{i-1}
    +\sum_{p=2}^i
        (D+CK^{p-1})K^{i-p}\le
    (D+C)iK^{i-1}\le
    (1+D+C)iK^i. 
\end{align*}

Set $B:=2(1+D+C)K$. Since $i\le 2^i$, the preceding bound gives $R_i\le B^i$. Taking $i=n$, we obtain
\[
e_n(t)\le e^{-t}Q_n(t) \le B^n e^{-t}\sum_{\ell=0}^{n-1}\frac{t^\ell}{\ell!}
=B^n\mathbb P(\operatorname{Poisson}(t)\le n-1).
\]
Applying the Poisson-tail estimate used in Theorem~\ref{thm:6.1}, for $0<\varepsilon\le1$ we obtain
\begin{align*}
t_{\mathrm{mean}}(\varepsilon) \le 2(n-1)+2n\log B+2\log\frac1\varepsilon \le C_{\mathrm{mean}} \left(n+\log\frac1\varepsilon\right),
\end{align*}
where $C_{\mathrm{mean}}:=2(1+\log B)$ depends only on $a,C,D$ and is independent of $n$ and $\varepsilon$.

\smallskip
\noindent
\textbf{Autocovariance convergence.} By the additional uniform variance assumption, choose $V\ge0$, independent of $n$, such that
\begin{align*}
    \var(A_i(0))\le V,
    \qquad
    \var(S_i(0))\le V
    \qquad (1\le i\le n).
\end{align*}
Write $C_{S_i}(\tau):=\cov_\pi(S_i(0),S_i(\tau))$. Define $K_0:=4(1+a+C+V)^2$, $L:=64(1+D+V)^2$, and $K_1:=K_0$, $K_{i+1}:=LK_iK_0$. Each source $A_i$ is an immigration-death process with stationary mean $a$ and stationary autocovariance
$ae^{-\tau}$. Its finite-time moment formulas give
\begin{align*}
    \max\left\{
        \sup_{t\ge0}\E[A_i(t)^2],
        \E_\pi[A_i^2]
    \right\}
    &\le K_0,\\
    \left|\E[A_i(t)]-a\right|
    &\le Ce^{-t}
    \le \sqrt{K_0}\,e^{-t/2},\\
    \sup_{\tau\ge0}
    \left|
        \cov(A_i(t),A_i(t+\tau))-ae^{-\tau}
    \right|
    &\le (a+2C+V)e^{-t}
    \le K_0e^{-t/2}.
\end{align*}

To control the constants in the covariance recursion,
we prove simultaneously that, for every $i$ and $t\ge0$,
\begin{align}
    \max\left\{
        \sup_{s\ge0}\E[S_i(s)^2],
        \E_\pi[S_i^2]
    \right\}
    &\le K_i,
    \label{eq:appendix-f-second-moment-bound}\\
    e_i(t)
    &\le \sqrt{K_i}\,e^{-t/2},
    \label{eq:appendix-f-auxiliary-mean-bound}\\
    \sup_{\tau\ge0}
    \left|
        \cov(S_i(t),S_i(t+\tau))-C_{S_i}(\tau)
    \right|
    &\le K_i e^{-t/2}.
    \label{eq:appendix-f-covariance-bound}
\end{align}
These bounds hold for $i=1$ because $S_1=A_1$ and $K_1=K_0$. Suppose they hold for some $i<n$. Define
\begin{align*}
    \lambda(t):=S_i(t)A_{i+1}(t),
    \qquad
    C_\lambda(\tau)
    :=\cov_\pi(\lambda(0),\lambda(\tau)).
\end{align*}
Independence of the two input processes implies
\begin{align*}
    \sup_{t\ge0}\E[\lambda(t)^2]
    \le K_iK_0,\quad\E_\pi[\lambda(0)^2]
    \le K_iK_0.
\end{align*}
In particular, $a^{i+1} =\E_\pi[\lambda(0)] \le \sqrt{K_iK_0}$. Also,
\begin{align*}
    \left|\E[\lambda(t)]-a^{i+1}\right|
    =
    \left|
        \E[S_i(t)]\E[A_{i+1}(t)]-a^i a
    \right|\le
    e_i(t)\E[A_{i+1}(t)]
    +a^i\left|\E[A_{i+1}(t)]-a\right|\le 2\sqrt{K_iK_0}\,e^{-t/2}.
\end{align*}

For the input covariance, the product covariance identity gives
\begin{align*}
    \cov(\lambda(t),\lambda(t+\tau))&=
    \cov(S_i(t),S_i(t+\tau))
    \cov(A_{i+1}(t),A_{i+1}(t+\tau))\\
    &\quad+
    \cov(S_i(t),S_i(t+\tau))
    \E[A_{i+1}(t)]\E[A_{i+1}(t+\tau)]\\
    &\quad+
    \cov(A_{i+1}(t),A_{i+1}(t+\tau))
    \E[S_i(t)]\E[S_i(t+\tau)].
\end{align*}
The moment and mean bounds imply
\begin{align*}
    \left|
        \E[S_i(t)]\E[S_i(t+\tau)]-a^{2i}
    \right|
    &\le 2K_i e^{-t/2},\\
    \left|
        \E[A_{i+1}(t)]\E[A_{i+1}(t+\tau)]-a^2
    \right|
    &\le 2K_0e^{-t/2}.
\end{align*}
Moreover, the finite-time and stationary covariances of
$S_i$ are bounded in absolute value by $K_i$, and those
of $A_{i+1}$ are bounded by $K_0$.
Subtracting the stationary counterpart of the product
covariance identity therefore yields
\begin{align}\label{eq:appendix-f-input-covariance-bound}
    \sup_{\tau\ge0}
    \left|
        \cov(\lambda(t),\lambda(t+\tau))
        -C_\lambda(\tau)
    \right|
    \le (2+3+3)K_iK_0e^{-t/2}
    =8K_iK_0e^{-t/2}.
\end{align}
The finite-time mean formula now gives
\begin{align*}
    e_{i+1}(t)
    &\le
    De^{-t}
    +2\sqrt{K_iK_0}
        \int_0^t e^{-(t-s)}e^{-s/2}\,ds\le
    (D+4)\sqrt{K_iK_0}\,e^{-t/2}\le
    \sqrt{K_{i+1}}\,e^{-t/2}.
\end{align*}
Here we used
\begin{align*}
    \int_0^t e^{-(t-s)}e^{-s/2}\,ds
    =2(e^{-t/2}-e^{-t})
    \le 2e^{-t/2},
\end{align*}
$K_iK_0\ge1$, and $L\ge(D+4)^2$. We next verify the second-moment bound. Set
\begin{align*}
    M_t^{(0)}
    :=\int_0^t e^{-(t-s)}\lambda(s)\,ds.
\end{align*}
Then
\begin{align*}
    \sqrt{\E[(M_t^{(0)})^2]}
    \le
    \int_0^t e^{-(t-s)}
        \sqrt{\E[\lambda(s)^2]}\,ds\le
    (1-e^{-t})\sqrt{K_iK_0}.
\end{align*}
Also,
\begin{align*}
    \E[S_{i+1}(0)^2]
    \le
    V+(a^{i+1}+D)^2\le
    \bigl(V+(1+D)^2\bigr)K_iK_0,
\end{align*}
and the mean formula implies
\begin{align*}
    \E[S_{i+1}(t)]
    \le \sqrt{K_iK_0}+D.
\end{align*}
By Corollary~\ref{cor:finite-time-moments}
and Cauchy--Schwarz,
\begin{align*}
    \E[S_{i+1}(t)^2]
    \le
    \E[S_{i+1}(t)]+
    \left(
        e^{-t}\sqrt{\E[S_{i+1}(0)^2]}
        +\sqrt{\E[(M_t^{(0)})^2]}
    \right)^2\le
    \bigl(V+(1+D)^2+1+D\bigr)K_iK_0\le K_{i+1}.
\end{align*}
The stationary moment formula similarly gives
\begin{align*}
    \E_\pi[S_{i+1}^2]
    \le a^{i+1}+\E_\pi[\lambda(0)^2]\le \sqrt{K_iK_0}+K_iK_0\le 2K_iK_0
    \le K_{i+1}.
\end{align*}
Thus \eqref{eq:appendix-f-second-moment-bound}
holds at level $i+1$. Finally, the input moment bounds give
\begin{align*}
    \var_\pi(\lambda (0))\le K_iK_0,\quad \left|\cov(S_{i+1}(0),\lambda(s))\right|
    \le \sqrt{VK_iK_0}
    \quad (s\ge0).
\end{align*}
For every $\tau\ge0$, splitting the integration domain
according to $u\le v$ and $v<u$ gives
\begin{align*}
    \int_0^t\int_0^{t+\tau}
        e^{-(t-u)}e^{-(t+\tau-v)}
        e^{-\min\{u,v\}/2}\,dv\,du\le
    2\int_0^t e^{-(t-s)}e^{-s/2}\,ds
    \le 4e^{-t/2}.
\end{align*}
Consequently, by
\eqref{eq:appendix-f-input-covariance-bound},
the contribution of the input covariance error to the
output covariance error is at most
$32K_iK_0e^{-t/2}$.

Subtracting the stationary covariance representation
from the finite-time covariance formula and bounding
the remaining terms yields
\begin{align*}
    &\sup_{\tau\ge0}
    \left|
        \cov(S_{i+1}(t),S_{i+1}(t+\tau))
        -C_{S_{i+1}}(\tau)
    \right|\\
    &\le
    (D+4)\sqrt{K_iK_0}\,e^{-t/2}
    +\bigl(V+\sqrt{K_iK_0}+D\bigr)e^{-2t}+
    2\sqrt{VK_iK_0}\,e^{-t}
    +2K_iK_0e^{-t}
    +32K_iK_0e^{-t/2}\\
    &\le
    \bigl(2D+V+2\sqrt{V}+39\bigr)K_iK_0e^{-t/2}\le
    K_{i+1}e^{-t/2}.
\end{align*}
The five terms bound, respectively, the mean error,
the initial variance correction, the two initial-input
covariance terms, the stationary past contributions,
and the input covariance error.
The final inequality follows from
$L\ge 2D+V+2\sqrt{V}+39$.
This completes the induction. Since $K_1=K_0$ and $K_{i+1}=LK_iK_0$, we obtain $K_n=L^{n-1}K_0^n$. Therefore,
\begin{align*}
    \sup_{\tau\ge0}
    \left|
        \cov(S_n(t),S_n(t+\tau))-C_{S_n}(\tau)
    \right|
    \le L^{n-1}K_0^n e^{-t/2}.
\end{align*}
Since this upper bound decreases with $t$, for $0<\varepsilon\le1$ we obtain
\[
t_{\mathrm{cov}}(\varepsilon) \le 2(n-1)\log L+2n\log K_0
   +2\log\frac1\varepsilon \le C_{\mathrm{cov}} \left(n+\log\frac1\varepsilon\right),
\]
where $C_{\mathrm{cov}}:=2\max\{1,\log(LK_0)\}$ depends only on $a,C,D,V$ and is independent of $n$ and $\varepsilon$.
\end{proof}

\section{Proof of Theorem~\ref{thm:fixed-time-aucotovariance-bounds-for-linear-chain-multiplication}}
\label{sec:appendix-e}

\begin{proof}
Define $C_i (t) := \cov_\pi (S_i (0), S_i (t))$, $D_i:=C_i(0)$. The base case is $S_1=A_1$, $C_1(t)=a e^{-t}$, $D_1=a$. By Corollary~\ref{lem : cox-process-variance},
\begin{align*}
    D_{i+1} = a^{i+1} + \int_0 ^{\infty} e^{-u} \{ (a^2 + a e^{-u} )C_i (u) + a^{2i+1} e^{-u} \} du.
\end{align*}
By Lemma~\ref{lem:6.12}, $e^{-t}D_i\leq C_i(t)\leq D_i$. We obtain
\begin{align*}
    a^{i+1} + \left(\frac{a^2}{2}+\frac{a}{3}\right)D_i + \frac12 a^{2i+1} \leq D_{i+1} \le a^{i+1} + \left(a^2+\frac{a}{2}\right)D_i + \frac12 a^{2i+1}.
\end{align*}
We now analyze this recursion.

\smallskip
\noindent
\textbf{(i) $0<a<1$.}
We normalize via $ R_n:=\frac{D_n}{a^n}$. Dividing the above inequalities by $a^{n+1}$, we obtain
\begin{align*}
    1+ \left(\frac {a}2+\frac13\right)R_n + \frac12(a)^n \leq R_{n+1} \leq 1+ \left(a +\frac12\right)R_n + \frac12 a^n.
\end{align*}
Since $a<1$, the forcing terms in the upper recurrence are uniformly bounded. The resulting upper bound therefore depends on whether the coefficient $a+1/2$ is less than, equal to, or greater than one. Hence:
\begin{alignat*}{2}
c_1 &\le R_n\le c_2,
&\qquad& 0<a<\frac12,\\
c_1 &\le R_n \le c_2n,
&\qquad& a=\frac12,\\
c_1 &\le R_n\le c_2\left(a+\frac12\right)^n,
&\qquad& \frac12<a<1.
\end{alignat*}
Therefore,
\begin{alignat*}{2}
c_1a^n &\le D_n\le c_2a^n,
&\qquad& 0<a<\frac12,\\
c_1 2^{-n} &\le D_n\le c_2 n2^{-n},
&\qquad& a=\frac12,\\
c_1a^n &\le D_n\le c_2a^n\left(a+\frac12\right)^n,
&\qquad& \frac12<a<1.
\end{alignat*}

\smallskip
\noindent
\textbf{(ii) $a=1$.}
The recursion becomes
\begin{align*}
    \frac32+\frac56D_n \leq D_{n+1} \leq \frac32+\frac32D_n.
\end{align*}
Thus, $c_1 \le D_n \le c_2(3/2)^n$.

\smallskip
\noindent
\textbf{(iii) $a>1$.}
We normalize via $R_n:=\frac{D_n}{a^{2n-1}}$. Dividing the recursion by $a^{2n+1}$, we get
\begin{align*}
    a^{-n} + \left(\frac12+\frac{1}{3a}\right)R_n + \frac12 \leq R_{n+1} \leq a^{-n} + \left(1+\frac{1}{2a}\right)R_n + \frac12.
\end{align*}
Hence,
\begin{align*}
    c_1 \le R_n \le c_2\left(1+\frac{1}{2a}\right)^n.
\end{align*}
Therefore,
\begin{align*}
    c_1a^{2n-1} \le D_n
    \le c_2a^{2n-1}\left(1+\frac{1}{2a}\right)^n.
\end{align*}
Since for each fixed $t\geq 0$,
\begin{align*}
    e^{-t}D_n\leq C_n(t)\leq D_n,
\end{align*}
the same upper and lower asymptotic bounds hold for $C_n(t)$. This proves the claim.
\end{proof}

\end{document}